\documentstyle[12pt, amstex, amssymb,eucal, xypic, latexsym]{amsart}

\def\gr{{\mathrm gr}}
\def\k{{\Bbbk}}

\def\g{{\frak g}}
\def\A{{\bf A}_{e}^{\text{op}}}
\def\h{{\frak h}}
\def\sl{{\frak sl}}

\def\Q{{\widehat{Q}_\chi}}

\def\z{{\frak z}}
\def\m{{\frak m}}
\def\n{{\frak n}}
\def\p{{\frak p}}
\def\Z{{\Bbb Z}}

\def\End{\mathop{\fam0 End}}
\def\codim{\mathop{\fam0 codim}\nolimits}
\def\Lie{\mathop{\fam0 Lie}}
\def\Ad{{\mathrm Ad}\,}
\def\Der{\mathop{\fam0 Der}}
\def\Mat{\mathop{\fam0 Mat}\nolimits}
\def\sl{\mathop{\frak sl}\nolimits}

\def\ad{{\mathrm ad\,}}

\newcommand{\into}{\,\hookrightarrow\,}
\newcommand{\map}{\longrightarrow}

\def\la{\langle}
\def\ra{\rangle}

\numberwithin{equation}{subsection}

\theoremstyle{plain}
\newtheorem{theorem}{Theorem}[section]
\newtheorem{corollary}{Corollary}[section]
\newtheorem{prop}{Proposition}[section]
\newtheorem{lemma}{Lemma}[section]

\theoremstyle{definition}

\newtheorem{conj}{Conjecture}[section]
\newtheorem{question}{Question}[section]

\theoremstyle{remark}

\newtheorem{rem}{Remark}[section]

\def\subtitle#1. {{\medskip\bf#1\par\nobreak\smallskip}}
\def\proclaim#1. {\medbreak\bgroup\noindent\bf#1. \it}

\def\endproclaim{\egroup
\ifdim\lastskip<\medskipamount\removelastskip\medskip\fi}
\newcount
=0
\def\citedef#1 {\advance\citation by1
  \expandafter\edef\csname#1\endcsname{{\the\citation}}
  \checkendcitedef}
\def\checkendcitedef#1{\ifx#1\endcitedef\else\citedef#1\fi}
\def\cite#1{\csname#1\endcsname}
\citedef Bac BoBr BK Bou BrK Ca C DvvD Dix GG GS G  Ja J J' J-pt
J4 J2 J3 J-wp KT KW Kh M McR Mc McS MS Mo P02  P03 RS Sh S   V
\endcitedef
\newtoks\nextauth
\newif\iffirstauth
\def\checkendauth#1{\ifx\endauth#1
        \iffirstauth\the\nextauth
        \else{} and \the\nextauth\fi,
    \else\iffirstauth\the\nextauth\firstauthfalse
        \else, \the\nextauth\fi
        \expandafter\auth\expandafter#1\fi}
\def\auth#1 #2 {\nextauth={#1 #2}\checkendauth}
\newif\ifinbook
\newif\ifbookref
\def\nextref#1 {\bookreffalse\inbookfalse
    \bibitem[\cite{#1}]{}
    \firstauthtrue
    \ignorespaces}
\def\paper#1{{\it#1,}}
\def\In#1{\inbooktrue In #1,}
\def\book#1{\bookreftrue{\it#1,}}
\def\journal#1{#1\ifinbook,\fi}
\def\bookseries#1{#1,}
\def\Vol#1{\ifbookref Vol. #1,\else\ifinbook Vol. #1,\else{\bf#1}\fi\fi
    \space\ignorespaces}

\def\publisher#1{#1,}
\def\Year#1{\ifbookref #1.\else\ifinbook #1,\else(#1)\fi\fi
    \space\ignorespaces}
\def\Pages#1{\ifinbook pp. #1.\else #1.\fi}

\begin{document}

\title{Enveloping algebras of Slodowy slices\\
and the Joseph ideal}

\author{Alexander Premet}
\thanks{\nonumber{\it Mathematics Subject Classification} (2000 {\it revision}).
Primary 17B35. Secondary 17B63, 17B81.}
\address{Department of Mathematics, The University of Manchester, Oxford Road,
M13 9PL, UK} \email{sashap@@ma.man.ac.uk}

\begin{abstract}
\noindent Let $G$ be a simple algebraic group over an
algebraically closed field $\k$ of characteristic $0$, and
$\g=\text{Lie}\,G$. Let $(e,h,f)$ be an $\sl_2$-triple in $\g$
with $e$ being a long root vector in $\g$. Let
$(\,\cdot\,,\,\cdot\,)$ be the $G$-invariant bilinear form on $\g$
with $(e,f)=1$ and let $\chi\in\g^*$ be such that $\chi(x)=(e,x)$
for all $x\in\g$. Let ${\mathcal S}$ be the Slodowy slice at $e$
through the adjoint orbit of $e$ and let $H$ be the enveloping
algebra of ${\mathcal S}$; see [\cite{P02}]. In this note we give
an explicit presentation of $H$ by generators and relations.  As a
consequence we deduce that $H$ contains an ideal of codimension
$1$ which is unique if $\g$ is not of type $\mathrm A$. Applying
Skryabin's equivalence of categories we then construct an explicit
Whittaker model for the Joseph ideal of $U(\g)$. Inspired by
Joseph's Preparation Theorem we prove that there exists a
homeomorphism between the primitive spectrum of $H$ and the
spectrum of all primitive ideals of infinite codimension in
$U(\g)$ which respects Goldie rank and Gelfand--Kirillov
dimension. We study highest weight modules for the algebra $H$ and
apply earlier results of Mili{\v c}i{\'c}--Soergel and Backelin to
express the composition multiplicities of the Verma modules for
$H$ in terms of some inverse parabolic Kazhdan--Lusztig
polynomials. Our results confirm in the minimal nilpotent case the
de Vos--van Driel conjecture on composition multiplicities of
Verma modules for finite ${\mathcal W}$-algebras. We also obtain
some general results on the enveloping algebras of Slodowy slices
and determine the associated varieties of related primitive ideals
of $U(\g)$. A sequel to this paper will treat modular aspects of
this theory.
\end{abstract}
\maketitle

\bigskip

\section{\bf Introduction}
\subsection{} Let $\k$ be an algebraically closed field of
characteristic $0$ and let $G$ be a simple algebraic group over
$\k$. Let $\g=\Lie\,G$ and let $(e,h,f)$ be an $\sl_2$-triple in
$\g$. Let $(\,\cdot\,,\,\cdot\,)$ be the $G$-invariant bilinear
form on $\g$ with $(e,f)=1$ and define $\chi=\chi_e\in\g^*$ by
setting $\chi(x)=(e,x)$ for all $x\in\g$. Let ${\cal O}_\chi$
denote the coadjoint orbit of $\chi$.

Let ${\mathcal S}_e=e+\text{Ker}\,\ad\,f$ be the Slodowy slice at
$e$ through the adjoint orbit of $e$ and let $H_\chi$ be the
enveloping algebra of ${\mathcal S}_e$; see [\cite{P02},
\cite{GG}, \cite{BrK}]. Recall that
$H_\chi={\End}_{\g}\,(Q_\chi)^{\text{op}}$ where $Q_\chi$ is a
generalised Gelfand--Graev module for $U(\g)$ associated with the
$\sl_2$-triple $(e,h,f)$. The module $Q_\chi$ is induced from a
one-dimensional module $\k_\chi$ over a nilpotent subalgebra
$\m_\chi$ of $\g$ such that $\dim \m_\chi=\frac{1}{2}\dim {\cal
O}_\chi$. The subalgebra $\m_\chi$ is $(\ad h)$-stable, all
weights of $\ad h$ on $\m_\chi$ are negative, and $\chi$ vanishes
on the derived subalgebra of $\m_\chi$. The action of $\m_\chi$ on
$\k_\chi=\k 1_\chi$ is given by $x(1_\chi)=\chi(x)1_\chi$ for all
$x\in\m_\chi$; see (2.1) for more detail.

Let $\z_\chi$ denote the stabiliser of $\chi$ in $\g$. Clearly,
$\z_\chi$ coincides with the centraliser ${\mathfrak c}_\g(e)$ of
$e$ in $\g$. The subalgebra $\z_\chi$ is $(\ad h)$-stable and it
follows from the $\sl_2$-theory that all weights of $\ad h$ on
$\z_\chi$ are nonnegative integers. Let $x_1,\ldots,x_r$ be a
basis of $\z_\chi$ such that $[h,x_i]=n_i x_i$ for some
$n_i\in\Z_+$. By [\cite{P02}, Theorem~4.6], to each basis vector
$x_i$ one can attach an element $\Theta_{x_i}\in H_\chi$ in such a
way that the monomials $\Theta_{x_1}^{i_1}\Theta_{x_2}^{i_2}\cdots
\Theta_{x_r}^{i_r}$ with $(i_1,i_2,\ldots,i_r)\in \Z_+^r$ form a
basis of $H_\chi$ over $\k$. We say that the monomial
$\Theta_{x_1}^{a_1}\Theta_{x_2}^{a_2}\cdots \Theta_{x_r}^{a_r}$
has {\it Kazhdan degree} $\sum_{i=1}^r a_i(n_i+2)$ and denote by
$H_\chi^k$ the span of all monomials as above of Kazhdan degree
$\le k$. According to [\cite{P02}, (4.6)], we then have
$H_\chi=\bigcup_{k\ge 0}\,H_\chi^k$ and $H_\chi^i\cdot
H_\chi^j\subseteq H_\chi^{i+j}$ for all $i,j\in\Z_+$. In other
words, $\{H_\chi^k\,|\,k\in\Z_+\}$ is an increasing filtration of
the algebra $H_\chi$. We call it the {\it Kazhdan filtration} of
$H_\chi$. The corresponding graded algebra $\text{gr}\,H_\chi$ is
a polynomial algebra in
$\text{gr}\,\Theta_{x_1},\text{gr}\,\Theta_{x_2},\ldots,\text{gr}\,
\Theta_{x_r}$ which identifies naturally with the function algebra
$\k[{\mathcal S}_e]$ on the special transverse slice ${\mathcal
S}_e=e+\text{Ker}\,\ad f$ endowed with its Slodowy grading.

According to [\cite{P02}, Theorem~4.6(iv)],
$$
[\Theta_{x_i},\Theta_{x_j}]\,=\,\Theta_{x_j}\circ\Theta_{x_i}-
\Theta_{x_i}\circ\Theta_{x_j}\,\equiv\,
\Theta_{[x_i,x_j]}+q_{ij}(\Theta_1,\ldots,\Theta_r)\ \ \,
\big({\mathrm mod}\ \,  H_{\chi}^{n_{i}+n_{j}}\big)
$$
where $q_{ij}$ is a polynomial in $r$ variables  with initial form
of total degree $\ge 2$. Using this result we prove in (2.3) that
there exists an associative $\k[t]$-algebra ${\mathcal H}_\chi$
free as a module over $\k[t]$ and such that
\[{\mathcal H}_{\chi}/(t-\lambda){\mathcal H}_{\chi}\cong
\left\{\begin{array}{ll} H_\chi &  \mbox{ if $\, \lambda\ne 0$},\\
U(\z_\chi) & \mbox{ if $\, \lambda=0$}
\end{array}
\right. \] as $\k$-algebras. Thus $H_\chi$ is a deformation of the
universal enveloping algebra $U(\z_\chi)$.

We have a certain degree of freedom in our choice of PBW
generators $\Theta_{x_i}$. In (2.2) we show that they can be
chosen such that the map $\Theta_{x_i}\mapsto
(-1)^{n_i}\Theta_{x_i}$, $\,1\le i\le r,$ extends to an
automorphism of the algebra $H_\chi$. This automorphism, denoted
by $\sigma$, will play an important r{\^o}le later on.

\subsection{}
Given $x\in\g$ we denote by $Z_G(x)$ the centraliser of $x$ in
$G$. It is well-known that $C(e):=Z_G(h)\cap Z_G(e)$ is a Levi
subgroup of $Z_G(e)$ and the centraliser $Z_G(e)$ decomposes as a
semidirect product of $C(e)$ and the unipotent radical
$R_u(Z_G(e))$. For $i\in\Z_+$ put
$\z_\chi(i)=\{x\in\z_\chi\,|\,\,[h,x]=ix\}$. It is well-known that
$\z_\chi(0)=\Lie\,C(e)$. Clearly, the group $C(e)$ preserves each
subspace $\z_\chi(i)$ of $\z_\chi$.

In [\cite{GG}], Gan and Ginzburg have found a different
realisation of the algebra $H_\chi$ which enables one to observe
that the reductive group $C(e)$ acts on $H_\chi$ as algebra
automorphisms. Moreover, this action of $C(e)$ preserves the
Kazhdan filtration of $H_\chi$; see (2.1) for more detail. In
Section~2 we show that there exists an injective
$C(e)$-equivariant linear map $\Theta\colon\,\z_\chi\rightarrow
H_\chi,\ x\mapsto\Theta_x,$ whose image generates $H_\chi$ as an
algebra, such that $\text{gr}\,\Theta(\z_\chi)\cong\,\z_\chi$ as
graded $C(e)$-modules. The subspace $\Theta(\z_\chi(0))$ can be
chosen to be a Lie subalgebra of $H_\chi$ with respect to the
commutator product in $H_\chi$. It follows that the map $\Theta$
can be selected in such a way that
$$[\Theta_x,\Theta_y]=\Theta_{[x,y]}\qquad\ \big(\forall\,x\in\z_\chi(0),\ \,\forall\,
y\in\z_\chi\big).$$ This is a consequence of Lemma~2.4 which
states, in particular, that the Lie algebra homomorphism ${\mathrm
ad}\circ\Theta\colon \z_\chi(0)\map{\mathrm Der}\,(H_\chi)$
coincides with the differential of the locally finite (rational)
action of $C(e)$ on $H_\chi$. Combined with Lemma~2.1 and Weyl's
theorem on complete reducibility this implies that every two-sided
ideal of $H_\chi$ is $\sigma$-stable; see Corollary~2.1.
\subsection{} For $\chi=(e,\,\cdot\,)$ we let ${\mathcal C}_\chi$
denote the category of all $\g$-modules on which $x-\chi(x)$ acts
locally nilpotently for all $x\in\m_\chi$. Given a $\g$-module $M$
we set $$\text{Wh}(M):=\,\{m\in M\,|\,x.m=\chi(x)m\
\,(\forall\,x\in\m_\chi)\}.$$ It should be mentioned here that the
algebra $H_\chi$ acts on $\text{Wh}(M)$ via a canonical
isomorphism $H_\chi\cong \big(U(\g)/N_\chi\big)^{\ad \m_\chi}$
where $N_\chi$ denotes the left ideal of $U(\g)$ generated by all
$x-\chi(x)$ with $x\in\m_\chi$. In the Appendix to [\cite{P02}],
Skryabin proved that the functors $V\rightsquigarrow
Q_\chi\otimes_{H_\chi}V$ and $M\rightsquigarrow \mbox{Wh}(M)$ are
mutually inverse equivalences between the category of all
$H_\chi$-modules and the category ${\mathcal C}_\chi$; see also
[\cite{GG}, Theorem~6.1].

Skryabin's equivalence implies that for any irreducible
$H_\chi$-module $V$ the annihilator
$\text{Ann}_{U(\g)}(Q_\chi\otimes_{H_\chi} V)$ is a primitive
ideal of $U(\g)$. By the Irreducibility Theorem, the associated
variety ${\mathcal VA}({\cal I})$ of any primitive ideal ${\cal
I}$ of $U(\g)$ is the closure of a nilpotent orbit in $\g^*$.
Generalising a classical result of Kostant on Whittaker modules we
show in Section~3 that for any irreducible $H_\chi$-module $V$ the
associated variety of $\text{Ann}_{U(\g)}(Q_\chi\otimes_{H_\chi}
V)$ contains the coadjoint orbit ${\cal O}_\chi$. In the most
interesting case where $V$ is a finite dimensional irreducible
$H_\chi$-module we prove that
$${\mathcal VA}\big(\text{Ann}_{U(\g)}(Q_\chi\otimes_{H_\chi} V)\big)
=\overline{\cal O}_\chi\quad\,\text{and}\,\quad
\text{Dim}(Q_\chi\otimes_{H_\chi} V)=\frac{1}{2}\dim
\overline{\cal O}_\chi$$ where $\text{Dim}(M)$ is the
Gelfand--Kirillov dimension of a finitely generated $U(\g)$-module
$M$. In  particular, this implies that for any irreducible finite
dimensional $H_\chi$-module $V$ the irreducible $U(\g)$-module
$Q_\chi\otimes_{H_\chi}V$ is holonomic.
\subsection{}
 Let
$\h$ be a Cartan subalgebra of $\g$, and let $\Phi$ be the root
system of $\g$ relative to $\h$.  Let
$\Pi=\{\alpha_1,\ldots,\alpha_\ell\}$ be a basis of simple roots
in $\Phi$ with the elements in $\Pi$ numbered as in [\cite{Bou}],
and let $\Phi^+$  be the positive system of $\Phi$ relative to
$\Pi$.  If $\g$ is not of type $\rm A$ or $\rm C$, there is a
unique long root in $\Pi$ linked with the lowest root
$-\widetilde{\alpha}$ on the extended Dynkin diagram of $\g$; we
call it $\beta$. For $\g$ of type $A_n$ and $C_n$ we set
$\beta=\alpha_n$. Choose root vectors $e_{\beta},e_{-\beta} \in
\g$ corresponding to roots $\beta$ and $-\beta$ such that
$(e_{\beta},[e_\beta,e_{-\beta}],e_{-\beta})$ is an $\sl_2$-triple
and put $h_\beta=[e_\beta,e_{-\beta}]$.

In this note we investigate the algebra $H_\chi$ in the case where
$(e,h,f)=(e_\beta,h_\beta,e_{-\beta})$. Then ${\cal
O}_{\chi}={\cal O}_{\rm min}$, the minimal nonzero nilpotent orbit
in $\g^*$. We let $H$ denote the minimal nilpotent algebra
$H_{\chi}$. One of our main objectives is to give a presentation
of $H$ by generators and relations.

The action of the inner derivation $\ad h$ gives $\g$ a short
$\Z$-grading
$$\g\,=\,\g(-2)\oplus\g(-1)\oplus\g(0)\oplus\g(1)\oplus\g(2),\qquad\,
\g(i)=\{x\in\g\,|\,\,[h,x]=ix\}$$  with $\g(1)\oplus\ \g(2)$ and
$\g(-1)\oplus \g(-2)$ being Heisenberg Lie algebras. One knows of
course that $\g(\pm 2)$ is spanned by $ e_{\pm \beta}$, that
$\z_\chi(i)=\g(i)$ for $i=1,2$, and that
 $\z_\chi(0)$  coincides with  the image of the Lie algebra homomorphism
$$\sharp\,\colon\, \g(0)\longrightarrow\, \g(0),\quad\ x\mapsto
\,x-\frac{1}{2}(x,h)\,h$$ whose kernel $\k h$ is a central ideal
of $\g(0)$. The graded component $\g(-1)$ has a basis
$z_1,\ldots,z_s, z_{s+1},\ldots, z_{2s}$ such that the $z_i$'s
with $1\le i\le s$ (resp. $s+1\le i\le 2s$) are root vectors for
$\h$ corresponding to negative (resp. positive) roots, and
$$[z_i,z_j]=[z_{i+s},z_{j+s}]=0,\qquad
[z_{i+s},z_j]=\delta_{ij}f,\qquad (1\le i,j\le s).$$ Moreover, in
the present case we can choose $\m_\chi$ to be the span of $f$ and
the $z_i$'s with $s+1\le i\le 2s$, an abelian  subalgebra of $\g$
of dimension $s+1=\frac{1}{2}\dim {\cal O}_{\rm min}$. We set
$z_i^*:=z_{i+s}$ for $1\le i\le s$ and $z_i^*:=-z_{i-s}$ for
$s+1\le i\le 2s$.

Let $C$ denote the Casimir element of $U(\g)$ corresponding to the
bilinear form $(\,\cdot\,,\,\cdot\,)$. This form is nondegenerate
on $\z_\chi(0)$, hence we can find bases $\{a_i\}$ and $\{b_i\}$
of $\z_\chi(0)$ such that $(a_i,b_j)=\delta_{ij}$. Set
$\Theta_{\rm Cas}:=\sum_i\Theta_{a_i}\Theta_{b_i}$, a central
element of the associative subalgebra of $H$ generated by the Lie
algebra $\Theta(\z_\chi(0))$. Obviously, we can regard $C$ as a
central element of $H$.

By a well-known result of Joseph,
 outside type $\rm A$ the universal
enveloping algebra $U(\g)$ contains a unique completely prime
primitive ideal whose associated variety is $\overline{\cal
O}_{\rm min}$; see [\cite{J}]. This ideal, often denoted ${\cal
J}_0$, is known as the Joseph ideal of $U(\g)$.

We are finally in a position to formulate one of the main results
of this note:
\begin{theorem}\label{main3}
The algebra $H$ is generated by the Casimir element $C$ and the
subspaces $\Theta(\z_\chi(i))$ for $i=0,1$, subject to the
following relations:

\smallskip

\begin{enumerate}
\item[(i)\,] $[\Theta_x,\Theta_y]=\Theta_{[x,y]}$ for all
$x,y\in\z_\chi(0)$;

\medskip

\item[(ii)\,] $[\Theta_x,\Theta_u]=\Theta_{[x,u]}$ for all
$x\in\z_\chi(0)$ and $u\in\z_\chi(1)$;

\medskip

\item[(iii)\,] $C$ is central in $H$;

\medskip

\item[(iv)\,]
$[\Theta_u,\Theta_v]\,=\,\frac{1}{2}(f,[u,v])\big(C-\Theta_{\rm
Cas}-c_0\big)+
\frac{1}{2}\sum_{i=1}^{2s}\,\big(\Theta_{[u,z_i]^\sharp}\,\Theta_{[v,z_i^*]^\sharp}+
\Theta_{[v,z_i^*]^\sharp}\,\Theta_{[u,z_i]^\sharp}\big)$

\medskip

\noindent for all $u,v\in\z_\chi(1)$, where $c_0$ is a constant
depending on $\g$.
\end{enumerate}
\smallskip

\noindent If $\g$ is not of type $\rm A$ then $c_0$ is the
eigenvalue of $C$ on the primitive quotient $U(\g)/{\cal J}_0$. If
$\g$ is of type ${\rm A}_n$, $n\ge 2$, then
$c_0=-\frac{n(n+1)}{4}$. If $\g$ is of type ${\rm A}_1$ then
 $H=\k[C]$.
\end{theorem}
We start proving this theorem in Section~2 where we show that (i)
and (ii) hold in $H$ for a suitable choice of $\Theta\colon
\z_\chi\rightarrow H$. In Section~4 we determine all of the
quadratic relation (iv) except the elusive constant $c_0$.

\subsection{} We first computed $c_0$ by brute force, but later it turned out that
there was a much better way to do it, based on a certain refined
version of Joseph's Preparation Theorem. This theorem which we
prove in Section~5 in our special case, enables us to link the
primitive ideals of $H$ directly with primitive ideals of $U(\g)$.

Let $\Delta$ denote the automorphism of the polynomial algebra
$\k[h]$ such that $\Delta(h)=h+1$. Let $\langle \Delta\rangle$
stand for the cyclic subgroup of $\text{Aut}(\k[h])$ generated by
$\Delta$. The skew group algebra $\k[h]*\langle\Delta\rangle$ has
$\{h^i\Delta^j\,|\,i\in\Z_+,\,j\in\Z\}$  as a $\k$-basis and
multiplication in $\k[h]*\langle\Delta\rangle$ has the property
that $\Delta\cdot h=(h+1)\cdot\Delta$.

Let ${\mathbf A}_e$ denote the Weyl algebra with standard
generators $z_1,\ldots,z_s,\partial_1,\ldots,\partial_s$, so that
$[\partial_i,z_j]=\delta_{ij}$ for $1\le i,j\le s$. Let ${\mathcal
A}_e:=\big(\k[h]*\langle\Delta\rangle\big)\otimes {\mathbf A}_e$,
a simple Noetherian algebra over $\k$, and identify
$\k[h]*\langle\Delta\rangle$ and ${\mathbf A}_e$ with subalgebras
of ${\mathcal A}_e$. Define an involution
$\tau\in\text{Aut}({\mathcal A}_e)$ by setting
$$\tau(z_i)=-z_i,\quad
\tau(\partial_i)=-\partial_i,\quad\tau(h)=h,\quad\tau(\Delta^k)=(-1)^k\Delta^k
\qquad(1\le i\le s,\,k\in\Z).$$ Then $\tau\otimes\sigma$ is an
automorphism of order two of the associative algebra ${\mathcal
A}_e\otimes H$.

Let $U(\g)_f$ denote the localisation of $U(\g)$ with respect to
the Ore set $\{f^i\,|\,i\in \Z_+\}$. By mapping $U(\g)_f$ into the
endomorphism algebra of the induced module $Q_\chi$ we are able to
identify $U(\g)_f$ with a subalgebra of ${\mathcal A}_e\otimes H$.
More precisely, we prove that
$$U(\g)_f\,=\,({\mathcal
A}_e\otimes H)^{\tau\otimes\sigma}\,=\,{\mathcal A}_e^\tau\otimes
H_+\oplus{\mathcal A}_e^\tau\Delta\otimes H_-$$ where
$H_{\pm}=\{x\in H\,|\,\sigma(x)=\pm x\}$. As mentioned in (1.1)
every two-sided ideal $I$ of $H$ is stable under the involution
$\sigma\in\text{Aut}(H)$. Hence $I=I_+\oplus I_-$ where
$I_{\pm}=I\cap H_{\pm}$. We identify $U(\g)$ with a subalgebra of
$U(\g)_f$ and set
$$\widetilde{I}\,:=\,U(\g)\cap ({\mathcal A}_e^\tau\otimes
I_+\oplus{\mathcal A}_e^\tau\Delta\otimes I_-).$$ Then
$\widetilde{I}$ is a two-sided ideal of $U(\g)$. By
Corollary~5.1(vi), the centre of $H$ identifies canonically with
$Z(\g)$, the centre of $U(\g)$. Let ${\mathcal
X}=\text{Prim}\,U(\g)$ and let ${\mathcal X}_{\rm inf}$ be the set
of all primitive ideals of infinite codimension in $\mathcal X$.
Given a prime Noetherian ring $R$ we let $\text{rk}(R)$ denote the
Goldie rank of $R$.

\begin{theorem}
Take ${\mathrm Prim}\,H$ with the Jacobson topology and take
${\mathcal X}_{\rm inf}$ with the topology induced by the Jacobson
topology of $\mathcal X$. Then the following hold:
\smallskip

\begin{enumerate}
\item [(i)] The map $I\mapsto\widetilde{I}$ induces a
homeomorphism $\varkappa\colon\, {\mathrm
Prim}\,H\stackrel{\sim}{\longrightarrow}\, {\mathcal X}_{\mathrm
inf}$.

\smallskip

\item[(ii)] For any $I\in\text{Prim}\,H$ we have  ${\mathrm
Dim}(U(\g)/\widetilde{I})\,=\,{\mathrm Dim}(H/I)+\dim\,
\overline{\cal O}_{\rm min}$.

\smallskip

\item [(iii)] If $I={\mathrm Ann}_H\, V$ where $V$ is a finite
dimensional irreducible $H$-module, then $\widetilde{I}={\mathrm
Ann}_{U(\g)}\big(Q_\chi\otimes_H V\big)$ and $\,\,{\mathrm
rk}(U(\g)/\widetilde{I})\,=\,\dim{\mathrm Wh}(Q_\chi\otimes_H
V)\,=\,\dim V.$

\smallskip

\item [(iv)] For any  $\cal I\in{\mathcal X}$ with ${\mathcal
VA}({\cal I})\,=\,\overline{{\cal O}}_{\mathrm min}$ there is a
finite dimensional irreducible $H$-module $V$ such that ${\cal
I}=\,{\rm Ann}_{U(\g)}\big(Q_\chi\otimes_H V\big)$.

\smallskip

\item [(v)] Let $V_1$ and $V_2$ be two finite dimensional
irreducible $H$-modules. Then $V_1\cong V_2$ as $H$-modules if and
only if $\,{\mathrm Ann}_{U(\g)}\big(Q_\chi\otimes_H
V_1\big)\,=\,{\mathrm Ann}_{U(\g)}\big(Q_\chi\otimes_H V_2\big).$

\smallskip

\item [(vi)] A prime ideal $I$ of $H$ is primitive if and only if
$I\cap Z(H)$ is a maximal ideal of $Z(H)$.

\end{enumerate}
\end{theorem}

 It follows from Theorem~1.2 that for any  homomorphism
$\eta\colon Z(\g)\rightarrow \k$ the map $\varkappa$ induces a
bijection between the isoclasses of finite dimensional irreducible
$H$-modules with central character $\eta$ and the primitive ideals
${\cal I}\in{\mathcal X}$ such that ${\cal I} \cap
Z(\g)\,=\,{\mathrm Ker}\,\eta$ and ${\mathcal VA}({\cal
I})\,=\,\overline{{\cal O}}_{\mathrm min}$ (recall that
$Z(\g)=Z(H)$). This result indicates that for any nilpotent $\chi$
it should be possible to interpret the number of isoclasses of
irreducible finite dimensional $H_\chi$-modules with a fixed
central character as the dimension of a cell representation of the
integral Weyl group of the character. In type $\mathrm A$ this
agrees with recent results of Brundan--Kleshchev.

Theorem~1.2(iii) relates the dimensions of irreducible finite
dimensional $H$-modules with Goldie-rank polynomials. We explore
this in (6.4) to obtain dimension formulae for all irreducible
finite dimensional representations of $H$ for $\g$ of type ${\rm
C}_n$ and ${\rm G}_2$. It is quite possible that {\it all}
Goldie-rank polynomials, properly scaled, will appear in dimension
formulae for ``nonrestricted Weyl modules'' over Lie algebras of
reductive groups in characteristic $p$ (we recall that in
characteristic $p$ a truncated version of $H_\chi$ is Morita
equivalent to the reduced enveloping algebra $U_\chi(\g)$; see
[\cite{P02}, (2.3), (2.6)]).

Theorem~1.2(vi) says that $H$ satisfies the
Dixmier--M{\oe}glin--Rentschler equivalence. Again it is possible
that this holds for any algebra $H_\chi$.

\subsection{} In the last section of this note we introduce and
study highest weight modules for the algebra $H$. Let
$\Phi_e=\{\alpha\in\Phi\,|\,\alpha(h)=0\mbox{ or } 1\}$, and put
$\Phi_e^\pm=\Phi_e\cap \Phi^\pm$ where $\Phi^-=-\Phi^+$. For
$i=0,1$ put
$\Phi_{e,i}^{\pm}=\{\alpha\in\Phi_e^\pm\,|\,\alpha(h)=i\}$. Note
that $\z_\chi$ is spanned by $\h_e:=\h\cap\g(0)^\sharp$, by root
vectors $e_\alpha$ with $\alpha\in\Phi_e$, and by $e$.  Let
$h_1,\ldots, h_{l-1}$ be a basis of $\h_e$, and let $\n^\pm(i)$ be
the span of all $e_\alpha$ with $\alpha\in\Phi^\pm_{e,i}$.
Clearly, $\n^+(0)$ and $\n^-(0)$ are maximal nilpotent subalgebras
of $\g(0)^\sharp$. Let $\{x_1,\dots ,x_t\}$ and
$\{y_1,\dots,y_t\}$ be bases of $\n^+(0)$ and $\n^-(0)$ consisting
of root vectors for $\h$. For $1\le i\le s$ let $\gamma_i$ (resp.
$\gamma^*_i$) denote the root of $z_i$ (resp. $z_i^*$), and put
$u_i=[e,z_i],\,$ $u_i^*=[e,z_i^*]$. Then  $\{u_1,\ldots,
u_s,u_1^*,\ldots,u_s^*\}$ is a $\k$-basis of $\z_\chi(1)$.

In general, $H$ is unlikely to possess a triangular decomposition
similar to that of $U(\g)$. Nevertheless, one can still define
Verma modules and highest weight modules for $H$. Given
$\lambda\in\h_e^*$ and $c\in\k$ we denote by $J_{\lambda,c}$ the
linear span in $H$ of all
$$\prod_{i=1}^t\Theta_{y_i}^{l_i}\cdot
\prod_{i=1}^s\Theta_{u_i}^{m_i} \cdot\prod_{i=1}^{\ell-1}
\big(\Theta_{h_i}-\lambda(h_i)\big)^{n_i}\cdot (C-c)^{n_\ell}\cdot
\prod_{i=1}^s\Theta_{u_i^*}^{r_i}\cdot
\prod_{i=1}^t\Theta_{x_i}^{q_i}$$ with
$\sum_{i=1}^{\ell}n_i+\sum_{i=1}^t r_i+\sum_{i=1}^s q_i>0$. Using
Theorem~1.1 we show in (7.1) that $J_{\lambda,c}$ is a left ideal
of $H$. We call the $H$-module $Z_H(\lambda,c):=\,H/J_{\lambda,c}$
the {\it Verma module of level} $c$ corresponding to $\lambda$. By
the above, $Z_H(\lambda,c)$ has a nice PBW basis. In (7.2) we show
that $Z_H(\lambda,c)$ contains a unique maximal submodule which we
denote $Z_H^{\rm max}(\lambda,c)$. Thus to every
$(\lambda,c)\in\h_e^*\times \k$ there corresponds an irreducible
highest weight $H$-module
$L_H(\lambda,c):=\,Z_H(\lambda,c)/Z_H^{\rm max}(\lambda,c)$. It is
fairly easy to show that $L_H(\lambda,c)\cong L_H(\lambda',c')$ if
and only if $(\lambda,c)=(\lambda',c')$ and that any irreducible
finite dimensional $H$-module is isomorphic to exactly one of
$L_H(\lambda,c)$ with $\lambda$ satisfying a natural integrality
condition.

To determine the composition multiplicities of the Verma modules
$Z_H(\lambda,c)$ we link them with $\g$-modules obtained by
parabolic induction from Whittaker modules for ${\frak sl}(2)$.
Let ${\frak s}_\beta\,=\,\k e_\beta\oplus\k h_\beta\oplus\k
f_\beta$ and put
$$\p_\beta:=\,{\frak s}_\beta+\h+\textstyle{\sum}_{\alpha\in\Phi^+}\,\,
\k e_\alpha,\quad\,
\n_\beta:=\,\textstyle{\sum}_{\alpha\in\Phi^+\setminus\{\beta\}}\,\,\k
e_\alpha,\quad\, \widetilde{\frak s}_\beta:=\,\h_e\oplus
{\mathfrak s}_\beta.$$  Let
$C_\beta=ef+fe+\frac{1}{2}h^2=2ef+\frac{1}{2}h^2-h$, a central
element of $U({\widetilde{\mathfrak s}}_\beta)$. Given
$\lambda\in\h_e^*$ and $c\in\k$ we denote by $I_\beta(\lambda,c)$
the left ideal of $U(\p_\beta)$ generated by $f-1,C_\beta-c$, all
$h-\lambda(h)$ with $h\in\h_e$, and all $e_\gamma$ with
$\gamma\in\Phi^+\setminus\{\beta\}$. Let
$Y(\lambda,c):=\,U(\p_\beta)/I_\beta(\lambda,c)$, a
$\p_\beta$-module with the trivial action of $\n_\beta$. Regarded
as an ${\frak s}_\beta$-module, $Y(\lambda,c)$ is isomorphic to a
Whittaker module for $\sl(2,\k)$. Now define
$$M(\lambda,c):=\,U(\g)\otimes_{U(\p_\beta)}
Y(\lambda,c).$$  Recall that each $z_i^*$ with $i\le s$ is a root
vector corresponding to $\gamma_i^*=-\beta-\gamma_i\in\Phi^+$. Let
$\delta=\frac{1}{2}(\gamma_1^*+\cdots+\gamma_s^*)$ and
$\rho=\frac{1}{2}\sum_{\alpha\in\Phi^+}\alpha.$  Since the
restriction of $(\,\cdot\,,\,\cdot\,)$ to $\h_e$ is nondegenerate,
for any $\eta\in\h_e^*$ there is a unique $t_\eta\in \h_e$ such
that $\varphi=(t_\eta,\,\cdot\,)$. Hence $(\,\cdot\,,\,\cdot\,)$
induces a bilinear form on $\h_e^*$ via
$(\mu,\nu):=(t_{\mu},t_{\nu})$ for all $\mu,\nu\in\h_e^*$. Given a
linear function $\varphi\in\h^*$ we denote by $\bar{\varphi}$ the
restriction of $\varphi$ to $\h_e$.
\begin{theorem}
Each $\g$-module $M(\lambda,c)$ is an object of the category
${\mathcal C}_\chi$. Furthermore, ${\mathrm Wh}(M(\lambda,c))\cong
Z_H(\lambda+\bar{\delta},c+(\lambda+2\bar{\rho},\lambda))$ as
$H$-modules.
\end{theorem}
Combined with Skryabin's equivalence  and the main results of
Mili{\v c}i{\'c}--Soergel [\cite{MS}] and Backelin [\cite{Bac}],
Theorem~1.3 shows that the composition multiplicities of the Verma
modules $Z_H(\lambda,c)$ can be computed with the help of certain
parabolic Kazhdan-Lusztig polynomials. This confirms in the
minimal nilpotent case the Kazhdan-Lusztig conjecture for finite
$\mathcal W$-algebras formulated by de~Vos and van~Driel in
[\cite{DvvD}]; see Remark~\ref{remlast} for more detail.

Apart from its relevance to the theory of primitive ideals this
work is a contribution to the rapidly growing theory of $\mathcal
W$-algebras. {\it Finite} $\mathcal W$-algebras are attached to
nilpotent elements of finite dimensional simple Lie algebras via
quantum Hamiltonian reduction. All finite $\mathcal W$-algebras of
type $\mathrm A$ were recently described by J.~Brundan and
A.~Kleshchev [\cite{BrK}] who identified them with shifted
truncated Yangians. It seems likely that their results can be
extended to some nilpotent elements in Lie algebras of types
${\mathrm B}$, ${\mathrm C}$ and ${\mathrm D}$. Hidden Yangian
symmetry of finite $\mathcal W$-algebras of type $\mathrm A$ was
first discovered, in some special cases, by E.~Ragoucy and
P.~Sorba [\cite{RS}].

{\it Affine} counterparts of finite $\mathcal W$-algebras have
been studied even more intensively. It should be mentioned here
that V.G.~Kac and M.~Wakimoto described {\it minimal nilpotent}
superconformal algebras in the context of vertex operators  and
quantum reduction; see [\cite{KW}] and the references therein. It
would be interesting to compare the algebras $H$ of this paper
with quasiclassical limits of vertex algebras of Kac--Wakimoto.

\smallskip

\noindent {\bf Acknowledgement.} I would like to thank Jonathan
Brundan, Victor Ginzburg, Anthony Joseph and Alexander Kleshchev
for interesting discussions and e-mail correspondence. Some
results of this work were announced in my talks at the AMS Summer
Research Conference in Snowbird, Utah in July 2004, at the
Oberwolfach meeting on enveloping algebras in March 2005, and at
the Luminy conference ``Geometry and Representations'' in April
2005.
\section{\bf Structural features of the algebras $H_\chi$}
\subsection{}
In this section we assume that $e$ is an arbitrary nilpotent
element in $\g$.  Decompose $\g$ into the weight spaces relative
to $\ad h$ giving a $\Z$-grading
$\g\,=\,\bigoplus_{i\in\Z}\,\g(i)$. Let $\chi$ be as in (1.1) and
denote by $\z_\chi$ the centraliser of $\chi$ in $\g$. It is
well-known  that $\z_\chi={\frak c}_{\g}(e)$ is a graded
subalgebra of $\p_{e}\,:=\,\bigoplus_{i\ge 0}\,\g(i)$, that is
$\z_{\chi}\,=\,\bigoplus_{i\ge 0}\,\z_{\chi}(i)$ where
$\z_{\chi}(i)=\z_{\chi}\cap\g(i)$. Choose a $\k$-basis
$x_1,\ldots, x_m$ of the parabolic subalgebra ${\frak p}_e$ with
$x_i\in\g(n_i)$ such that $x_1, \ldots, x_r$ span $\z_\chi$. Let
${\cal O}={\cal O}_\chi$ and let $d$ denote half of the dimension
of $\cal O$.

Define the skew-symmetric bilinear form $\la\cdot\,,\cdot \ra$ on
the subspace $\g(-1)$ by setting $\la x,y\ra=(e,[x,y])$ for all
$x,y\in\g(-1)$. As $\z_\chi\subset\p_e$, this form is
nondegenerate. Choose a  basis $z_{1},\ldots z_{s}, z_{s+1},\ldots
z_{2s}$ of $\g(-1)$ such that
$$\la z_{i+s},z_j\ra=\delta_{ij},\qquad \la z_i,z_j\ra=\la z_{i+s},z_{j+s}\ra=0
\qquad (1\le i,j\le r)$$ and denote by $\g(-1)^0$ the linear span
of $z_{s+1},\ldots, z_{2s}$. Let $\m_\chi\,=\,\g(-1)^0\oplus\,
\sum_{i\le 2}\,\g(i)$,  a nilpotent subalgebra of $\g$ of
dimension $d$; see [\cite{P02}] for example. Since $\chi$ vanishes
on the derived subalgebra of $\m_\chi$ the  ideal $N_\chi$ of
$U(\m_\chi)$ generated by all $x-\chi(x)$ with $x\in\m_\chi$ has
codimension one in the enveloping algebra $U(\m_\chi)$. Let
$\k_\chi\,=\,U(\m_\chi)/N_\chi$, a one-dimensional left
$U(\m_\chi)$-module, and let $1_\chi$ stand for the image of $1$
in $\k_\chi$. We denote by $Q_\chi$ the induced $\g$-module
$U(\g)\otimes_{\,U(\m_\chi)} \k_\chi$ and set
$$H_\chi:={\End}_{\g\,}(Q_\chi)^{\text{op}}.$$   It is proved in
[\cite{P02}] that the algebra $H_\chi$ is a filtered deformation
of the graded coordinate ring $\k[{\mathcal S}_e]$ .

In what follows we will rely on a different realisation of
$H_\chi$ found by W.L.~Gan and Ginzburg [\cite{GG}]. Let
$\n_\chi\,=\,\bigoplus_{i\le -1}\,\g(i)$ and
$\n_\chi'\,=\,\bigoplus_{i\le -2}\,\g(i)$. Clearly, $\n_\chi$ and
$\n_\chi'$ are nilpotent subalgebras of $\g$  and $\n'_\chi$ is an
ideal of $\n_\chi$. Since $\n_\chi'\subseteq \m_\chi$ we may view
$\k_\chi$ as an $\n'_\chi$-module. Let
$\widehat{Q}_\chi\,=\,U(\g)\otimes_{\,U(\n'_\chi)}\k_\chi$, an
induced $\g$-module and the quotient of  $U(\g)$ by the left ideal
${\mathcal I}_\chi$ generated by all  $x-\chi(x)$ with
$x\in\n'_\chi$. The representation of $U(\g)$ in $\End(Q)$ will be
denoted by $\widehat{\rho}_\chi$. Since $\chi$ vanishes on
$[\n_\chi,\n'_\chi]\subseteq \bigoplus_{i\le -3}\,\g(i)$, the left
ideal ${\mathcal I}_\chi$ is stable under the adjoint action of
$\n_\chi$ on $U(\g)$. Therefore, $\ad\,\n_\chi$ acts on
$\widehat{Q}_\chi$. The fixed point space
$\widehat{Q}_\chi^{\,\ad\,\n_\chi}$ carries a natural algebra
structure given by $(x+{\mathcal I}_\chi) (y+{\mathcal
I}_\chi)=xy+{\mathcal I}_\chi$ for all $x+{\mathcal I}_\chi,
y+{\mathcal I}_\chi\in \widehat{Q}_\chi$; see [\cite{GG}, p.~244]
for more detail. We furnish $Q_\chi^{\,\ad\,\m_\chi}$  and
$\widehat{Q}_\chi^{\,\ad\,\n'_\chi}$ with algebra structures  in a
similar fashion. It is well known (and easily seen) that
$$H_\chi\cong {Q}_\chi^{\,\ad\,\m_\chi} \quad\,  \mbox {  and
}\quad\, {\End}_{\g}(\widehat{Q}_\chi)^{\text{op}}\cong
\widehat{Q}_\chi^{\,\ad \n'_\chi}$$ as algebras. As
$\n'_\chi\subseteq \m_\chi$, there is a natural $\g$-module
epimorphism $\widehat{Q}_\chi\twoheadrightarrow Q_\chi$. As
$\m_\chi \subseteq \n_\chi$,  it induces an algebra map
$\eta\colon\,\widehat{Q}_\chi^{\,\ad\,\n_\chi}\map H_\chi$. By
[\cite{GG}, Theorem~4.1], $\eta$ is an isomorphism of algebras.
Henceforth  we will make no distinction between $H_\chi$ and
$\widehat{Q}_\chi^{\,\ad \n_\chi}$ and view the latter as a
subalgebra of $\End_{\g}(\widehat{Q})^{\text{op}}$.

Given  $({\bf a},{\mathbf b})\in\Z_+^m\times \Z_+^{2s}$ we set
$x^{\bf a}z^{\bf b}:= x_1^{a_1}\cdots x_m^{a_m}z_{1}^{b_1}\cdots
z_{2s}^{b_{2s}}$. By the PBW theorem, the monomials $x^{\bf
a}z^{\bf b}\otimes 1_\chi$ with $({\bf a},{\bf b})\in
\Z_+^{m}\times\Z_+^{2s}$ form a $\k$-basis of  $\widehat{Q}_\chi$.
For $k\in\Z_+$ we denote by $\widehat{Q}_{\chi}^{\,k}$ the linear
span of all $x^{\mathbf a}z^{\mathbf b}\otimes 1_\chi$ with
\begin{equation}\label{GK1}
|({\bf a}, {\bf b})|_e:=\sum_{i=1}^m a_i(n_i+2)+
\sum_{i=1}^{2s} b_i\le k.
\end{equation}
We let $H_{\chi}^k$ denote the subspace of $H_\chi$ consisting of
all $h\in H_\chi$ with $h(1_\chi)\in \Q^{\,k}$. By [\cite{P02}]
and [\cite{GG}], the subspaces $\{H_{\chi}^k\,|\,k\in\Z_+\}$ form
an increasing filtration of the algebra $H_\chi$ and the
corresponding graded algebra $\gr\,H_\chi$ is isomorphic to a
polynomial ring in $r$ variables with free homogeneous generators
of degree $n_1+2,\ldots,n_r+2$. The elements $x$ in
$\Q^{\,k}\setminus \Q^{\,k-1}$ and  $H_{\chi}^k\setminus
H_{\chi}^{k-1}$ are said to have {\it Kazhdan degree} $k$, written
$\deg_{e}(x)=k$. It is immediate from [\cite{P02},  Theorem~4.6]
that in our present realisation the algebra $H_\chi$ has a
distinguished generating set $\Theta_1,\ldots,\Theta_r$ such that
$\gr\,\Theta_1, \ldots, \gr\,\Theta_r$ generate $\gr\,H_\chi$ and
\begin{equation}\label{GK2}
\Theta_{k}(1_\chi)\,= \,\Big(x_k+
\sum_{0< |({\mathbf i},{\mathbf j})|_e\le n_k+2}\,
\lambda^{k}_{{\mathbf i},{\mathbf j}}\,x^{\mathbf i}z^{\mathbf j}
\Big)\otimes 1_{\chi},\,
\quad 1\le k\le r,
\end{equation}
where  $\lambda^{k}_{{\mathbf i},{\mathbf j}}\in\k$ and
$\lambda^{k}_{\bf{i},\bf{j}}=0$  if either $|({\bf i},{\bf
j})|_e>n_k+2$ or $|({\bf i},{\bf j})|_e=n_k+2$ and $|{\bf
i}|+|{\bf j}|=1$ or ${\bf i}\ne {\bf 0}$, $\,{\bf j}={\bf 0}$, and
$i_j=0$ for $j>r$. The monomials $\Theta_1^{a_1}\cdots
\Theta_r^{a_r}$ with $(a_1,\ldots, a_r) \in\Z_+^r$ form a PBW
basis of $H_\chi$.

\subsection{}\label{2.2}
Given a subset $X$ of $\g$ we denote by $Z_G(X)$   the closed
subgroup of $G$ consisting of all $g\in G$ with $(\Ad g)( x)=x$
for all $x\in X$. Let $P_e$ denote the parabolic subgroup of $G$
with $\Lie\,P_e=\p_e$. There exists a $1$-parameter subgroup
$\lambda_e\colon\, \k^\times\rightarrow G$ optimal for the
$G$-unstable vector $e$ and such that:

\smallskip

\begin{itemize}

\item[$\bullet$\,] $(\Ad \lambda_e(t))_{\vert
\g(i)}\,=\,t^i\,\text{id}$ for all $t\in \k^\times$ and $i\in\Z$;

\medskip

\item[$\bullet$\,] $Z_G(e)\subset P(e)$, $\,R_u(Z_G(e))\subset
R_u(P_e)$,  $\,Z_G(e)\,=\big(Z_G(e)\cap Z_G(\lambda_e)\big)
R_u(Z_G(e))$;

\medskip

\item[$\bullet$\,] $C(e):=Z_G(e)\cap Z_G(\lambda_e)$ is  a
reductive group, and $\Lie\, C(e)=\z_\chi(0);$

\end{itemize}

\medskip

\noindent see [\cite{Ca}, Chapter~5] and [\cite{P03}]. Let $\Ad
C(e)$ denote the image of $C(e)$ in the adjoint group $\Ad
G=(\text{Aut}\,\g)^\circ$. Put $\sigma=\Ad \lambda_e(-1)$, an
element of order $\le 2$ in $\Ad G$. Clearly, $\sigma$ lies in the
centre of $\Ad C(e)$ and $\sigma(x)=(-1)^i\, x$ for all
$x\in\g(i)$ and $i\in\Z$.
\begin{lemma} \label{L}
The element $\sigma$ belongs to any maximal torus of $\Ad C(e)$.
\end{lemma}
\begin{pf}
Let $T_0$ be a maximal torus of $\Ad C(e)$, $\widetilde{T}_0$ the
inverse image of $T_0$ in $G$, and  $L=Z_G(\widetilde{T}_0)$. Then
$L$ is a Levi subgroup of $G$ and $e$ is a distinguished nilpotent
element in ${\frak l}=\Lie\,L$. The construction in [\cite{P03}]
shows that all weights of $\Ad \lambda_e(\k^\times)$ on $\frak l$
are even. Then $\sigma$ acts trivially on $\frak l$,  yielding
$\g^{\,\sigma t}\supseteq \frak l$ for all $t\in T_0$. As $\k$ is
infinite, there is $t_0\in T_0$ such that $\g^{\,\sigma t_0}=\frak
l$. Let $\cal C$ denote the conjugacy class of the image of
$\sigma t_0$ in the component group $Z_{\Ad G}(e)/Z_{\Ad
G}(e)^\circ\cong (\Ad C(e))/(\Ad C(e))^\circ$. As $\Ad G$ is a
group of adjoint type, the $G$-conjugacy class of the pair $(L,e)$
corresponds under Sommers' bijection  to the $G$-conjugacy class
of the pair $(e,\cal C)$; see [\cite{S}, \cite{McS}, \cite{P03}].
As $L$ is a Levi subgroup in $G$, {\it loc. cit.} also shows that
${\cal C}=\{1\}$. But $T_0\subseteq (\Ad C(e))^\circ$ and $t_0\in
T_0$. So we get $\sigma\in Z\big((\Ad C(e))^\circ\big)$. As $(\Ad
C(e))^\circ$ is a reductive group, the torus $T_0$ is
self-centralising in $(\Ad C(e))^\circ$. Hence $\sigma\in T_0$
completing the proof.
\end{pf}

We now fix a maximal torus $T_e$ in $\Ad C(e)$ and assume (without
loss of generality) that all $z_i$ with $i\le 2s$ and $x_j$ with
$j\le m$  are weight vectors with respect to $T_e$. By
Lemma~\ref{L}, $\sigma\in T_e$. Note that $C(e)$ preserves both
$\n'_\chi$ and ${\mathrm Ker}\,\chi$. Since $C(e)$ acts on $U(\g)$
as algebra automorphisms, it  preserves the left ideal ${\mathcal
I}_\chi$ and thus acts on $\Q$. This action is compatible with
that of $\g$, i.e.
\begin{equation}\label{GK2'}
g \circ {\widehat{\rho}}_{\chi} (x)\circ g^{-1}=\,\widehat{\rho}_\chi\big((\Ad g)(x)\big)\ \ \,
\qquad \big(\forall\, g\in C(e),  \, \, x\in\g\big).
\end{equation}
Since $C(e)$ preserves $\n_\chi$ too, it acts on
$H_\chi=\Q^{\,\ad\,\n_\chi}$ as algebra automorphisms.  Since
$g(1_\chi)=1_\chi$ for all $g\in C(e)$, the action of $C(e)$ on
$\Q$ and $H_\chi$ is filtration preserving, hence locally finite.
Since $Z(G)$ acts trivially on $U(\g)$,  there is a natural action
of $\Ad C(e)$ on $\widehat{Q}_\chi$ and $H_\chi$. It should be
mentioned that
\begin{equation}\label{GK2''}
\sigma(x^{\bf a}z^{\bf b}\otimes 1_\chi)\,=\,(-1)^{|({\bf a},{\bf b})|_e}\,
x^{\bf a}z^{\bf b}\otimes 1_\chi
\end{equation}
for all $({\bf a},{\bf b})\in\Z_+^m\times \Z_+^{2s}$.
\begin{lemma}\label{L1} Each generator
$\Theta_k\in H_\chi$ can be chosen to be a weight vector for $T_e$
of the same weight as $x_k$.
\end{lemma}
\begin{pf} Let
$\gamma_k$ denote the $T_e$-weight of $x_k$. If $\gamma_k\ne 0$, we
assume without
loss of generality that $\lambda^k_{{\bf 0},{\bf 0}}=0$.
Let $t\in T_e$ and $\,\Theta_t:=t(\Theta_k)-\gamma_k(t)\Theta_k$, an element
in $H_\chi$. Since all $x^{\bf i}z^{\bf j}\otimes 1_\chi\in\Q$ are weight vectors for $T_e$,
we deduce from
(\ref{GK2}) and (\ref{GK2'}) that $\Theta_t(1_\chi)$
is a linear combination of $x^{\bf a}z^{\bf b}\otimes 1_\chi$ with either ${\bf b}\ne {\bf 0}$
or $a_j\ne 0$ for some $j>r$. Then $\Theta_t=0$ for all
$t\in T_e$,  by [\cite{P02}, Lemma~4.5], and the result follows.
\end{pf}
\subsection{} We now consider the linear map
$\Theta\colon\z_\chi\,\longrightarrow\, H_\chi,\,$ $x\mapsto \Theta_x,$
such that $\Theta_{x_i}=\Theta_i$ for all $i$. Thanks to Lemma~\ref{L1},
$\Theta$ is an injective homomorphism of $T_e$-modules.
Although $\Theta$ is not a Lie algebra homomorphism, in general,
it follows from  [\cite{P02}, Theorem~4.6(iv)] and (\ref{GK2''}) and Lemma~\ref{L1}
that
\begin{equation}\label{GK3}
[\Theta_{x_i},\Theta_{x_j}]\,\equiv\,
\Theta_{[x_i,x_j]}+q_{ij}(\Theta_1,\ldots,\Theta_r)\ \ \,
\big({\mathrm mod}\ \,  H_{\chi}^{n_{i}+n_{j}}\big)
\end{equation}
where $q_{ij}$ is a polynomial in $r$ variables  with initial form
of total degree $\ge 2$.

\begin{rem} As $C(e)$ is a reductive group, each $C(e)$-module $H^k_\chi$ is completely
reducible. From this it follows that there exists a unitriangular
polynomial substitution
$$F\colon\,(\Theta_1,\ldots,\Theta_r)\,\,\longmapsto\,\big(F_1(\Theta_1,\cdots,\Theta_r),
\ldots,F_r(\Theta_1,\ldots,\Theta_r)\big)$$
which satisfies the following conditions:
\smallskip

\begin{itemize}
\item[$\bullet\, $]
$\deg_e\,F(\Theta_i)=\deg_e\,\Theta_i=n_i+2$ for all $i\le r$;

\medskip

\item[$\bullet\, $]
the linear map $\Theta_F\colon\,\z_\chi\,\rightarrow\, H_\chi$ with
$\Theta_F(x_i)=F(\Theta_i)$  for all $i$ is an injective homomorphism of
$C(e)$-modules, and $\z_\chi\cong \,\gr\,\Theta_F(\z_\chi)$
as graded $C(e)$-modules;

\medskip

\item[$\bullet\, $] an analogue of (\ref{GK3}) holds for
$F(\Theta_1),\ldots, F(\Theta_r)$ and the subspaces
$\Theta_F(\z_\chi)$, and $\gr\,\Theta_F(\z_\chi)$  generate the
algebras $H_\chi$ and $\gr\, H_\chi$, respectively.
\end{itemize}
\end{rem}

\begin{prop} \label{P1}
There exists an associative $\k[t]$-algebra ${\mathcal H}_\chi$  free
as a module over $\k[t]$ and such that
\[{\mathcal H}_{\chi}/(t-\lambda){\mathcal H}_{\chi}\cong
\left\{\begin{array}{ll} H_\chi &  \mbox{ if $\, \lambda\ne 0$},\\
U(\z_\chi) & \mbox{ if $\, \lambda=0$}
\end{array}
\right. \] as
$\k$-algebras.
In other words, the enveloping algebra $U(\z_\chi)$ is a contraction of
$H_\chi$.
\end{prop}
\begin{pf}
Consider the algebra $H(R)=R\otimes H_\chi$ over the ring of
Laurent polynomials $R=k[t,t^{-1}]$ obtained from $H_\chi$ by
extension of scalars, and identify $H_\chi$ with the subspace
$\k\otimes H_\chi$ of $H(R)$. Define an invertible $R$-linear
transformation $\pi$ on $H(R)$ by setting
$$\pi(\Theta_1^{k_1}\cdots\Theta_r^{k_r})=t^{n_1k_1+\ldots+n_rk_r}
\Theta_1^{k_1}\cdots\Theta_r^{k_r} \quad \ \  \forall\,
(k_1,\ldots, k_r)\in\Z_+^r$$ and extending to $H(R)$ by
$R$-linearity. We view $\pi$ as an isomorphism from $H(R)$ onto a
new $R$-algebra $H(R,\pi)$ with underlying $R$-module $R\otimes
H_\chi$ and with associative product  given by $(x\cdot y)_\pi:=\,
\pi^{-1}\big(\pi(x)\cdot\pi(y)\big)$ for all $x,y\in R\otimes
H_\chi$. We denote by ${\mathcal H}_\chi$ the free
$\k[t]$-submodule of $H(R,\pi)$ generated by
$\Theta_1^{a_1}\cdots\Theta_r^{a_r}$ with
$(a_1,\ldots,a_r)\in\Z_+^r$. It follows from (\ref{GK1}) and
(\ref{GK2}) that
$$\deg_e\big(\Theta_1^{k_1}\cdots\Theta_r^{k_r}\big)\,=\,
\sum_{i=1}^r n_{i}k_{i}+2\sum_{i=1}^r k_i.$$ In view of (\ref{GK3}) this yields
$$(\Theta_i\cdot\Theta_j-\Theta_j\cdot\Theta_i)_\pi\,=\,\pi^{-1}\big(t^{n_i+n_j}
[\Theta_i,\Theta_j]\big)\,\equiv\,\Theta_{[x_i,x_j]}\ \, \,
\big({\mathrm mod}\, \,  t{\mathcal H}_{\chi}\big)$$ (because the
initial form of $q_{ij}$ has total degree $\ge 2$ and $\deg_e
q_{ij}(\Theta_1,\ldots,\Theta_r) =n_i+n_j+2$ if $q_{ij}\ne 0$).
Using induction on the Kazhdan degree of $\Theta_1^{k_1}\cdots
\Theta_r^{k_r}$ and the commutativity of $\gr\,H_\chi$ we now
deduce that $(\Theta_i\cdot {\mathcal H}_\chi)_\pi
\subseteq{\mathcal H}_\chi$ for all $i$. So ${\mathcal H}_\chi$ is
a $\k[t]$-subalgebra of $H(R,\pi)$.

If $\lambda\ne 0$ then the homomorphism $\k[t]\rightarrow k$
taking $t$ to $\lambda$ extends to a homomorphism $R\rightarrow
k$. The isomorphism $\pi^{-1}$ injects $(t-\lambda)H(R,\pi)$ onto
$(t-\lambda)H(R)$. Because ${\mathcal H}_\chi\cap
(t-\lambda)H(R,\pi)=(t-\lambda){\mathcal H}_\chi$ and
$\,H_\chi\cap(t-\lambda)H(R)=0$, we have
$${\mathcal H}_\chi/(t-\lambda){\mathcal H}_\chi\cong H(R,\pi)/(t-\lambda)H(R,\pi)
\cong H(R)/(t-\lambda)H(R)\cong H_\chi,$$ by the theorem on
isomorphism. Now put $\overline{\mathcal H}_\chi:={\mathcal
H}_\chi/t{\mathcal H}_\chi$ and identify the generators
$\Theta_i=\Theta_{x_i}$ of ${\mathcal H}_\chi$ with their images
in $\overline{\mathcal H}_\chi$. It is immediate from our earlier
remarks that these images satisfy the relations
$[\Theta_{x_i},\Theta_{x_j}]=\Theta_{[x_i,x_j]}$ for all $i,j$. By
the universality property of  the enveloping algebra $U(\z_\chi)$,
there exists an algebra homomorphism $\phi\colon\,U(\z_\chi)
\twoheadrightarrow \overline{\mathcal H}_\chi$ with
$\phi(x_i)=\Theta_i$ for all $i$. Since ${\mathcal H}_\chi$ is a
free $\k[t]$-module, the monomials
$\Theta_1^{a_1}\cdots\Theta_r^{a_r}$ with $(a_1,\ldots,
a_r)\in\Z_+^r$ are linearly independent in $\overline{\mathcal
H}_\chi$. As a consequence, $\phi$ is an isomorphism.
\end{pf}
\subsection{} \label{2.4}
Let ${\bf A}_e$ denote the associative algebra over $\k$ generated
by $z_1,\ldots,z_s, z_{s+1},\ldots z_{2s}$ subject to the
relations $[z_{i+s},z_j]=\delta_{ij}$ and
$[z_i,z_j]=[z_{i+s},z_{j+s}]=0$ where $1\le i,j\le s$. Clearly,
${\bf A}_e \cong {\mathbf A}_s(\k)$, the $s^{\text{th}}$ Weyl
algebra over $\k$. If  $s=0$ then ${\mathbf A}_e=\k$.

Let $i\mapsto i^*$ denote the involution on the set of indices
$\{1,\ldots,s,s+1,\ldots, 2s\}$ such that $i^*=i+s$ for $i\le s$
and $i^*=i-s$ for $i>s$. For $1\le i\le 2s$ define
$z_i^*:=(-1)^{p(i)}\,z_{i^*}$ where
\[
p(i)=\left\{
\begin{array}{rr}
0&  \text{\ if \  $i\le s$},\\
1&  \text{\ if \  $i>s$}.
\end{array}
\right.
\]
Note that $z_i=(-1)^{p(i^*)}\,z_{i^*}^*$ for $1\le i\le 2s$ and
$z_i^*=z_{i+s}\,\,$, $z_{i+s}^*=-z_i$ for $1\le i\le s$.
It is worth remarking that the following relation holds in $U(\g)$:
\begin{equation}\label{GK4}
\sum_{i=1}^{2s}\,z_iz^*_i\,=\,-\sum_{i=1}^{2s}\,z_i^*z_i\,\, \equiv\, s\
\big(\text{mod}\ \, {\mathcal I}_\chi\big).
\end{equation}
As the form $\la\cdot\,,\,\cdot\ra$
is $\z_\chi(0)$-invariant and $\la z_i^*,z_j\ra=\delta_{ij}$ for $1\le i,j\le 2s$,
for all $x\in\z_\chi(0)$ we have
\begin{equation}\label{GK5}
[x,z_k^*]\,=\,\sum_{i=1}^{2s}\,\la z_k^*, [x,z^*_i] \ra \,z_i
\,=\,-\sum_{i=1}^{2s}\,\la [x,z^*_i] , z_k^*\ra \,z_i.
\end{equation}

Each $h\in H_\chi$ is determined by its effect on the canonical
generator $1_\chi\in\widehat{Q}_\chi$. Since the vector
$h(1_\chi)$ can be uniquely expressed as $h(1_\chi)=
\Big(\sum_{\,{\bf i}\in\Z_+^{2s}}\, u_{\bf i}\cdot z^{\bf
i}\Big)\otimes 1_\chi$ with $u_{{\bf i}}\in U(\p_e),$ one obtains
a natural linear injection
\begin{equation}
\widetilde{\mu}\colon\,H_\chi\map\, U(\p_e)\otimes{\mathbf
A}_e^{\mathrm op},\qquad\, \widetilde{\mu}(h)\,=\,{\sum}_{\,{\bf
i}\in\Z_+^{2s}} u_{\bf i}\otimes z^{\bf i}.
\end{equation}

As the form $\la\cdot\,,\,\cdot\ra$ is $C(e)$-invariant, the group $C(e)$ acts
on $\A$ as automorphisms. As $C(e)$ also acts on $U(\p_e)$, it acts as automorphisms
on the algebra $U(\p_e)\otimes\A$, via $g(u\otimes a)=g(u)\otimes g(a)$ with
the obvious choices of $g$, $u$, $a$.
\begin{prop}\label{L2}
The map $\widetilde{\mu}\colon\,H_\chi\into\,
U(\p_e)\otimes{\mathbf A}_e^{\mathrm op}$ is a $C(e)$-equivariant
algebra homomorphism.
\end{prop}
\begin{pf} Let $\mathcal Z$ denote the linear span of all $z^{\bf i}\otimes 1_\chi$ with
${\bf i}\in\Z_+^{2s}$. We identify $\mathcal Z$ with the space of
the left regular representation of ${\mathbf A}_e$ via $z^{\bf
i}\otimes 1_\chi \mapsto z^{\bf i}$. Now $\widehat{\rho}_\chi$
induces a representation of $U(\n_\chi)$ in ${\End}({\mathcal
Z})$, say $\psi_0$. Since $\g(-1)\subset\n_\chi$ and $\g(i)\subset
\text{Ker}\,\chi$ for all $i\le -3$, the definition of
$\Q^{\,\ad\,\n_\chi}$ and induction on $k$ show that
$$\widehat{\rho}_\chi(z_1\cdots z_k)\big(h(1_\chi)\big)\,=\,
{\sum}_{\,{\bf i}\in\Z_+^{2s}}\, u_{\bf i}\cdot
\widehat{\rho}_\chi\big(z^{\bf i}\cdot z_1\cdots z_k\big)(
1_\chi)$$ for all $z_1,\ldots,z_k\in\g(-1)$. Now let $h'$ be
another element in $H_\chi$ and suppose that $h'(1_\chi)=
\Big({\sum}_{\,{\bf i}\in\Z_+^{2s}}\, u'_{\bf i}\cdot z^{\bf
i}\Big)\otimes 1_\chi$  where $u'_{\bf i}\in U(\p_e)$. Then
\begin{eqnarray*}
(h\cdot h')(1_\chi)&=&h'\big(h(1_\chi)\big)\,\,=\,\,\sum_{\,{\bf i}}\,
\widehat{\rho}_\chi(u_{\bf i})\cdot \widehat{\rho}_\chi(z^{\bf i})
\big(h'(1_\chi )\big)\\
&=&\sum_{\,{\bf i}}\,
\sum_{\,{\bf j}}\, u_{\bf i}\cdot
u'_{\bf j}\cdot \widehat{\rho}_\chi\big(
z^{\bf j}\cdot z^{\bf i}\big)(1_\chi).
\end{eqnarray*}
It remains to note that the map $z^{\bf i}\otimes 1_\chi\mapsto
z^{\bf i}$ mentioned above identifies $\psi_0\big(U(\n_\chi)\big)$
with the image of ${\mathbf A}_e$ in its left regular
representation. The $C(e)$-equivariance of $\widetilde{\mu}$ is
immediate from the definitions.
\end{pf}
\begin{rem}
Composing  $\widetilde{\mu}$ with the natural projection
$U(\p_e)\otimes{\mathbf A}_e^{\rm op}\twoheadrightarrow
U(\g(0))\otimes {\mathbf A}_e^{\rm op}$ one obtains an algebra
homomorphism
$$\mu\colon\, H_\chi\,\longrightarrow\,U(\g(0))\otimes{\mathbf A}_e^{\rm op}$$
which will be referred to as
the {\it Miura map}. In the special case where $e$ is even this map has already
appeared
in [\cite{P02}, (7.1)] (note that for $e$ even we have $\A=\k$). It can  be proved that
the map $\mu$
is always injective (this will not be required  in the present note).
\end{rem}

The adjoint action of $\z_\chi(0)$ on $\g$ induces Lie algebra
maps $\z_\chi(0)\,\longrightarrow\, \Der\,({\mathbf
A}_e^{\text{op}})$, $\z_\chi(0)\,\longrightarrow\,
\Der\big(U(\p_e)\otimes {\mathbf A}_e^{\text{op}}\big)$ and
$\z_\chi(0)\,\longrightarrow\, \Der\,(H_\chi)$ (of course, the
same maps can be obtained by differentiating the respective
actions of $C(e)$ on $\A$, $U(\p_e)\otimes \A$ and $H_\chi$). By
abuse of notation, the image of $x\in\z_\chi(0)$ under each of
these  maps will be denoted by $\ad \,x$.

\subsection{} In what follows we will need
explicit formulae for the generators $\Theta_k$ of small Kazhdan degree.
The reader will observe strong similarity between our formulae  and the expressions
for conserved fields of low conformal weight found
by Kac and Wakimoto [\cite{KW}] in the context of vertex algebras and quantum
reduction.

\smallskip

\begin{lemma}\label{L3}
If $v\in\z_\chi(0)$ then it can be assumed that $$\Theta_v(1_\chi)=
\Big(v+\frac{1}{2}\sum_{i=1}^{2s}\,z_i\,[v,z_i^*]\Big)\otimes 1_\chi
\,=\,\Big(v+\frac{1}{2}\sum_{i=1}^{2s}\,[v,z_i^*]\,z_i\Big)\otimes 1_\chi.$$
\end{lemma}
\begin{pf}
It follows from (\ref{GK2}) and  (\ref{GK2''})
that there exit a scalar $\beta$ and a symmetric matrix $A=(\alpha_{ij})$
of order $2s$ such that
$$\Theta_v(1_\chi)=\big(v+
\frac{1}{2}\sum_{i,j=1}^{2s}\,\alpha_{ij}\,z_iz_j+\beta\big)\otimes
1_\chi.$$  Since $A$ is symmetric and $\big[z_k^*,\,v+
\frac{1}{2}\sum_{i,j=1}^{2s}\,\alpha_{ij}\,z_iz_j+
\beta\big]\in {\mathcal I}_\chi$ for all $k$, it must be that
$[v,z_k^*]=\sum_{j=1}^{2s}\,\alpha_{kj}z_j$. Therefore, after a proper adjustment of
$\beta$ we get
$$\Theta_v(1_\chi)=\Big(v+\frac{1}{2}\sum_{i=1}^{2s}\,z_i\,[v,z_i^*]\Big)\otimes
1_\chi.$$ Since  ${\mathcal I}_\chi$ is $(\ad\, v)$-stable, (\ref{GK4}) yields
$\sum_{i=1}^{2s}z_i\,[v,z_i^*]\,\equiv\,
\sum_{i=1}^{2s}\,[v,z_i^*]\,z_i\ \, \big(\text{mod}\ {\mathcal I}_\chi\big)$.
This completes the proof.
\end{pf}
From now on we always assume that the generators
$\Theta_v$ with $v\in\z_\chi(0)$ are chosen in accordance
with Lemma~\ref{L3}. This has the following advantage:
\begin{lemma}\label{L4} The restriction of $\Theta$ to $\z_\chi(0)$
is a Lie algebra homomorphism, i.e.
$$[\Theta_u,\Theta_v]=
\Theta_{[u,v]}\quad\qquad\, \big(\forall\ u,v\in\z_\chi(0)\big).$$
Moreover, the Lie algebra homomorphism  ${\mathrm ad}\circ\Theta\colon
\z_\chi(0)\map{\mathrm Der}\,(H_\chi)$ coincides with the differential of
the rational action
$C(e)\map {\mathrm Aut}(H_\chi)$.
\end{lemma}
\begin{pf}
We are going to use the injective homomorphism $\widetilde{\mu}$
from (\ref{2.4}). Let $x\in\z_\chi(0)$.
Computing in $\A$ and applying (\ref{GK5}) we get
\begin{eqnarray}
\Big[\frac{1}{2}\sum_{i=1}^{2s}\,[x,z_i^*]\,z_i,\,z\Big]&=
&-\frac{1}{2}\sum_{i=1}^{2s}\big(\la [x,z_i^*],z\ra\,z_i+\la z_i,z\ra \,
[x,z_i^*]\big)\nonumber \\
&=&\ \ \frac{1}{2}\big([x,z]+[x,z]\big)\,=\,[x,z]\label{GK6}
\end{eqnarray}
for all $z\in\g(-1)$. Hence $\ad x\,=\,\ad\big(\frac{1}{2}\sum_{i=1}^{2s}\,[x,z_i^*]
z_i\big)$ as derivations of $\A$. Then
$$[\tilde{\mu}(\Theta_x),\,\widetilde{\mu}(h)]\,=
\,\sum\,\big([x,u^{\bf i}]\otimes z^{\bf i}+u^{\bf i}\otimes [x,z^{\bf i}]\big)
\,=\,\widetilde{\mu}\big((\ad\,x)(h)\big)$$ for all $h\in H_\chi$.
As $\widetilde{\mu}$ is injective, it must be that $[\Theta_x,h]\,=\,(\ad\,x)(h)$, i.e.
the adjoint action of $\Theta(\z_\chi(0))$ coincides with
the differential of the action of $C(e)$ on $H_\chi$. Also,
\begin{eqnarray}\label{GK7}
\sum_{i=1}^{2s}\,[v,z_i^*][u,z_i]
&\stackrel{(\ref{GK5})}{=\!=\!=}
&-\sum_{i,j=1}^{2s}\,\la z_j^*,[u,z^*_i]\ra [v,z_i]\, z_j\nonumber\\
&\stackrel{(\ref{GK5})}{=\!=\!=}&-\sum_{i,j=1}^{2s}\,\la z_i^*,[u,z^*_j]\ra [v,z_i]\, z_j
=-\sum_{j=1}^{2s}\,[v,[u,z_j^*]]\,z_j
\end{eqnarray} as elements in $\A$, for all $u,v\in\z_\chi(0)$.
It follows that
\begin{eqnarray*}
\big[\widetilde{\mu}(\Theta_u),\,
\widetilde{\mu}(\Theta_v)\big]
&\stackrel{(\ref{GK6})}{=\!=\!=}&
[u,v]\otimes 1+\frac{1}{2}\otimes \sum_{i=1}^{2s}\,\big([u,[v,z_i^*]]\, z_i+
[v,z_i^*]\,[u,z_i]\big)\\
&\stackrel{(\ref{GK7})}{=\!=\!=}&
[u,v]\otimes 1+
\frac{1}{2}\otimes \sum_{i=1}^{2s}\,[u,[v,z_i^*]]\, z_i
-\frac{1}{2}\otimes\sum_{i=1}^{2s}\,[v,[u,z^*_i]]\,z_i\\
&=&[u,v]\otimes
1+\frac{1}{2}\otimes\sum_{i=1}^{2s}\,[[u,v],z_i^*]\,z_i\,\,=\,\,
\widetilde{\mu}(\Theta_{[u,v]}).
\end{eqnarray*}
But then $[\Theta_u,\Theta_v]\,=\,\Theta_{[u,v]}$ for all
$u,v\in\z_\chi(0)$,  as stated.
\end{pf}
\begin{corollary}\label{C}
Any two-sided ideal of $H_\chi$ is $\sigma$-stable.
\end{corollary}
\begin{pf}
Let $I$ be a two-sided ideal of $H_\chi$. Clearly, $I$ is
invariant under the adjoint action of $\Theta(\z_\chi(0))$. By
Lemma~\ref{L4}, $I$ is then stable under the differential of the
$C(e)$-action on $H_\chi$. Since $C(e)^\circ$ is a connected
reductive group and the action of $C(e)$ on $H_\chi$ is filtration
preserving, Weyl's theorem on complete reducibility shows that all
subspaces $I\cap H_\chi^k$ are $C(e)^\circ$-stable. Since
$Z(G)\subseteq C(e)$ acts trivially on $H_\chi$ and $(\Ad
C(e))^\circ$ coincides with the image of $C(e)^\circ$ in $\Ad
C(e)\cong C(e)/Z(G)$, Lemma~\ref{L} shows that $I$ is
$\sigma$-stable, as claimed.
\end{pf}
Given $n$ elements $x_1,x_2,x_3,\ldots, x_n$ in a Lie algebra we denote by
$[x_1x_2x_3\ldots x_n]$ the commutator $[\ldots[[x_1,x_2],x_{3}],\ldots,  x_n]$.
\begin{lemma}\label{L5}
If  $v\in\z_\chi(1)$ then the generator $\Theta_v\in H_\chi$ has
the following property:
$$\Theta_v(1_\chi)=\Big(v+\sum_{i=1}^{2s}\,[v,z_i^*]\,z_i+\frac{1}{3}
\sum_{i,j=1}^{2s}\,[vz_i^*z_j^*]\,z_jz_i+z_v \Big)\otimes 1_\chi$$
where $z_v\,=\,\frac{1}{3} \sum_{i=1}^{2s}\big(\sum_{k=1}^{2s}
\big\la z_k,[v,[z_k^*,z_i^*]]\big\ra\big)z_i$. Moreover,
$$[\Theta_u,\Theta_v]=\Theta_{[u,v]} \qquad\ \ \,\big(\forall\,
u\in\z_\chi(0)\big).$$
\end{lemma}
\begin{pf} Let $h_v=\sum_{i}\,[v,z_i^*]\,z_i+\frac{1}{3}\sum_{i,j}\,[vz_i^*z_j^*]\,z_jz_i+z_v$,
an element in $U(\g)$.
By anticommutativity and the Jacobi identity, we have that
$$
[z_k^*,[vz_i^*z_j^*]]\,=\,[z_j^*,[vz_i^* z_k^*]]+[v[z_k^*z_j^*]z_i^*]+[v,[z_i^*[z_k^*,z_j^*]]].
$$
Since $(e,[v,x])=0$ for all $x\in\g$ this yields
\begin{eqnarray}\label{GK8}
\big\la z_k^*, [v z_i^* z_j^*]\big\ra\,=\,\big\la z_j^*, [v z_i^* z_k^*]\big\ra -
\big\la z_i^*, [v z_k^* z_j^*]\big\ra + \big\la z_i^*, [v z_j^* z_k^*]\big\ra
\end{eqnarray}\label
where $1\le i,j,k\le 2s$. Computing in $U(\g)$ modulo ${\mathcal I}_\chi$ we now get:
\begin{eqnarray*}
\Big[z_k^*,\,\sum_{i,j}\,[vz_i^*z_j^*]\,z_jz_i\Big]&\equiv & \sum_{ij}\,
\la z_k^*,[vz_i^*z_j^*]\ra z_jz_i+\sum_{i}\,[vz_i^*z_k^*]\,z_i+
\sum_{i}\,[vz_k^*z_i^*]\,z_i\\
&=&\sum_{ij}\Big(\,\big\la z_j^*, [v z_i^* z_k^*]\big\ra -
\big\la z_i^*, [v z_k^* z_j^*]\big\ra + \big\la z_i^*, [v z_j^* z_k^*]\big\ra
\Big)z_jz_i\\
&+&\sum_{i}\Big([vz_i^*z_k^*]\,z_i+[vz_k^*z_i^*]\,z_i\Big)\\
&\equiv&\sum_{i}\Big([vz_i^*z_k^*]\,z_i-z_i\,[vz_k^*z_i^*]+z_i\,[vz_i^*z_k^*]\Big)\\
&\equiv&3\sum_{i}\Big([vz_i^*z_k^*]\,z_i-\frac{1}{3}\big\la z_i,[vz_k^*z_i^*]\big\ra
+\frac{1}{3}\big\la z_i,[vz_i^*z_k^*]\big\ra\Big)\\
&\equiv&3\Big(\sum_{i}\,([vz_i^*z_k^*]\,z_i\Big)-3[z_k^*,z_v].
\end{eqnarray*}
As a consequence,
\begin{eqnarray*}
[z_k^*,\,h_v]&\equiv&[z_k^*,\,v]+\sum_{i}\,[z_k^*,[v,z_i^*]]\,z_i+[v,z_k^*]+
[z_k^*,z_v]\\
&+&\sum_i\,[vz_i^*z_k^*]\,z_i -[z_k^*,z_v]
\equiv 0\ \,\, (\text{mod}\ {\mathcal I}_\chi)
\end{eqnarray*}
for all $k$. It is easy to see that
$[z,\,h_v]
\equiv 0\ \, (\text{mod}\ {\mathcal I}_\chi)$ for all $z\in \n'_\chi$.
By Lemma \ref{L2}, $\Theta_v$ is
a $(-1)$-eigenvector for $\sigma$. In conjunction with (\ref{GK2}), [\cite{P02},
Lemma~4.5],
and the computation above
this shows that $\Theta_v(1_\chi)=h_v\otimes 1_\chi$.

Now let $u$ be any element in $\z_\chi(0)$ and put
$\Theta':=[\Theta_u,\Theta_v]-\Theta_{[u,v]}$. It is immediate
from  (\ref{GK3}) that $\Theta'$ is a polynomial in $\Theta_{x_i}$
with $x_i\in\z_\chi(0)$. Since $\sigma(\Theta')=-\Theta'$ by Lemma
\ref{L2}, this polynomial must be zero. So
$[\Theta_u,\Theta_v]=\Theta_{[u,v]}$ necessarily holds, completing
the proof.
\end{pf}
\section{\bf Associated varieties and Gelfand--Kirillov dimension}
\subsection{}  At present  very little is known  about finite dimensional
representations of the algebras $H_\chi$. In view of
Proposition~\ref{P1} this can be partly explained by the lack of
detailed information on the structure of the centraliser ${\frak
c}_{\g}(e)$. Besides, if $e$ is not even then there is no obvious
reason for $H_\chi=H_{\chi_e}$ to possess such representations. On
the other hand, the evidence collected so far suggests that each
algebra $H_\chi$ has infinitely many isoclasses of finite
dimensional irreducible representations and dimension formulae for
those have roughly the same format as the Weyl dimension formula
for $H_0=U(\g)$; see (6.4). It is therefore natural to ask:
\begin{question}
{\em Is it true that for any nonzero $h\in H_\chi$ there exists a
finite dimensional irreducible representation $\rho$ of $H_\chi$
such that $\rho(h)\ne 0$?}
\end{question}
Let ${\mathcal C}_\chi$ denote the category of all $\g$-modules on
which $x-\chi(x)$ acts locally nilpotently for each $x\in\m_\chi$.
Given a $\g$-module $M$ we denote by $\text{Wh}(M)$ the subspace
of $M$ consisting of all $m\in M$ such that $x.m=\chi(x)m$ for all
$x\in\m_\chi$. Of course, for $M\in{\mathcal C}_\chi$ we have
$\text{Wh}(M)=0$ if and only if $M=0$. Let $H_\chi\mbox{-mod}$
denote the category of all left $H_\chi$-modules. In the Appendix
to [\cite{P02}], Skryabin proved that the functor
\begin{eqnarray} \label{S1}
H_\chi\mbox{-mod}\,\map \,{\mathcal C}_\chi,\qquad\quad
V\longmapsto Q_\chi\otimes_{H_\chi}V,
\end{eqnarray}
is an equivalence of categories. The inverse equivalence is given by
the functor
\begin{eqnarray}\label{S2}
{\mathcal C}_\chi\,\map\,H_\chi\mbox{-mod},\qquad\quad
M\longmapsto \mbox{Wh}(M);
\end{eqnarray}
see also [\cite{GG}, Sect.~6]. Skryabin's result implies that  the
$\g$-module $Q_\chi\otimes_{H_\chi} V$ is simple if and only if so
is the $H_\chi$-module $V$. By the Irreducibility Theorem, the
associated variety of the annihilator in $U(\g)$ of any simple
$\g$-module coincides with the closure of a nilpotent orbit in
$\g^*$; see [\cite{BoBr}, \cite{J2}, \cite{KT}, \cite{V},
\cite{G}]. Our goal in this section is to determine the associated
varieties of the annihilators ${\rm Ann}_{U(\g)}\, M$ for all
$M\in{\mathcal C}_\chi$ with $\dim\,\text{Wh}(M)<\infty$. Such
modules are in 1-1 correspondence with the finite dimensional
representations of $H_\chi$.
\subsection{}
We recall the definition of the Gelfand--Kirillov dimension
of a finitely generated $U(\g)$-module $M$.
Firstly note that there exists a finite
dimensional subspace $M_0\subseteq M$ such that $M\,=\,\bigcup_{n\ge 0}\,
U_nM_0$ where $U_n$ stands for the $n$th component of the
standard filtration of $U(\g)$. It is known that for all $n\gg 0$
the dimension of $U_nM_0$ is a polynomial in $n$. The {\it Gelfand--Kirillov
dimension}
of $M$, denoted $\text{Dim}(M)$, is defined as the degree of this polynomial.
The key point in this definition  is that  $\text{Dim}(M)$
is independent of the choice of $M_0$; see [\cite{Ja}, p.~134] for more detail.

Now let $I$ be a two-sided ideal of the universal enveloping
algebra $U(\g)$. The subspaces $I_{n}:=I\cap U_n$ with $n\in\Z_+$
form an increasing filtration of $I$ satisfying $U_m I_{n}
\subseteq I_{m+n}$ for all $m,n\in\Z_+$. The associated graded
algebra $\gr\,I\hookrightarrow \gr\,U(\g)\cong S(\g)$ is therefore
identified with a homogeneous ideal of the symmetric algebra
$S(\g)$ stable under the adjoint action of $G$. The {\it
associated variety} ${\mathcal VA}(I)$ of the ideal $I$ is defined
as the maximal spectrum  of the affine algebra $S(\g)/\gr\,I$. It
is immediate from the definition that ${\mathcal VA}(I)$ is a
Zariski closed, conical, $G$-invariant  subset of $\text{Max}\,
S(\g)\,=\,\g^*$. For $M$ as above we have
\begin{eqnarray}\label{va}
\dim {\mathcal VA}\big({\mathrm Ann}_{U(\g)}M\big)\,\le\, 2\,\text{Dim}(M);
\end{eqnarray}
see [\cite{Ja}, (10.7) and (17.11)]. The $\g$-module $M$ is called
{\it holonomic} if the equality holds here, that is
$\dim {\mathcal VA}\big({\mathrm Ann}_{U(\g)}M\big)\,=\, 2\,\text{Dim}(M).$
\subsection{} We are now in a position to
state and prove the main result of this section.
\begin{theorem} \label{S3}
Let $M\in{\mathcal C}_\chi$ and $I={\mathrm Ann}_{U(\g)} M$.
Then the following hold:
\begin{enumerate}
\item[(i)] ${\cal O}_\chi\subset {\mathcal VA}(I).$
\smallskip
\item[(ii)] If  $\dim {\mathrm Wh}(M)<\infty$, then ${\mathrm
Dim}(M)=\frac{1}{2}\dim{\cal O}_\chi$ and ${\mathcal VA}(I)$
coincides with the Zariski closure of ${\cal O}_\chi$. In
particular, $M$ is a holonomic $\g$-module.
\end{enumerate}
\end{theorem}
\begin{pf} (1) Let $\scriptstyle{\top}$ denote the anti-involution
of the algebra $U(\g)$ such that $x^{\scriptstyle{\top}}=-x$ and
$(uv)^\top=u^\top v^\top$ for all $x\in\g$ and all $u,v\in U(\g)$.
Let $M^*$ denote the $\g$-module dual to $M$. It is easy to see
that ${\mathrm Ann}_{U(\g)} M^*\,=\,I^\top$.
Since $\scriptstyle{\top}$ preserves the standard filtration of $U(\g)$ and acts
as a scalar operator on  each factor space $U_n/U_{n-1}$ we have
$\gr\,{\mathrm Ann}_{U(\g)} M^*\,=\,\gr\,(I^\top)\,=\,\gr\,I$.  Consequently,
${\mathcal VA}\big({\mathrm Ann}_{U(\g)} M^*\big)\,=\,{\mathcal VA}(I)$.

Pick any nonzero $m\in {\mathrm Wh}(M)$ and view it as a linear function on
$M^*$ via $m(f)=f(m)$ for all $f\in M^*$. Then
$$m(x.f)=(x.f)(m)=-f(x.m)=-\chi(x)f(m)=-\chi(x)m(f)\qquad\quad\, (\forall x\in\m_\chi).$$
This shows that $m$ is a dual $(\m_\chi,-\chi)$-Whittaker vector
of the $\g$-module $M^*$; see [\cite{M}, p.~221]. Thanks to
Matumoto's theorem [\cite{M}] we are now able to deduce that the
associated variety of ${\mathrm Ann}_{U(\g)}M^*$ contains $-\chi$.
As ${\mathcal VA}\big({\mathrm Ann}_{U(\g)}M^*\big)={\mathcal
VA}(I)$ is conical and $G$-stable this yields ${\cal
O}_\chi\subset {\mathcal VA}(I)$ proving (i).

\smallskip

\noindent (2) From now on suppose that $M_0:={\mathrm Wh}(M)$ is
finite dimensional. Let $2d=\dim{\cal O}_\chi$. By the ${\frak
sl}(2)$-theory, $r=\dim \z_\chi=\dim \g(0)+\dim\g(-1)=\dim
\g(0)+2s$. Hence $d=\dim \m_\chi=\dim\p_e-\dim\g(0)-s=m-r+s$. Let
$Y_1,\ldots, Y_d$ be a basis of $\m_\chi$ with $Y_i\in \g(-l_i-2)$
for some $l_i\ge -1$ and choose $X_i\in\g(l_i)$ with $1\le i\le d$
such that $\chi([Y_i,X_j])=\delta_{ij}$. No generality will be
lost by assuming further that $X_i=z_i$ for $i\le s$ and
$X_{s+j}=x_{r+j}$ for $1\le j\le m-r$, where $z_i$ and $x_j$ are
basis vectors introduced in (2.1). For ${\bf
a}=(a_1,...,a_d)\in\Z_+^d \ $ put
$$|{\bf a}|\,=\,\sum_{i=1}^d
\,a_i,\qquad {\mathrm wt}\,{\bf a}\,=\,-\sum_{i=1}^s\,a_i+
\sum_{i=1}^{m-r}\ n_{r+i}a_{s+i},\qquad
X^{\bf a}=X_1^{a_1}\cdots X_d^{a_d}\in U_{|\bf a|}.$$

\smallskip

\noindent (3) Let $\{m_i\}$ be a basis of $M_0$. Since the
$U(\g)$-module $Q_\chi$ is generated by $1_\chi$  it follows from
(\ref{S1}) and (\ref{S2}) that the $U(\g)$-module $M$ is generated
by $M_0$. As explained in [\cite{P02}, p.~53] the vectors $X^{\bf
a}(m_i)$ with ${\bf a}\in\Z_+^d$ and $i\le\dim M_0$ are linearly
independent. Therefore,
$$\dim U_nM_0\,\ge\, (\dim M_0)\cdot\text{Card}\,\{{\bf a}\in\Z_+^d\,|\,\,n\ge |{\bf a}|\}
\,=\,(\dim M_0)\cdot {{n+d}\choose{d}}.$$ For all $n\gg 0$ the  LHS is a polynomial
in $n$ of degree ${\rm Dim}(M)$, while the RHS is a polynomial in $n$ of degree $d$.
This yields ${\mathrm Dim}(M)\ge d$.

\medskip

\noindent
(4) Now put $N=\max\,\{n\in\Z\,|\,\,\g(n)\ne 0\}$ and let  $M_{d,j}$ denote
the subspace of $M$ spanned by all vectors $X^{\bf a}(m_i)$ with $|{\bf a}|\le j$.
We claim that
\begin{eqnarray*}
U_k M_0\,\subseteq \,M_{d, \,(N+2)k}\qquad\quad\  (\forall
k\in\Z_+).
\end{eqnarray*}
For ${\bf k}=(k_1,\ldots, k_r)\in\Z_+^r$ set
$${\mathrm wt}\,{\bf k}=\sum_{i=1}^rn_ik_i,\qquad
x^{\bf k}=x_1^{k_1}\cdots x_r^{k_r}\in U(\g),\qquad \Theta^{\bf
k}= \Theta_1^{k_1}\cdots \Theta_r^{k_r}\in H_\chi.$$ Note that
${\rm wt}\,{\bf k}\ge 0$. Given ${\bf a}\in \Z_+^d$ and ${\bf
b}\in \Z_+^r$ put $|({\bf a};{\bf b})|_e:={\mathrm wt}({\bf
a})+{\mathrm wt}({\bf b})+ 2|{\bf a}|+2|{\bf b}|$. Using the
formula on [\cite{P02}, p.~27] and the isomorphism $M\cong
Q_\chi\otimes_{H_\chi}M_0$ it is easy to observe that
$$X^{\bf a}\, x^{\bf b}(m_i)=\,
\big(X^{\bf a}\,\Theta^{\bf b}+
\sum_{|({\bf i};{\bf j})|_e=|({\bf a};{\bf b})|_e,\, |{\bf i}+|{\bf j}|>|{\bf a}|+|{\bf b}|}
\mu_{{\bf i},{\bf j}} \,X^{\bf i}\Theta^{\bf j}\ \  +
\sum_{|({\bf i};{\bf j})|_e< \,|({\bf a};{\bf b})|_e}
\mu_{{\bf i},{\bf j}} \,X^{\bf i}\Theta^{\bf j}\big)(m_i)
$$
for some $\mu_{{\bf i},{\bf j}}\in \k$. Since the subspace
$U_kM_0$ is spanned by all $X^{\bf a}x^{\bf b}(m_i)$
with $|{\bf a}|+|{\bf b}|\le k$, it is contained in the span of all $X^{\bf i}(m_j)$ such that
\begin{eqnarray*}
{\mathrm wt}\,{\bf i}&\le& |({\bf i};{\bf j})|_e\,\le\,\,
\max_{|{\bf a}|+|{\bf b}|\le k}\,|({\bf a};{\bf b})|_e\,\le\, \,2k+
\max_{|{\bf a}|+|{\bf b}|\le k}\,({\mathrm wt}\,{\bf a}+{\mathrm wt}\,{\bf b})\\
&\le &2k+\max_{|{\bf a}|+|{\bf b}|\le k}\,(N|{\bf a}|+N|{\bf b}|)\,\le\,\,(N+2)k.
\end{eqnarray*}
The claim follows. Since all vectors $X^{\bf i}(m_j)$ are linearly independent,
we derive:
\begin{eqnarray*}
\dim U_kM_0 &\le&\dim M_{d,\,(N+2)k}\,=\,
(\dim M_0)\cdot\text{Card}\,\{{\bf i}\in\Z_+^d\,|\,\,(N+2)k\ge |{\bf i}|\}\\
&=&(\dim M_0)\cdot {{Nk+2k+d}\choose{d}}.
\end{eqnarray*}
Since the RHS is a polynomial in $k$ of degree $d$, we get
${\mathrm Dim}(M)\le d$. In conjunction with part (3) this shows
that ${\mathrm Dim}(M)=d=(\dim {\cal O}_\chi)/2$.

\medskip

\noindent
(5) Since the $H_\chi$-module
${\mathrm Wh}(M)$  is finite dimensional,
Skryabin's equivalence of categories described in (3.1) implies that
the $\g$-module $M$ has a  composition series
$M=M_1\supset M_2\supset\ldots\supset M_l\supset M_{l+1}=0$
such that $M_i/M_{i+1}\in{\mathcal C}_\chi$ and
$\dim {\mathrm Wh}(M_{i}/M_{i+1})<\infty$ for all $i\le l$. Set
$J_i:={\mathrm Ann}_{U(\g)}\,(M_i/M_{i+1})$.
It is immediate from  the discussion in [\cite{J}, (17.7)] that
\begin{eqnarray}\label{sqr}\sqrt{\gr \,I}\,\,=\,
\,\bigcap_{i=1}^l\,\sqrt{\gr \,J_i}.\end{eqnarray}
On the other hand, it follows from (\ref{va}) and parts (1) and  (4) of this proof
that for any $\g$-module $N\in{\mathcal C}_\chi$ with $\dim {\mathrm Wh}(N)<\infty$
one has
$$\dim {\cal O}_\chi\le \dim {\mathcal VA}({\mathrm Ann}_{U(\g)}\, N)
\le 2\,{\mathrm Dim}(N)= \dim{\cal O}_\chi.$$ In conjunction with
the Irreducibility Theorem mentioned in (3.1) this shows that
${\mathcal VA}(M_i/M_{i+1})$ coincides with the Zariski closure
$\overline{\cal O}_\chi$ for all $i\le l$. But then (\ref{sqr})
yields ${\mathcal VA}(I)\,=\,\overline{\cal O}_\chi$, completing
the proof.
\end{pf}
 \subsection{} Recall that a two-sided ideal $I$ of $U(\g)$ is called
{\it completely prime} (respectively, {\it primitive}) if
$U(\g)/I$ is an domain (respectively, if $I$ is the annihilator of
a simple $\g$-module). For $n\in\Z_+$ the set
$Y_n:=\{\psi\in\g^*\,|\,\dim ({\mathrm Ad}^*\,G)\,\psi=n\}$ is
locally closed in the Zariski topology of  $\g^*$. A (locally
closed) subset of $\g^*$ is called a {\it sheet} if it coincides
with an irreducible component of one of the locally closed sets
$Y_n$. It is well-known that each sheet is $G$-invariant and
contains a unique nilpotent coadjoint orbit (such an orbit may lie
in several sheets, however).

\begin{conj}  Let $e$ be  an arbitrary  nilpotent element in $\g$ and  let
$\chi=\chi_e$
be the corresponding linear function on $\g$.

\smallskip
\begin{itemize}
\item[1.] The algebra $H_\chi$ contains an ideal of codimension $1$.

\smallskip

\noindent \item[2.]  The ideals of codimension $1$ in $H_\chi$ are
finite in number if and only if ${\cal O}_\chi$ is a sheet in
$\g^*$.

\smallskip

\noindent \item[3.]  For any ideal $I$ of codimension $1$  in
$H_\chi$ the ideal $\widetilde{I}={\mathrm
Ann}_{U(\g)}\big(Q_\chi\otimes_{H_\chi}\,H_\chi/{I})$ of $U(\g)$
is completely prime.
\end{itemize}
\end{conj}

\smallskip

\noindent
Our last conjecture provides
a hypothetical converse to Theorem~\ref{S3}(ii).
It indicates that each category
${\mathcal C}_\chi$ is potentially very important
for the theory of primitive ideals.
\begin{conj}  Let $\chi$ be as above and let $I$ be a primitive ideal
of $U(\g)$ whose associated variety equals $\overline{\cal
O}_\chi$. Then there exists a simple $\g$-module $M\in{\mathcal
C}_\chi$ with $\dim\,{\mathrm Wh}(M)<\infty$ such that $I={\mathrm
Ann}_{U(\g)}\,M.$
\end{conj}

\smallskip

\section{\bf Minimal nilpotent algebras: a quadratic relation}
\subsection{}
Let $\h$ be a Cartan subalgebra of $\g$, and let $\Phi$ be the
root system of $\g$ relative to $\h$. Let
$\{e_\alpha\,|\,\alpha\in\Phi\}\cup
\{h_\alpha\,|\,\alpha\in\Phi\}$ be a Chevalley system in $\g$ with
each triple $(e_\alpha, h_\alpha, e_{-\alpha})$ being an ${\frak
sl}_2$-triple in $\g$. Let $\Pi=\{\alpha_1,\ldots,\alpha_\ell\}$
be a basis of simple roots in $\Phi$ with the elements in $\Pi$
numbered as in [\cite{Bou}], and let
$\{\varpi_1,\ldots,\varpi_\ell\}$ be the corresponding system of
fundamental weights in $\h^*$. Let $\Phi^+$ and $\Phi^-$ be the
positive and the negative system of $\Phi$ relative to $\Pi$,
respectively, and let $P$ denote the lattice of integral weights
in $\h^*$. As usual, given $\lambda,\mu\in P$ we write $\lambda\ge
\mu$ if and only if $\lambda-\mu$ is a sum of positive roots. Let
$P^+=\{\sum_i a_i\varpi_i\,|\,\,i\in\Z\}$, the set of dominant
weights, and $\rho=\varpi_1+\cdots
+\varpi_\ell=\frac{1}{2}\sum_{\alpha\in \Phi^+}\alpha$, the
half-sum of positive roots.  Let $W$ be the Weyl group of the root
system $\Phi$; it is generated by reflections $s_\alpha$ with
$\alpha\in\Phi$. The {\it dot action} of $W$ on $\h^*$ is defined
by setting $w\!\centerdot\! \lambda=w(\lambda+\rho)-\rho$ for all
$w\in W$ and $\lambda\in\h^*$.

If $\g$ is not of type $\rm A$ or $\rm C$, there is a unique long
root in $\Pi$ linked with the lowest root
$-\widetilde{\alpha}\in\Phi^-$ on the extended Dynkin diagram of
$\g$; we call it $\beta$. For $\g$ of type $A_n$ and $C_n$ put
$\beta=\alpha_n$. In this paper, we will be mostly concerned with
the ${\frak sl}_2$-triple $(e,h,f)=(e_\beta,h_\beta,e_{-\beta})$.
Recall that the invariant form $(\,\cdot\,,\,\cdot\,)$ on $\g$ has
the property that $(e,f)=1$. This entails $(h,h)=2$. It is
well-known the restriction of $(\,\cdot\,,\,\cdot\,)$ to $\h$ is
nondegenerate and induces a $W$-invariant scalar product on the
$\mathbb Q$-span of $P$ in $\h^*$. More precisely, for all
$\lambda,\mu\in\h^*$ we have $(\lambda,\mu)=(t_\lambda, t_\mu)$
where $t_\lambda,t_\mu\in\h$ are such that
$\lambda=(t_\lambda,\,\cdot\,)$ and $\mu=(t_\mu,\,\cdot\,)$. Put
$\la \lambda,\alpha\ra=2(\lambda,\alpha)/(\alpha,\alpha)$ for all
$\lambda\in \h^*$ and $\alpha\in \Phi$. Since
$(\,\cdot\,,\,\cdot\,)$ is a multiple of the Killing form of $\g$,
there is a constant $c\in\k^\times$ such that
$\beta(x)=c(h_\beta,x)$ for all $x\in\h$. The equality
$\beta(h_\beta)=2=(h_\beta,h_\beta)$ now shows that $c=1$ and
$t_\beta=h_\beta$. Hence $(\gamma,\gamma)=2$ for all long roots
$\gamma\in\Phi$.

From now on we assume that $\chi=(e,\,\cdot\,)$ where $e=e_\beta$.
It is well-known that the adjoint action of $h=h_\beta$ on $\g$
gives rise to a short $\Z$-grading
$$\g\,=\,\g(-2)\oplus\g(-1)\oplus\g(0)\oplus\g(1)\oplus\g(2)$$
with $\g(i)=\{x\in\g\,|\,\,[h,x]=ix\}$ for all $i\in\Z$ and with
$\g(1)\oplus \g(2)\cong \,\g(-1)\oplus \g(-2)$ isomorphic to a
Heisenberg Lie algebra. We also have that $\g(2)=\k e,\,$
$\g(-2)=\k f$, and $\g(0)={\frak c}_\g(h)$. The ${\frak
sl}_2$-theory implies that
$${\frak c}_\g(e)=\g(0)^\sharp\oplus \g(1)\oplus\g(2)$$ where
$\g(0)^\sharp=\{x\in\g(0)\,\vert \,\,[x,e]=0\}$. More importantly
for our later deliberations, $\g(0)^\sharp$ is the orthogonal
complement to $\k h$ in $\g(0)$ and hence coincides with image of
the Lie algebra endomorphism
$$\sharp\,\colon\, \g(0)\longrightarrow\, \g(0),\quad\ x\mapsto
\,x-\frac{1}{2}(x,h)\,h.$$  In particular, $\g(0)^\sharp$ is an
ideal of codimension $1$ in the Levi subalgebra $\g(0)$. It is
well-known that outside type $\rm A$ the centre of $\g(0)$
coincides with $\k h$ and $\g(1)$ is an irreducible $\ad
\g(0)^\sharp$-module. As a consequence, if $\g$ is not of type
$\rm A$, then $\g(0)^\sharp=[\g(0),\g(0)]$ is the only ideal if
codimension $1$ in $\g(0)$ and ${\frak c}_\g(e)=[{\frak
c}_\g(e),{\frak c}_\g(e)]$ is a perfect Lie algebra.

Note that in the present case $\g(0)^\sharp=\Lie C(e).$ Put
$\h_e:=\,\h\cap{\frak c}_\g(e)$, a Cartan subalgebra in
$\g(0)^\sharp$. It can be assumed without loss of generality that
$\h_e=\Lie T_e$ where $T_e$ is as in Section~2. We can  choose
$z_1,\ldots,z_s, z_{s+1},\ldots, z_{2s}$ to be root vectors for
$\h$. Moreover, we can (and will) assume that the $z_i$'s with
$1\le i\le s$ are root vectors in $\g(-1)$ corresponding to {\it
negative} roots $\gamma_i\in\Phi^-$. Then each $z_i^*$ with $1\le
i\le s$ is a root vector in $\g(-1)$ corresponding to
$\gamma_i^*:=-\beta-\gamma_i\in\Phi^+$.
\subsection{} Set $H:=H_\chi$ and identify $H$ with $\tilde{\mu}(H).$
Given $a\in\z_\chi(0)$ and $w\in\z_\chi(1)$ we define the
following elements of ${\mathbf A}_e$:
$$\psi_a\,:=\,\frac{1}{2}\sum_{i=1}^{2s}\,[a,z_i^*]z_i,\qquad\
\varphi_w\,:=\,\frac{1}{3}\sum_{i,j=1}^{2s}\, [w
z_i^*z_j^*]z_jz_i+z_w.$$ Recall that
$\sum_{i=1}^{2s}[w,z_i^*]z_i=-\sum_{i=1}^{2s}[w,z_i]z_i^*$. The
computation used in the proof of Lemma~2.5 shows that
\begin{eqnarray}\label{MN1}
[z_k^*,\varphi_w]\,=\,\sum_{i=1}^{2s}\,[wz_i^*z^*_k]z_i.\end{eqnarray}

Notice that in $U(\p_e)\otimes {\mathbf A}_e^{\rm op}$ we have
$\sum_{i=1}^{2s}[w,z_i^*]\otimes
z_i=-\sum_{i=1}^{2s}[w,z_i]\otimes z_i^*$, and
$$[a\otimes f,b\otimes g]=[a,b]\otimes gf-ba\otimes [f,g]\qquad
\big(\forall\, a,b\in U(\p_e), \quad \forall\, f,g\in {\mathbf
A}_e\big).$$ Keeping this in mind it is straightforward to see
that for all $u,v\in\z_\chi(1)$,
\begin{eqnarray}\label{MN2}
\big[\Theta_u,\Theta_v\big]&=&\big[u\otimes
1+\sum_{i=1}^{2s}\,[u,z_i^*]\otimes z_i+ 1\otimes
\varphi_u,v\otimes 1+\sum_{i=1}^{2s}\,[v,z_i^*]\otimes z_i +
1\otimes \varphi_v\big]\nonumber\\ &=& [u,v]\otimes
1+\sum_{i=1}^{2s}\,[[u,v],z_i^*]\otimes
z_i+\sum_{i,j=1}^{2s}\,[u,z_i]\otimes [z_i^*,\varphi_v]\nonumber\\
&-&\sum_{i,j=1}^{2s}\,[v,z_i]\otimes
[z_i^*,\varphi_u]+\sum_{i,j=1}^{2s}\,\big[[u,z_i^*],[v,z_j^*]\big]\otimes
z_jz_i\nonumber\\
 &+&\sum_{i=1}^{2s}\, [v,z_i^*][u,z_i]\otimes 1
-1\otimes [\varphi_u,\varphi_v].\end{eqnarray} In view of
(\ref{MN1}) and (\ref{GK8}) we have that
\begin{eqnarray}\label{MN3}
&\displaystyle{\sum_{i=1}^{2s}}&\!\![x,z_i]\otimes
[z_i^*,\varphi_y]\,=\,\sum_{i,j=1}^{2s}\,[x,z_i]\otimes[y
z_j^*z_i^*]z_j\, =\,\sum_{i,j,\,k}\,[x,z_i]\otimes \langle
z_k^*,[yz_j^*z_i^*] \rangle z_k z_j\nonumber\\
&=&\sum_{i,j,\,k}\,[x,\langle z_i^*,[yz_j^*z_k^*] \rangle
z_i]\otimes  z_k z_j+\sum_{i,j,\,k}\,[x,z_i]\otimes z_k \langle
z_j^*,[y,[z_i^*,z_k^*]] \rangle z_j\nonumber\\
&=&\sum_{i,j=1}^{2s}\,[x,[yz_j^*z_i^*]]\otimes
z_iz_j+\sum_{i=1}^{2s}\,[x,z_i^*]\otimes z_i[y,f]
\end{eqnarray}
for all $x,y\in\z_\chi(1)$. Combining this with (\ref{MN2}) and
taking into account that
$$\sum_{i,j=1}^{2s}\big[[u,z_i^*],[v,z_j^*]\big]\otimes z_iz_j=
\sum_{i,j=1}^{2s}\big[[u,z_i^*],[v,z_j^*]\big]\otimes
z_jz_i-\sum_{i=1}^{2s}\,\big[[u,z_i],[v,z_j^*]\big],$$ we now
obtain
\begin{eqnarray}\label{MN4}
\big[\Theta_u,\Theta_v\big]&=& [u,v]\otimes
1+\sum_{i=1}^{2s}\,[[u,v],z_i^*]\otimes
z_i-\sum_{i,j=1}^{2s}\,[v,[uz_j^*z_i^*]]\otimes
z_iz_j\nonumber\\
&-&\sum_{i=1}^{2s}\,[v,z_i^*]\otimes z_i[u,f]
+\sum_{i,j=1}^{2s}\,[u,[vz_j^*z_i^*]]\otimes
z_iz_j+\sum_{i=1}^{2s}\,[u,z_i^*]\otimes z_i[v,f]\nonumber\\
&+&\sum_{i,j=1}^{2s}\,\big[[u,z_i^*],[v,z_j^*]\big]\otimes z_jz_i
 +\sum_{i=1}^{2s}\, [v,z_i^*][u,z_i]\otimes 1
-1\otimes [\varphi_u,\varphi_v]\nonumber\\
 &=& (f,[u,v])\Big(e\otimes
1+\sum_{i=1}^{2s}\,[[e,z_i^*]\otimes z_i+
\sum_{i,j=1}^{2s}\,[ez_i^*z_j^*]\otimes
z_jz_i\Big)\nonumber\\
&-&\sum_{i,j=1}^{2s}\big[[u,z_i^*],[v,z_j^*]\big]\otimes z_jz_i
+\sum_{i=1}^{2s}\big([u,z_i^*]\otimes z_i[v,f]-[v,z_i^*]\otimes
z_i[u,f]\big)\nonumber\\
&+&\sum_{i=1}^{2s}[u,z_i][v,z_i^*]\otimes 1 -1\otimes
[\varphi_u,\varphi_v].\end{eqnarray}

Next observe that
$[a,z_i^*]^\sharp=[a,z_i]-\frac{1}{2}(h,[a,z_i^*])h,$ and
$$\sum_{i=1}^{2s}\,(h,[a,z_i^*])z_i\,=\,\sum_{i=1}^{2s}\,(e,[f,[a,z_i^*]])z_i\,=\,\sum_{i=1}^{2s}\,
\langle z_i^*,[a,f] \rangle z_i=[a,f]$$ for all $a\in\z_\chi(1)$.
Since $\sum_{i=1}^{2s}z_i^*z_i=s$ (as elements in ${\bf A}_e$), it
follows from (\ref{MN3}) that
\begin{eqnarray}
&\displaystyle{\sum_{i,j=1}^{2s}}&\!\![x,z_i]^\sharp\otimes
[[y,z_i^*]^\sharp,z_j^*]z_j\,=\,\sum_{i,j=1}^{2s}\,[x,z_i]^\sharp\otimes
[yz_i^*z_j^*]z_j
+\frac{s}{2}\sum_{i=1}^{2s}(h,[y,z_i^*])[x,z_i]^\sharp\otimes
1\nonumber\\
&=&\sum_{i,j=1}^{2s}\,[x,z_i]\otimes
[yz_i^*z_j^*]z_j-\frac{1}{2}\sum_{i,j=1}^{2s}(h,[x,z_i^*])h\otimes
[yz_i^*z_j^*]z_j +\frac{s}{2}[x,[y,f]]^\sharp\otimes 1\nonumber\\
&=&\frac{s}{2}[x,[y,f]]^\sharp\otimes
1+\sum_{i,j=1}^{2s}\,[x,z_i]\otimes
[yz_j^*z_i^*]z_i+\sum_{i=1}^{2s}\,[x,z_i^*]\otimes
[y,f]z_i\nonumber\\
&+&\frac{1}{2}\sum_{i=1}^{2s}\,h\otimes [[y,[x,f]],z_j^*]z_j\,=\,
\frac{s}{2}[x,[y,f]]^\sharp\otimes
1+\frac{1}{2}\sum_{i=1}^{2s}\,h\otimes
[[y,[x,f]],z_i^*]z_i\nonumber\\
&+&\sum_{i,j=1}^{2s}\,[x,[yz_j^*z_i^*]]\otimes
z_iz_j+\sum_{i=1}^{2s}\,[x,z_i^*]\otimes
z_i[y,f]+\sum_{i=1}^{2s}\,[x,z_i^*]\otimes [y,f]z_i\nonumber
\end{eqnarray}
for all $x,y\in\z_\chi(1)$. But then
\begin{eqnarray}\label{MN5}
&\displaystyle{\sum_{i,j=1}^{2s}}&\!\![x,z_i]^\sharp\otimes
[[y,z_i^*]^\sharp,z_j^*]z_j+
\sum_{i,j=1}^{2s}\,[y,z_i^*]^\sharp\otimes
[[x,z_i]^\sharp,z_j^*]z_j\nonumber\\
&=&\sum_{i,j=1}^{2s}\,[x,z_i]^\sharp\otimes
[[y,z_i^*]^\sharp,z_j^*]z_j-
\sum_{i,j=1}^{2s}\,[y,z_i]^\sharp\otimes
[[x,z_i^*]^\sharp,z_j^*]z_j\nonumber\\
&=&(f,[x,y])\Big(\frac{s}{2}h^\sharp\otimes
1-\frac{1}{2}h\otimes\sum_{i=1}^{2s}\,[h,z_i]^*z_i+\sum_{i,j=1}^{2s}[ez_i^*z_j^*]\otimes
z_jz_i\Big)\nonumber\\
&-&2\sum_{i,j=1}^{2s}\,\big[[x,z_i^*],[y,z_j^*]\big]\otimes
z_jz_i-\sum_{i=1}^{2s}\,\big[[x,z_i],[y,z_i^*]\big]\otimes
1\nonumber\\
&+&2\sum_{i=1}^{2s}\big([x,z_i^*]\otimes z_i[y,f]-[y,z_i^*]\otimes
z_i[x,f]\big)-[[x,y],f]\otimes 1.
\end{eqnarray}
We used the fact that $$\sum_{i=1}^{2s}\,[x,z_i^*]\otimes
[[y,f],z_i]\,=\,\sum_{i=1}^{2s}\,[x,z_i]\otimes\langle
z_i^*,[y,f]\rangle=-[x,[y,f]]\otimes 1.$$

To ease notation, set
$$A(u,v):=[\Theta_u,\Theta_v]-\frac{1}{2}\sum_{i=1}^{2s}\,
\big(\Theta_{[u,z_i]^\sharp}\,\Theta_{[v,z_i^*]^\sharp}+
\Theta_{[v,z_i^*]^\sharp}\,\Theta_{[u,z_i]^\sharp}\big).$$ Note
that
\begin{eqnarray*}
&\displaystyle{\sum_{i=1}^{2s}}&\big(\Theta_{[u,z_i]^\sharp}\,\Theta_{[v,z_i^*]^\sharp}+
\Theta_{[v,z_i^*]^\sharp}\,\Theta_{[u,z_i]^\sharp}\big)\,=\,\sum_{i=1}^{2s}\big(
[u,z_i]^\sharp [v,z_i^*]^\sharp+
[v,z_i^*]^\sharp [u,z_i]^\sharp\big)\otimes 1\\
&+&\sum_{i,j=1}^{2s}\,[u,z_i]^\sharp\otimes[[v,z_i^*]^\sharp,z_j^*]z_jz_i
+\sum_{i,j=1}^{2s}\,[v,z_i^*]^\sharp\otimes[[u,z_i]^\sharp,z_j^*]z_jz_i+1\otimes\psi(u,v),
\end{eqnarray*}
where
$\psi(u,v)=\frac{1}{4}\sum_{i=1}^{2s}\big(\psi_{[v,z_i^*]^\sharp}\,\psi_{[u,z_i]^\sharp}+\psi_{[u,z_i]^\sharp}
\,\psi_{[v,z_i^*]^\sharp}\big)\in{\mathbf A}_e$. Since
$$\sum_{i=1}^{2s}(x,z_i^*)z_i\,=\,-\sum_{i=1}^{2s}\,(x,z_i)z_i^*=\,[x,f]\qquad\
\big(\forall \,x\in\z_\chi(1)\big),$$ we have that
\begin{eqnarray*}
&\displaystyle{\sum_{i=1}^{2s}}&\!\!\big( [u,z_i]^\sharp
[v,z_i^*]^\sharp+ [v,z_i^*]^\sharp
[u,z_i]^\sharp\big)\,=\,\sum_{i=1}^{2s}\,\big( [u,z_i] [v,z_i^*]+
[v,z_i^*][u,z_i]\big)\\
&-&\sum_{i=1}^{2s}\,\big( (u,z_i)h [v,z_i^*]+
(v,z_i^*)h[u,z_i]-\frac{1}{4}(u,z_i)(v,z_i^*)h^2-\frac{1}{4}(v,z_i^*)(u,z_i)h^2\big)\\
&=&-\frac{1}{2}(e,[u,v])h^2+2\sum_{i=1}^{2s}\,\big[u,z_i]
[v,z_i^*]-\sum_{i=1}^{2s}\,\big[ [u,z_i], [v,z_i^*]\big].
\end{eqnarray*}
Since $h^\sharp=0$ and
$h\otimes\sum_{i=1}^{2s}\,[h,z_i^*]z_i=-sh\otimes 1$, we now
combine the above with (\ref{MN4}) and (\ref{MN5}) to deduce that
\begin{eqnarray}\label{MN6}
A(u,v)&=&\frac{1}{2}(e,[u,v])\Big(\frac{4e+h^2-(s+2)h}{2}\otimes
1+2\sum_{i=1}^{2s}\,[e,z_i^*] \otimes z_i\nonumber\\
&+&\sum_{i,j=1}^{2s}\,[ez_i^*z_j^*]\otimes z_jz_i\Big)+1\otimes
\big([\varphi_v,\varphi_u]-\frac{1}{2}\psi(u,v)\big).
\end{eqnarray}

\subsection{} Recall that $\z_\chi(0)=\,\g(0)^\sharp$ is an ideal of codimension
$1$ in the Levi subalgebra $\g(0)={\mathfrak c}_{\g}(h)$ of $\g$.
Let $\{a_i\}$ and $\{b_i\}$ be dual bases of $\z_\chi(0)$ with
respect to the restriction of $(\,\cdot\,,\,\cdot\,)$ to
$\z_\chi(0)$, and let $C_0=\sum_{i}\,a_ib_i$ be the corresponding
Casimir element of $U(\z_\chi(0))$. Set $\Theta_{\rm
Cas}\,:=\,\sum_{i}\,\Theta_{a_i}\Theta_{b_i}$, an element of $H$.
Although $\Theta_{\rm Cas}$ is not central in $H$, Lemma~\ref{L4}
shows that it commutes with all operators $\Theta_x$ for
$x\in\z_\chi(0)$. Since the skew-symmetric form
$\langle\,\cdot\,,\,\cdot\,\rangle$ is invariant under
$\z_\chi(0)$ and $\k h$ is orthogonal to $\g(0)^\sharp$ with
respect to $(\,\cdot\,,\,\cdot\,)$, we have
\begin{eqnarray*}
\sum_{i,j=1}^{2s}\,[ez_i^*z_j^*]^\sharp\otimes
z_jz_i&=&\sum_{i,j,k}\,(b_k,[ez_i^*z_j^*]^\sharp]) a_k\otimes
z_jz_i\,=\,\sum_{i,j,k}\,(b_k,[ez_i^*z_j^*]]) a_k\otimes z_jz_i\\
&=&\sum_{i,j,k}\,([b_kz_j^*z_i^*],e)a_k\otimes
z_jz_i\,=\,-\sum_{i,j,k}\,\langle z_i^*, [b_kz_j^*]\rangle
a_k\otimes z_jz_i\\
&=&-\sum_{i,j,k}\,\langle z_j^*, [b_kz_i^*]\rangle a_k\otimes
z_jz_i=-\sum_{k,i}\,a_k\otimes [b_k,z_i^*]z_i.
\end{eqnarray*}
Interchanging the r{\^o}les of $\{a_i\}$ and $\{b_i\}$ we now
obtain
$$\sum_{i,j=1}^{2s}\,[ez_i^*z_j^*]^\sharp\otimes
z_jz_i\,=\,-\sum_{k,i}\,b_k\otimes [a_k,z_i^*]z_i.$$ Next observe
that
\begin{eqnarray*}
\sum_{i,j=1}^{2s}\,[ez_i^*z_j^*]^\sharp\otimes
z_jz_i&=&\sum_{i,j=1}^{2s}\,[ez_i^*z_j^*]\otimes z_jz_i-
\frac{1}{2}\sum_{i,j=1}^{2s}\,(h,[ez_i^*z_j^*])h\otimes
z_jz_i\\
&=&\sum_{i,j=1}^{2s}\,[ez_i^*z_j^*]\otimes
z_jz_i+\frac{1}{2}\sum_{i,j=1}^{2s}\,([[h,z_j],z_i^*],e)h\otimes
z_j^*z_i\\
&=&\sum_{i,j=1}^{2s}\,[ez_i^*z_j^*]\otimes
z_jz_i+\frac{sh}{2}\otimes 1.
\end{eqnarray*}
It follows that
$$
\sum_{i,j=1}^{2s}\,[ez_i^*z_j^*]\otimes
z_jz_i\,=\,-\frac{sh}{2}\otimes
1-\frac{1}{2}\Big(\sum_{i,j}\,a_i\otimes
[b_i,z_j^*]z_j+\sum_{i,j}\,b_i\otimes [a_i,z_j^*]z_j\Big).
$$ As a result,
\begin{eqnarray}\label{MN7}
\Theta_{\rm Cas}&=&\sum_{i}\,a_ib_i\otimes
1+\frac{1}{2}\sum_{i,j}\,\big(a_i\otimes[b_i,z_j^*]z_j+b_i\otimes
[a_i,z_j^*]z_j\big)+1\otimes c'_0\nonumber\\
&=&\Big(-\frac{sh}{2}+\sum_{i}a_ib_i\Big)\otimes
1-\sum_{i,j=1}^{2s}\,[ez_i^*z_j^*]\otimes z_jz_i+1\otimes c'_0,
\end{eqnarray}
where
$c'_0\,=\,\frac{1}{4}\sum_{i,j,k}\,[b_i,z_j^*]z_j[a_i,z_k^*]z_k\in{\mathbf
A}_e$.
\subsection{} Let $C$ denote the Casimir element of $U(\g)$
corresponding to the invariant form $(\,\cdot\,,\,\cdot\,)$.
Clearly, $C$ induces a $\g$-endomorphism of $\widehat{Q}_\chi$,
and hence can be viewed as a central element of $H$. To determine
$\widetilde{\mu}(C)$ we first observe that the (ordered) sets
$$\{e,h,f\}\cup\{a_i\}\cup\{[e,z_i^*]\,|\,1\le i\le
2s\}\cup\{z_i\,|\,1\le i\le 2s\}$$ and
$$\{f,\frac{h}{2},e\}\cup\{b_i\}\cup\{z_i\,|\,1\le i\le
2s\}\cup\{[e,z_i^*]\,|\,1\le i\le 2s\} $$ are dual bases of $\g$
with respect to $(\,\cdot\,,\,\cdot\,)$. Since
$$
\sum_{i=1}^{2s}\,[z_i,[e,z_i^*]]\,=\,\sum_{i=1}^{2s}\big([[z_i,e],z_i^*]+[e,[z_i,z_i^*]]\big)
\,=\,-\sum_{i=1}^{2s}\big([[z_i^*,e],z_i]+[e,[z_i^*,z_i]]\big),
$$
we have $\sum_{i=1}^{2s}\,[z_i,[e,z_i^*]]\,=\,-sh.$ As
$\widehat{\rho}_\chi(f)1_\chi=(e,f)1_\chi=1_\chi$, it follows that
\begin{eqnarray*}
C(1_\chi)&=&\Big(2e+\frac{h^2}{2}-h+\sum_{i}a_ib_i+2\sum_{i=1}^{2s}\,[e,z_i^*]z_i+\sum_{i=1}^{2s}\,
[z_i,[e,z_i^*]]\Big)\otimes 1_\chi\\
&=&\Big(2e+\frac{h^2}{2}-(s+1)h+\sum_{i}a_ib_i+2\sum_{i=1}^{2s}\,[e,z_i^*]z_i\Big)\otimes
1_\chi.
\end{eqnarray*}
As a consequence,
\begin{eqnarray}\label{MN8}
\widetilde{\mu}(C)\,=\,\Big(2e+\frac{h^2}{2}-(s+1)h+\sum_{i}a_ib_i\Big)\!\otimes
1+2\sum_{i=1}^{2s}\,[e,z_i^*]\otimes z_i.
\end{eqnarray}
In view of (\ref{MN7}) and the convention of 4.1 this yields
\begin{eqnarray}\label{MN9}
C-\Theta_{\rm Cas}&=&\Big(2e+\frac{h^2-(s+2)h}{2}\Big)\otimes
1\nonumber\\
&+&2\sum_{i=1}^{2s}\,[e,z_i^*]\otimes z_i
+\sum_{i,j=1}^{2s}\,[ez_i^*z_j^*]\otimes z_jz_i-1\otimes c'_0.
\end{eqnarray}
\begin{prop}\label{P2}
Let $u,v\in\z_\chi(1)$. Then the following relation holds in $H$:
$$[\Theta_u,\Theta_v]\,=\,\frac{1}{2}(f,[u,v])\big(C-\Theta_{\rm Cas}-c_0\big)+
\frac{1}{2}\sum_{i=1}^{2s}\,\big(\Theta_{[u,z_i]^\sharp}\,\Theta_{[v,z_i^*]^\sharp}+
\Theta_{[v,z_i^*]^\sharp}\,\Theta_{[u,z_i]^\sharp}\big),
$$ where $c_0\in \k$ is a constant depending on $\g$.
\end{prop}
\begin{pf} Set
$B(u,v):=A(u,v)-\frac{1}{2}(f,[u,v])\big(C-\Theta_{\rm Cas}\big)$,
an element in $H$. From (\ref{MN6}) and (\ref{MN9}) we deduce that
$\widetilde{\mu}\big(B(u,v)\big)\,=\,1\otimes b(u,v)$ for some
$b(u,v)\in {\mathbf A}_e$. In conjunction with [\cite{P02},
Lemma~4.5] this shows that $b(u,v)\in \k$ for all
$u,v\in\z_\chi(1)$. In other words,
$b\colon\,\z_\chi(1)\times\z_\chi(1)\,\rightarrow\,\k$,
$\,(u,v)\mapsto\, b(u,v),$ is a bilinear form on $\z_\chi(1)$. In
view of Lemma~\ref{L5}, the $\z_\chi(0)$-invariance of
$\langle\,\cdot\,,\,\cdot\,\rangle$, and the Jacobi identity, this
form is invariant under $\z_\chi(0)$. On the other hand, it is
well-known (and easily seen) that if $\g$ is not of type $\mathrm
A$, then $\z_\chi(1)$ is an irreducible $\z_\chi(0)$-module, and
if $\g$ is of type $\mathrm A$ and $\z_\chi(1)\ne \,0$, then
$\z_\chi(1)\cong M\oplus M^*$ where $M$ is an irreducible
$\z_\chi(0)$-module such that $M\not\cong M^*$. Therefore, in all
cases $b=c_0(f,[\,\cdot\,,\,\cdot\,])$ for some $c_0\in\k$. This
completes the proof.
\end{pf}

Let $\{x_1,\ldots, x_n\}$ and $\{u_1,\ldots, u_{2s}\}$ be  bases
of $\z_\chi(0)$ and $\z_\chi(1)$, respectively, and let $H^+$
denote the $\k$-span of
$$\big\{\Theta_{x_1}^{i_1}\cdots
\Theta_{x_n}^{i_n}\cdot\Theta_{u_1}^{j_1}\cdots\Theta_{u_{2s}}^{j_{2s}}
\cdot(C-c_0)^l\,|\,\,\, \textstyle{\sum}\,i_k+\textstyle{\sum}\,
{j_k}+l\ge 1\big\},$$ a subspace of codimension $1$ in $H$; see
[\cite{P02}, Theorem~4.6(ii)].
\begin{corollary}\label{C1}
The subspace $H^+$ is a two-sided ideal of the minimal nilpotent
algebra $H$. If $\g$ is not of type $\mathrm A$, then $H^+$ is the
only ideal of codimension $1$ in $H$.
\end{corollary}
\begin{pf}
We need to show that $h\cdot h'\in H^+$ for all $h,h'\in H^+$.
Since $C-c_0\in Z(H)$, it suffices to show that $\Theta_x\cdot
H^+\subset H^+$ for all $x\in\z_\chi(0)\cup\z_\chi(1)$. Using
Lemma~\ref{L4} it is easy to observe that the span of all
$\Theta_{x_1}^{i_1}\cdots \Theta_{x_n}^{i_n}$ with $\sum_{k=1}^n
i_k\ge 1$ is stable under the left multiplications by $\Theta_x$
with $x\in \z_\chi(0)$. Thus we may assume that $x\in\z_\chi(1)$.
Lemma~\ref{L5} (in conjunction with [\cite{P02}, Theorem~4.6(ii)])
now provides a further reduction: it suffices to show that
$\Theta_{x}\cdot\big(\Theta_{u_{j_1}}\cdots\Theta_{u_{j_N}}\big)\in
H^+$ for all $x\in\z_\chi(1)$ and all $j_1,\ldots,
j_N\in\{1,\ldots,2s\}$. This follows from Proposition~\ref{P2} and
Lemma~\ref{L5} by induction on $N$.

Suppose $\g$ is not of type $\mathrm A$. Then only one node on the
extended Dynkin diagram of $\g$ is linked with the node
corresponding to $-\tilde{\alpha}$. From this it is immediate that
the derived subalgebra of ${\mathfrak c}_{\g}(h_{\tilde{\alpha}})$
is semisimple and has codimension $1$ in ${\mathfrak
c}_{\g}(h_{\tilde{\alpha}})$. Since the roots $\beta$ and
$-\tilde{\alpha}$ lie in the same $W$-orbit, the subalgebras
${\mathfrak c}_{\g}(h_{\tilde{\alpha}})$ and $\g(0)={\mathfrak
c}_{\g}(h)$ are conjugate under $\Ad G$. It follows that
$\z_\chi(0)=[\g(0),\g(0)]$ is semisimple.

Let $I$ be any ideal of codimension $1$ in $H$. Then $xy-yx\in I$
for all $x,y\in H$. Since $\z_\chi(0)$ is semisimple, we have
$\z_\chi(0)=[\z_\chi(0),\z_\chi(0)]$. In view of Lemma~\ref{L4}
this implies that $\Theta_x\in I$ for all $x\in\z_\chi(0)$. We
have already mentioned that in the present case $\z_\chi(1)$ is an
irreducible $\z_\chi(0)$-module. So
$\z_\chi(1)=[\z_\chi(0),\z_\chi(1)]$ holds, yielding $\Theta_u\in
I$ for all $u\in \z_\chi(1)$; see Lemma~\ref{L5}. Since $I$ is a
subalgebra of $H$ containing $\z_\chi(0)\cup\z_\chi(1)$,
Proposition~\ref{P2} entails $C-c_0\in I$. As a consequence,
$H^+\subseteq I$. Since $H_\chi^+$ has codimension $1$ in $H$, we
conclude that $I=H^+$, as desired.
\end{pf}

\smallskip

\section{\bf Primitive ideals and Goldie rank}

\subsection{} We retain the assumptions of the previous section
and denote by ${\mathcal C}$ the category ${\mathcal C}_\chi$ for
$\chi=(e_\beta,\,\cdot\,)$. According to (3.1), given
$M\in{\mathcal C}$ one has $M\cong {Q_\chi}\otimes_{H} M_0$ where
$M_0={\mathrm Wh}(M)$. Let $\{m_i\,|\,i\in I\}$ be a basis of
$M_0$. It is immediate from [\cite{P02}, p.~52] that the
vectors$$\{h^l z_1^{i_1}\cdots z_s^{i_s}\otimes m_j\,|\,j\in I;
\,\, l,i_i,\ldots, i_s\in\Z_+\}$$ form a basis of $Q_\chi\otimes_H
M_0$ over $\k$. We can thus identify $M$ and
$\k[h,z_1,\ldots,z_s]\otimes M_0$ as vector spaces over $\k$.
Recall that $[z_i,z_j]=0$ and $[h,z_i]=-z_i$ for all $i,j\le s$.

Let $\Delta$ denote the automorphism of the polynomial algebra
$\k[h]$ such that $\Delta(h)=h+1$. Clearly, $\Delta^{-1}(h)=h-1$.
Let $\langle \Delta\rangle$ stand for the cyclic subgroup of
$\text{Aut}(\k[h])$ generated by $\Delta$. The skew group algebra
$\k[h]*\langle\Delta\rangle$ has the set
$\{h^i\Delta^j\,|\,i\in\Z_+,\,j\in\Z\}$  as a $\k$-basis and the
following relations hold: $$\Delta^i\cdot
h^k=(h+i)^k\cdot\Delta^i\qquad(i\in\Z,\ \,k\in\Z_+).$$

Let ${\mathcal A}_e:=\big(\k[h]*\langle\Delta\rangle\big)\otimes
{\mathbf A}_e$, an associative algebra over $\k$, and identify the
Weyl algebra ${\mathbf A}_e$ and the skew group algebra
$\k[h]*\langle\Delta\rangle$ with the subalgebras $\k\otimes
{\mathbf A}_e$ and $\big(\k[h]*\langle\Delta\rangle\big)\otimes\,
\k$ of ${\mathcal A}_e$, respectively. Define an involution
$\tau\in\text{Aut}({\mathcal A}_e)$ by setting
$$\tau(z_i)=-z_i,\quad
\tau(\partial_i)=-\partial_i,\quad\tau(h)=h,\quad\tau(\Delta^k)=(-1)^k\Delta^k
\qquad(1\le i\le s,\,k\in\Z).$$

The polynomial algebra $k[h,z_1,\dots,z_s]= \k[h][z_1,\ldots,z_s]$
has a natural ${\mathcal A}_e$-module structure such that $h\, .
f(h,z_1,\ldots,z_s)=hf(h,z_1,\ldots,z_s)$ and
$\Delta^k.f(h,z_1,\ldots,z_s)=f(h+k,z_1,\ldots,z_s)$. As this
${\mathcal A}_e$-module is faithful, we will identify ${\mathcal
A}_e$ with a subalgebra of $\End\big(\k[h,z_1,\dots,z_s]\big)$.
Since $\Delta^{2k}\in{\mathcal A}_e^\tau$ for all $k\in \Z$ and
$\Delta\otimes z_i, \Delta\otimes\partial_i\in{\mathcal A}_e^\tau$
for all $i\le s$, it is easy to see that $\k[h,z_1,\ldots,z_s]$
remains irreducible when restricted to the fixed point algebra
${\mathcal A}_e^\tau$.

Let $I$ be any two-sided ideal of $H$. Since $I$ is
$\sigma$-stable by Corollary~\ref{C}, the ideal ${\mathcal
A}_e\otimes I$ of the algebra ${\mathcal A}_e\otimes H$ is
invariant under the involution $\tau\otimes\sigma$ of ${\mathcal
A}_e\otimes H$. Hence $\tau\otimes\sigma$ indices an automorphism
of order two on the algebra ${\mathcal A}_e\otimes (H/I)$. We
mention for further references that $\big({\mathcal A}_e\otimes
(H/I)\big)^{\tau\otimes\sigma}$ contains the images in ${\mathcal
A}_e\otimes(H/I)$ of ${\mathcal A}_e^\tau\otimes\k$,
$1\otimes\Theta_a$, and $\Delta^{-1}\otimes \Theta_u$ for all
$a\in\z_\chi(0)$ and $u\in\z_\chi(1)$.

\subsection{} Recall that in our present setting the element $f$ is
a root vector for $\h$ corresponding to $-\beta$. Since $\ad f$ is
locally nilpotent on $U(\g)$, the set $S_f:=\{f^i\,|\,i\in\Z_+\}$
is an Ore set in $U(\g)$; see [\cite{Ja}, (11.2)]. We denote by
$U(\g)_f$ the localization $S^{-1}_f U(\g)$. Since $f$ commutes
with the $z_i$'s and $fh^m=(h+2)^mf$ for all $m\in \Z_+$, it
follows that $f$ acts on $M=\k[h,z_1,\ldots,z_s]\otimes M_0$ as
$\Delta^2\otimes\text{id}_{M_0}$ (one should keep in mind that
$\chi(f)=1$). In particular, $f$ acts invertibly on $M$. From this
it follows that the localization $S_f^{-1}M$ can be identified
with $M$; see [\cite{Ja}, (11.4)]. As a result, the action of
$U(\g)$ on $M$ extends uniquely to a representation of $U(\g)_f$
in $\End M$. We note for completeness that the enveloping algebra
$U(\g)$ is canonically identified with a subalgebra of $U(\g)_f$.
\begin{theorem}\label{S5}
Let $M\in{\mathcal C}$ and identify $M$ with
$\k[h,z_1,\ldots,z_s]\otimes M_0$ where $M_0={\mathrm Wh}(M)$. Let
$\tilde{\rho}\colon\,U(\g)_f\,\rightarrow\,\End M$ and
$\rho\colon\,H\rightarrow\End M_0$ be the representations of
$U(\g)_f$ and $H$ induced by the action of $\g$ on $M$. Then the
following hold:
\smallskip
\begin{enumerate}
\item[(i)\ ] $ \ (\Delta\otimes{\mathrm
id}_{M_0})\circ\tilde{\rho}(U(\g)_f)\circ(\Delta\otimes{\mathrm
id}_{M_0})^{-1}\,=\,\tilde{\rho}(U(\g)_f)$;
\smallskip
\item[(ii)\ ] $\ \tilde{\rho}(f)=\Delta^2\otimes{\mathrm
id}_{M_0}$;
\smallskip
\item[(iii)\ ] \ ${\mathcal A}_e\otimes\rho(H)\,=\,
\tilde{\rho}(U(\g)_f)\bigoplus\,\tilde{\rho}(U(\g)_f)(\Delta\otimes{\mathrm
id}_{M_0});$
\smallskip
\item[(iv)\ ] \ $\tilde{\rho}(U(\g)_f)\,=\,\big({\mathcal
A}_e\otimes \rho(H)\big)^{\tau\otimes\sigma}$.
\end{enumerate}
\end{theorem}
\begin{pf}
Put $\text{id}=\text{id}_{M_0}$. We have already mentioned that
$\tilde{\rho}(f)=\Delta^2\otimes\text{id}$, showing (ii). Since
$z_ih^m=(h+1)^mz_i$ for all $m\in \Z_+$ and $i\le 2s$, it follows
that
$$\tilde{\rho}(h)=h\otimes\text{id},\quad\tilde{\rho}(z_i)=(\Delta\circ
z_i)\otimes\text{id},\quad
\tilde{\rho}(z_{i+s})=(\Delta\circ\partial_i)\otimes\text{id}\qquad\
\, (1\le i\le s).$$ All these are in ${\mathcal
A}_e^{\tau}\otimes\text{id}\subset\big({\mathcal A}_e\otimes
\rho(H)\big)^{\tau\otimes\sigma}$. Now let $a\in\z_\chi(0)$ and
write $[a,z_i]=\sum_{i=1}^{2s}\,\mu_{ij}\,z_j$ for $1\le i\le 2s$,
where $\mu_{ij}\in\k$. Note that for $1\le i,j\le s$ we have
$[az_iz_j]=\mu_{i,j+s}f$, and
\begin{eqnarray*}
a z_1^{k_1}\cdots z_s^{k_s}&=&\sum_{i=1}^s k_i \,z_1^{k_1}\cdots
z_i^{k_i-1}\cdots z_s^{k_s}[a,z_i] +\frac{k_i(k_i-1)}{2}
\sum_{i=1}^s
z_1^{k_1}\cdots z_i^{k_i-2}\cdots z_s^{k_s}[az_i^2]\\
&+&k_ik_j\sum_{1\le i<j\le s} z_1^{k_1}\cdots z_i^{k_i-1}\cdots
z_j^{k_j-1}\cdots z_s^{k_s}[az_iz_j] +z_1^{k_1}\cdots z_s^{k_s}a
\end{eqnarray*}
for all $k_i\in\Z_+$. Since $\langle\,\cdot\,,\,\cdot\,\rangle$ is
$\z_\chi(0)$-invariant, it must be that
$\mu_{i+s,\,j+s}=-\mu_{ji}$, $\mu_{i,\,j+s}=\mu_{j,\,i+s}$, and
$\mu_{i+s,\,j}=\mu_{j+s,\,i}$ where $1\le i,j\le s$. In view of
Lemma \ref{L3} and the fact that $z_{i+s}\in\m_\chi$ for all $i\le
s$, we must have that $\tilde{\rho}(a)(1\otimes
m)=1\otimes\rho(\Theta_a)(m)-\frac{1}{2}\sum_{i=1}^s[a,z_{i+s}]z_i\otimes
1$ for all $m\in M_0.$ In conjunction with our earlier remarks
(and the fact that $[a,h]=0$) this yields
\begin{eqnarray}\label{eqq}
\tilde{\rho}(a)&=&1\otimes\Big(\rho(\Theta_a)+\sum_{i=1}^s\,\mu_{ii}\Big)+\Big(
 \sum_{i,j=1}^s\,\mu_{ij}\,z_j\partial_i\Big)\otimes\text{id} \nonumber\\
 &+&\Big(\sum_{i=1}^s\,\frac{\mu_{i,\,i+s}}{2}\,\partial_i^2+
 \sum_{1\le i<j\le
 s}\mu_{i,\,j+s}\,\partial_i\partial_j\Big)\otimes\text{id}\nonumber\\
 &-&
\Big(\sum_{i=1}^s\,\frac{\mu_{i,\,i+s}}{2}\,z_i^2+\sum_{1\le i<
j\le s} \mu_{i,j+s}\,z_iz_j\big)\otimes\text{id}.
\end{eqnarray}
This shows that $\tilde{\rho}(a)$ commutes with
$\Delta\otimes\text{id}$ and lies in $\big({\mathcal A}_e\otimes
\rho(H)\big)^{\tau\otimes\sigma}$. Since all operators
$(\Delta^2\circ z_iz_j)\otimes\text{id},$\,
$(\Delta^2\circ\partial_i\partial_j)\otimes\text{id}$, and
$(\Delta^2\circ z_j\partial_i)\otimes\text{id}$ are in
$\tilde{\rho}\big(U(\n_\chi)\big)$ by our remarks earlier in the
proof, and since
$\tilde{\rho}(f^{-1})=\Delta^{-2}\otimes\text{id}$, it also
follows that $1\otimes \rho(\Theta_a)$ is in
$\widetilde{\rho}\big(U(\g)_f\big)$ for all $a\in\z_\chi(0)$.

Now let $u\in\z_\chi(1)$. First note that
$[u,z_i]=[u,z_i]^\sharp+\frac{1}{2}(u,z_i)h$ and $uh^m=(h-1)^mu$
for all $i\le s$ and $m\in\Z_+$. Next observe that
\begin{eqnarray*}
uz_1^{k_1}\cdots z_s^{k_s}&=&\sum_{i=1}^s\,z_1^{k_1}\cdots
z_i^{k_i-1}\cdots z_s^{k_s}[u,z_i]
+\sum_{i=1}^s\,\frac{k_i(k_i-1)}{2}z_1^{k_1}\cdots
z_i^{k_i-2}\cdots z_s^{k_s}[uz_i^2]\\
&+& \sum_{1\le i<j\le s}k_ik_j\,z_1^{k_1}\cdots z_i^{k_i-1}\cdots
z_j^{k_j-1}\cdots z_s^{k_s}[uz_iz_j]\\
&+&\sum_{1\le i<j<l\le s}k_ik_jk_l\,z_1^{k_1}\cdots
z_i^{k_i-1}\cdots
z_j^{k_j-1}\cdots z_l^{k_l-1}\cdots z_s^{k_s}[uz_iz_jz_l]\\
&+&\sum_{1\le i<j\le s}\frac{k_i(k_i-1)}{2}k_j\,z_1^{k_1}\cdots
z_i^{k_i-2}\cdots
z_j^{k_j-1}\cdots z_s^{k_s}[uz_i^2z_j]\\
&+&\sum_{1\le i<j\le s}\,k_i\frac{k_j(k_j-1)}{2}\,z_1^{k_1}\cdots
z_i^{k_i-1}\cdots
z_j^{k_j-2}\cdots z_s^{k_s}[uz_iz_j^2]\\
&+&\sum_{i=1}^s\,\frac{k_i(k_i-1)(k_i-2)}{6}\,z_1^{k_1}\cdots
z_i^{k_i-3}\cdots z_s^{k_s}[uz_i^3]+z_1^{k_1}\cdots z_s^{k_s}u,
\end{eqnarray*}
and
\begin{eqnarray*}
z_1^{k_1}\cdots z_s^{k_s}h&=&hz_1^{k_1}\cdots
z_s^{k_s}-[h,z_1^{k_1}\cdots z_s^{k_s}]\,=\,hz_1^{k_1}\cdots
z_s^{k_s}\\
&+&(k_1+\cdots k_s)z_1^{k_1}\cdots z_s^{k_s}
\,=\,(h+\sum_{i=1}^s\,z_i\partial_i)(z_1^{k_1}\cdots z_s^{k_s}).
\end{eqnarray*}
This shows that for all $m\in M_0$ we have \begin{eqnarray*} &
&\sum_{i=1}^s\, z_1^{k_1}\cdots z_i^{k_i-1}\cdots
z_s^{k_s}\,\tilde{\rho}([u,z_i])(1\otimes
m)\,=\,\\
&=&\Big(\big(h+\sum_{i=1}^s\,z_i\partial_i\big)\circ
\big(\sum_{i=1}^s\,\frac{(u,z_i)}{2}\partial_i\big)+\sum_{i=1}^s\,
\partial_i\otimes\rho\big(\Theta_{[u,z_i]^\sharp}\big)\Big)(z_1^{k_1}\cdots
z_s^{k_s}\otimes m)\\
&-&\frac{1}{2}\sum_{i,j=1}^s\,z_1^{k_1}\cdots z_i^{k_i-1}\cdots
z_s^{k_s}\tilde{\rho}\big([[u,z_i]^\sharp,z_{j+s}]z_j\big)(1\otimes
m).
\end{eqnarray*}
In view of Lemma \ref{L5} we now get
\begin{eqnarray*}
& &z_1^{k_1}\cdots z_s^{k_s}\,\tilde{\rho}(u)(1\otimes
m)\,=\,\big(\text{id}\otimes \rho(\Theta_u)\big)(z_1^{k_1}\cdots
z_s^{k_s}\otimes m)\\&-&z_1^{k_1}\cdots
z_s^{k_s}\Big(\sum_{i=1}^s\,z_i\tilde{\rho}\big([u,z_{i+s}]\big)(1\otimes
m)+\frac{1}{3}\sum_{i,j=1}^s\,\tilde{\rho}\big([uz_{i+s}z_{j+s}]\big)(z_jz_i\otimes
m)\\
&+&\frac{2}{3}\sum_{i=1}^s\,\tilde{\rho}\big([uz_{i+s}z_i]\big)(1\otimes
m)+\tilde{\rho}(z_u)(1\otimes m)\Big).
\end{eqnarray*}
Together with the above remarks this says that $\tilde{\rho}(u)$
is a linear combination of the following operators:
\begin{eqnarray*}
& &\big(h+ \textstyle{\sum_{i=1}^s}\,z_i\partial_i\big)\circ
\big(\textstyle{\sum_{i=1}^s}\,\frac{(u,z_{i})}{2}\,\partial_i\big)
\circ\Delta^{-1}\otimes\text{id},\
\textstyle{\sum_{i=1}^s}\,\big(\partial_i\circ\Delta^{-1}\big)
\otimes\rho\big(\Theta_{[u,z_i]^\sharp}\big),\\
& &\big(h+\textstyle{\sum_{i=1}^s}\, z_i\partial_i\big)\circ
\big(\sum_{i=1}^s\,\frac{(u,z_{i+s})}{2}\,z_i\big)\circ\Delta^{-1}
\otimes\text{id},\ \sum_{i=1}^s\,\big(
z_i\circ\Delta^{-1}\big)\otimes\rho\big(\Theta_{[u,z_{i+s}]^\sharp}\big),\\
&&(\partial_i\partial_j\partial_k\circ\Delta^{-1})\otimes\text{id},\
(z_iz_jz_k\circ\Delta^{-1})\otimes \text{id},\ (
z_iz_j\partial_k\circ\Delta^{-1})\otimes\text{id},\\
& &(z_k\partial_i\partial_j\circ\Delta^{-1})\otimes\text{id},\
(\partial_i\circ\Delta^{-1})\otimes\text{id},\
(z_i\circ\Delta^{-1})\otimes\text{id},\
\Delta^{-1}\otimes\rho(\Theta_u),
\end{eqnarray*}
where $i,j,k\le s.$ But then $\tilde{\rho}(u)\in\big({\mathcal
A}_e\otimes\rho(H)\big)^{\sigma\otimes\tau}$ and
$(\Delta\otimes\text{id})\,\tilde{\rho}(u)
(\Delta\otimes\text{id})^{-1}-\tilde{\rho}(u)\in{\mathcal
A}_e^\tau\otimes\,\text{id}$. As ${\mathcal
A}_e^\tau\otimes\,\text{id}\subset\tilde{\rho}(U(\g)_f)$ by our
earlier remarks, this yields
$$(\Delta\otimes\text{id})\,\tilde{\rho}(x)\,
(\Delta\otimes\text{id})^{-1}\in
\tilde{\rho}(U(\g)_f)\qquad(\forall\,x\in\textstyle{\bigcup}_{i\le
1}\,\g(i)).$$ Since $\g(2)=[\g(1),\g(1)]$, we obtain (i) and the
inclusion $\tilde{\rho}(U(\g)_f)\subseteq\big({\mathcal
A}_e\otimes\rho(H)\big)^{\sigma\otimes\tau}$. Notice that
$\Delta^{-1}\otimes\rho(\Theta_u)\in{\mathcal
A}_e^\tau\otimes\text{id}+\tilde{\rho}(U(\g)_f)$ for all
$u\in\z_\chi(1)$ and $1\otimes\rho(\Theta_a)\in{\mathcal
A}_e^\tau\otimes\text{id}+\tilde{\rho}(U(\g)_f)$ for all
$a\in\z_\chi(0)$; see (\ref{eqq}). Since the algebra $H$ is
generated by the elements $\Theta_x$ with
$x\in\z_\chi(0)\cup\z_\chi(1)$, by Proposition~\ref{P2}, we get
$$1\otimes\rho(H)\subset
\tilde{\rho}(U(\g)_f)+\tilde{\rho}(U(\g)_f)(\Delta\otimes\text{id}).$$
Since ${\mathcal A}_e\,=\,{\mathcal A}_e^\tau+{\mathcal
A}_e^\tau\,\Delta$ and
$\tilde{\rho}(U(\g)_f)\subseteq\big({\mathcal
A}_e\otimes\rho(H)\big)^{\sigma\otimes\tau}$ we derive (iii). Then
(iv) follows, completing the proof.
\end{pf}
\subsection{} By [\cite{P02}, (6.2)], the restriction of
the representation $\tilde{\rho}_\chi\colon\,U(\g)\rightarrow
\End(Q_\chi)$ to the centre $Z(\g)$ of $U(\g)$ is injective. Since
any nonzero two-sided ideal of $U(\g)$ intersects with $Z(\g)$ by
[\cite{Dix}, (4.2.2)], it follows that $\tilde{\rho}_\chi$ is a
faithful representation of $U(\g)$. Note that $Q_\chi\in\mathcal
C$ and $\text{Wh}(Q_\chi)$ is canonically identified with $H$
(viewed as the left regular $H$-module). By (5.1), $Q_\chi\cong
Q_\chi\otimes_H H$ is then identified with ${\mathcal A}_e\otimes
H$ as vector spaces over $\k$. Since $\tilde{\rho}$ is faithful,
it extends to a faithful representation of $U(\g)_f$ in
$\End(Q_\chi)$ (one should take into account that
$\tilde{\rho}_\chi(f)$ is invertible). Since $\text{Wh}(Q_\chi)$
is identified with the left regular $H$-module,
Theorem~\ref{S5}(iv) yields $\tilde{\rho}_\chi(U(\g)_f)\subseteq
{\mathcal A}_e\otimes H$.  Applying Theorem~\ref{S5} in this
situation we now obtain:
\begin{corollary}\label{C5} Set $\widehat{\Delta}=\Delta\otimes{\mathrm id}_H$
and identify $U(\g)_f$ with its image in ${\mathcal A}_e\otimes
H$. Then the following hold:
\smallskip
\begin{itemize}
\item[(i)\ ] $ \ \widehat{\Delta}\, U(\g)_f\,
\widehat{\Delta}^{-1} \,=\, U(\g)_f$;
\smallskip
\item[(ii)\ ] \ $f=\,\widehat{\Delta}^2$;
\smallskip
\item[(iii)\ ] \ ${\mathcal A}_e\otimes H \,=\,
U(\g)_f\bigoplus\,U(\g)_f\,\widehat{\Delta};$
\smallskip
\item[(iv)\ ] \ $\ U(\g)_f\,=\,\big({\mathcal A}_e\otimes
H\big)^{\tau\otimes\sigma}$.
\smallskip
\item[(v)\ ]  \ $\big({\mathcal A}_e\otimes H\big)^{\n_\chi}
\,=\,\k[\Delta,\Delta^{-1}]\otimes H\,=\,
\big(U(\g)^{\n_\chi}\big)_f\bigoplus\,\big(U(\g)^{\n_\chi}\big)_f\,\widehat{\Delta};$
\smallskip \item[(vi)\ ] \ $Z(H)\,\cong\, Z(U(\g)_f)\,\cong\,Z(
\g).$
\end{itemize}
\end{corollary}
\begin{pf}
Parts~(i)--(iv) follow from Theorem~\ref{S5} and the discussion
above. For (v), note that $(\Delta\circ z_i)\otimes
1,\,(\Delta\circ
\partial_i)\otimes 1\in \n_\chi$ for all $i\le s$ (see the
beginning of the proof of Theorem~\ref{S5}).  From this it is
immediate that $\big({\mathcal A}_e\otimes H\big)^{\n_\chi}
\,=\,\k[\Delta,\Delta^{-1}]\otimes H$. Since $f\in
Z(U(\g)^{\n_\chi})$ we have
$\big(U(\g)_f\big)^{\n_\chi}\,=\,\big(U(\g)^{\n_\chi}\big)_f$,
hence the result.

For, (vi), we first recall that the action of $1\otimes\sigma$ on
${\mathcal A}_e\otimes H$ is induced by the adjoint action of
$1\otimes\Theta(\z_\chi(0))\subset 1\otimes H$. Consequently,
$$Z(H)\,\cong\, Z({\mathcal A}_e\otimes H)\,=\, \k\otimes
(H^\sigma\cap Z(H))\,=\, \big(Z({\mathcal A}_e\otimes
H)\big)^{\tau\otimes\sigma}\subseteq \,Z(U(\g)_f),$$ by (iv). On
the other hand,  $Z(U(\g)_f)\subseteq
Z((U(\g)^{\n_\chi})_f)\subseteq \k[\Delta,\Delta^{-1}]\otimes
Z(H),$ by (v). Since $[h,\Delta]=-\Delta$, it must be that
$Z(U(\g)_f)\subseteq\, \k\otimes Z(H)=Z({\mathcal A}_e\otimes H)$.
Therefore, $Z(H)\cong Z(U(\g)_f)$.

It remains to show that $Z(U(\g)_f)=Z(\g)$. This is easy and must
be well-known, but we could not find a good reference. One can
argue as follows: Let $K(\g)=\text{Fract}\,U(\g)$ be the Lie field
of $\g$. By [\cite{Dix}, (4.3.2)], the centre of $K(\g)$ coincides
with $\text{Fract}\,Z(\g)$. Let $z\in Z(U(\g)_f)$. Then
$z=f^{-d}\,c$ for some $c\in U(\g)$ and $d\in\Z_+$. Since
$Z(U(\g)_f)\subset Z(K(\g))$, there are $a,b\in Z(\g)$ such that
$af^d=bc$. By a classical result of Kostant, $U(\g)\cong
Z(\g)\otimes {\mathsf H}(\g)$ as $Z(\g)$-modules, where ${\mathsf
H}(\g)$ denotes the subspace of $U(\g)$ spanned by the powers of
the nilpotent elements of $\g$; see [\cite{Dix}, (8.2.4) and
(8.5.5)]. Choose a basis $\{u_i\}$ in $\mathsf{H}(\g)$ with
$u_1=f^d$. Write $c=\sum_{i} z_iu_i$ with $z_i\in Z(\g)$. Then
$a=bz_1$ and $z_i=0$ for $i\ne 1$, forcing
$z=f^{-d}c=f^{-d}z_1f^d=z_1\in Z(\g)$. The result follows.
\end{pf}
\subsection{} Recall that a two-sided ideal $I$ of an associative
ring $R$ is called {\it prime} if $I\ne R$ and for any two
two-sided ideals $J_1, J_2$ the inclusion $J_1J_2\subseteq I$
implies that either $J_1\subseteq I$ or $J_2\subseteq I$. We let
$\text{Spec}\,R$ denote the set of all prime ideals of $R$. Any
primitive ideal of $R$ is prime. The ring $R$ is termed {\it
prime} if $(0)\in\text{Spec}\,R$.

According to [\cite{Dix}, (3.6.17)], if $R$ is Noetherian and $S$
is an Ore set in $R$, then the mapping $I\mapsto S^{-1}I$ sets up
a bijection  between the subset
$(\text{Spec}\,R)_S:=\{I\in\text{Spec}\,R\,|\,I\cap S=\emptyset\}$
of $\text{Spec}\,R$ and $\text{Spec}\,S^{-1}R$. On the other hand,
the complement to $\big(\text{Spec}\,U(\g)\big)_{S_f}=\,\{I\in
\text{Spec}\,U(\g)\,\,|\,f^n\not\in I\,\, \, \text{for all
}\,n\in\Z_+\}$ in $\text{Spec}\,U(\g)$ consists of all prime
ideals of finite codimension in $U(\g)$; see [\cite{Ja},
Lemma~13.17] for example. Thus the mapping $I\mapsto S_{f}^{-1}I$
is a bijection between the set of all prime ideals
 of infinite codimension in $U(\g)$ and the set of all prime  ideals of
 $U(\g)_f$.

By Goldie's theorem, the set $S$ of all regular elements of a
prime Noetherian ring $R$ is an Ore set in $R$ and the
localisation ${\mathcal Q}(R):=S^{-1}R$ is isomorphic to the
matrix algebra $\text{Mat}_n(K)$ over a noncommutative field $K$.
Both $K$ and $n$ can be described intrinsically, hence are
uniquely determined by $R$. They are called the {\it Goldie field}
and the {\it Goldie rank} of $R$, respectively. We write
$n=\text{rk}(R)$.
\subsection{}
Corollary~\ref{C5}(iv) allows us to identify $U(\g)_f$ with the
subalgebra $\big({\mathcal A}_e\otimes H\big)^{\tau\otimes\sigma}$
of ${\mathcal A}_e\otimes H$. Since ${\mathcal A}_e\,=\,{\mathcal
A}_e^\tau\oplus{\mathcal A}_e^\tau\,\Delta$ and
$\tau(\Delta)=-\Delta$, we obtain the decomposition
$$U(\g)_f=\,{\mathcal A}_e^\tau\otimes
H_+\,\textstyle{\bigoplus}\, \,{\mathcal A}_e^\tau\Delta\otimes
H_-$$ where $H_+=\{x\in H\,|\,\sigma(x)=x\}$ and $H_-=\{x\in
H\,|\,\sigma(x)=-x\}$. Let $I$ be any two-sided ideal of $H$. We
let $\text{Dim}(H/I)$ denote the Gelfand--Kirillov dimension of
the factor algebra $H/I$. Since $I$ is $\sigma$-stable
(Corollary~\ref{C}), we have $I=I_+\oplus I_{-}$ where
$I_\pm=I\cap H_\pm$. Put
\begin{eqnarray}\label{eqS5}
\widetilde{I}_f:=\,{\mathcal A}_e^\tau\otimes
I_+\,\textstyle{\bigoplus}\, \,{\mathcal A}_e^\tau\Delta\otimes
I_-\quad\,\text{and}\,\quad
\,\widetilde{I}\,:=\,\,\widetilde{I}_f\cap\, U(\g).
\end{eqnarray}
Clearly, $\widetilde{I}_f$ and $\widetilde{I}$ are two-sided
ideals of $U(\g)_f$ and $U(\g)$, respectively. Let ${\cal
O}_{\mathrm min}$ denote the minimal nilpotent orbit $(\Ad^*
G)\cdot\chi$ in $\g^*$. Recall that $$\dim {\cal O}_{\mathrm
min}\,=\,2\dim\m_\chi\,=\,2(s+1).$$
\begin{theorem}\label{main1} The following are true:
\begin{enumerate}

\item[(i)] The map $I\mapsto \widetilde{I}_f$ sets up a bijection
between the set of all two-sided ideals of $H$ and the set of all
two-sided ideals of $U(\g)_f$. For any two-sided ideal $I$ of $H$
one has $${\mathrm Dim}(U(\g)_f/\widetilde{I}_f)=\,{\mathrm
Dim}(H/I)+\dim {\cal O}_{\mathrm min}.$$

\item[(ii)] The map $I \mapsto \widetilde{I}$ is a bijection
between ${\mathrm Spec}\,H$ and the set of all prime ideals of
infinite codimension in $U(\g)$. Furthermore,
$$\qquad\ \, {\mathrm Dim}(U(\g)/\widetilde{I})=\,{\mathrm
Dim}(H/I)+\dim {\cal O}_{\mathrm min}\quad\ (\forall\,I\in{\mathrm
Spec}\,H).$$
\end{enumerate}
\end{theorem}
\begin{pf}
(a) Let $\Xi$ denote the set of all quadruples $(m,n,{\mathbf
i},\mathbf{j})$ with $m\in\Z_+$, $n\in\Z$, and ${\mathbf
i},{\mathbf j}\in\Z_+^s$. Order the elements in $\Xi$
lexicographically. Given $\xi=(m,n,{\mathbf i},\mathbf{j})\in\Xi$
define $a_\xi\in {\mathcal A}_e$ by setting
$a_\xi=h^m\Delta^nz_1^{i_1}\cdots
z_s^{j_s}\,\partial_1^{j_1}\cdots
\partial_s^{j_s}$, where
${\mathbf i}=(i_1,\ldots, i_s)$ and ${\mathbf j}=(j_1,\ldots,
j_s)$. Any nonzero $x\in {\mathcal A}_e\otimes H$ can be written
uniquely as $x=\sum_{\,\xi\in\Xi(x)}\,a_\xi\otimes h_\xi$ for some
nonzero $h_\xi\in H$. Here $\Xi(x)$ is a finite subset of $\Xi$
depending on $x$.

Let $\cal I$ be any two-sided ideal of $U(\g)_f=\big({\mathcal
A}_e\otimes H\big)^{\tau\otimes\sigma}$. Recall that the action of
$\sigma$ on $H$ is induced by the adjoint action of
$\Theta(\z_\chi(0))$. Since $1\otimes\Theta(\z_\chi(0))\subset
\big({\mathcal A}_e\otimes H\big)^{\tau\otimes\sigma}$, the ideal
$\cal I$ is stable under the involution $1\otimes\sigma$ of
${\mathcal A}_e\otimes H$. It follows that ${\cal I}={\cal
I}_+\oplus {\cal I}_-$ where ${\cal I}_\pm=\{x\in{\cal
I}\,|\,(1\otimes\sigma)(x)=\pm x\}$. Let $x\in \cal I_+\cup\cal
I_-$ and let $\xi_0=(m_0,n_0,{\bf a},{\bf b})$ be the maximal
element in $\Xi(x)$.  Then there exists a polynomial
$f_x(t)\in\k[t]$ such that
$$f_x(\ad\,h)\circ\Big(\prod_{i=1}^s(\ad\, \Delta\otimes
\partial_i)^{a_i}\circ\prod_{i=1}^s(\ad \,\Delta\otimes
z_i)^{b_i}\Big)\circ(\ad\,
\widehat{\Delta}^{2m_0})(x)\in\,\k^\times\Delta^{N+n_0}\otimes
h_{\xi_0}$$ where $N=2m_0+\sum_{i=1}^s(a_i+b_i)$. As
$\widehat{\Delta}^2$ is invertible, it follows that
$\Delta^{\varepsilon(x)}\otimes h_{\xi_0}\in{\cal I}_\pm$, where
$\varepsilon(x)=0$ if $x\in\cal I_+$ and $\varepsilon(x)=1$ if
$x\in\cal I_-$. But then $a_{\xi_0}\Delta^{-\varepsilon(x)}\otimes
1\in{\mathcal A}_e^\tau\otimes\k\subset\big({\mathcal A}_e\otimes
H\big)^{\tau\otimes\sigma}$ yielding $a_{\xi_0}\otimes
h_{\xi_0}\in \cal I_\pm$. Continuing this process one eventually
observes that $a_\xi\otimes h_\xi\in\cal I_\pm$ for all
$\xi\in\Xi(x)$. This implies that there is a graded subspace
$I=I_+\oplus I_-$ in $H=H_+\oplus H_-$ such that ${\cal
I}={\mathcal A}_e^\tau\otimes I_+\bigoplus\, {\mathcal
A}_e^\tau\Delta\otimes I_-$. Since $1\otimes
H_+\cup\,\Delta\otimes H_-\subset\big({\mathcal A}_e\otimes
H\big)^{\tau\otimes\sigma}$ and $\widehat{\Delta}^2$ is
invertible, it follows that $I$ is a two-sided ideal of $H$. As a
result, ${\cal I}=\widetilde{I}_f$ showing that the map $I\mapsto
\widetilde{I}_f$ is surjective. The injectivity of this map is
obvious.

\smallskip

\noindent (b) Let $\Xi_{\ge 0}$ be the subset of $\Xi$ consisting
of all $(m,n,{\bf i},{\bf j})$ with $n\ge 0$, and let ${\mathcal
A}_e'$ denote the $\k$-span of all $a_\xi$ with $\xi\in\Xi_{\ge
0}$. Clearly, ${\mathcal A}_e'$ is a $\tau$-invariant subalgebra
of ${\mathcal A}_e$. For $l\in\Z_+$ we let ${\mathcal A}_{e,l}'$
denote the $\k$-span of all $a_\xi$ with $\xi=(m,n,{\bf i},{\bf
j})\in\Xi_{\ge 0}$ such that $m+n+\sum_{k=1}^s(i_k+j_k)\le l$. It
is easy to see that $\{{\mathcal A}_{e,l}'\,\,|\,\,l\in\Z_+\}$ is
an increasing $\tau$-invariant filtration in ${\mathcal A}_e'$ and
the corresponding graded algebra $\text{gr}\,{\mathcal A}_e'$ is
isomorphic to the graded polynomial algebra $\k[X_1,\ldots,
X_{2s+2}]$ with all $X_i$ having degree $1$. Note that
$\ad(\Delta^2)$ is locally nilpotent on ${\mathcal A}_e'$ and the
algebra ${\mathcal A}_e$ identifies with the localisation
$({\mathcal A}_e')_{\Delta ^2}$.

Let $I$ be a two-sided ideal of $H$. It follows from part~(a) of
this proof that the two-sided ideal
$\widetilde{I}_f\oplus\widetilde{I}_f\widehat{\Delta}$ of
${\mathcal A}_e\otimes H$ coincides with ${\mathcal A}_e\otimes
I$, so that $\widetilde{I}_f\,=\big({\mathcal A}_e\otimes
I\big)^{\tau\otimes\sigma}$. Therefore, the involution
$\tau\otimes\sigma$ acts on the algebra ${\mathcal A}_e\otimes
(H/I)\cong ({\mathcal A}_e\otimes H)/({\mathcal A}_e\otimes I)$ in
such a way that the quotient $U(\g)_f/\widetilde{I}_f$ identifies
with $\big({\mathcal A}_e\otimes (H/I)\big)^{\tau\otimes\sigma}$.
Since $f=\widehat{\Delta}^2$, it is straightforward to see that
$$\big({\mathcal A}_e\otimes
(H/I)\big)^{\tau\otimes\sigma}\cong\big(({\mathcal
A}_e')_{\Delta^2}\otimes
(H/I)\big)^{\tau\otimes\sigma}\cong\big(({\mathcal A}_e'\otimes\,
(H/I))_{\bar{f}}\big)^{\tau\otimes\sigma}\cong\big(({\mathcal
A}_e'\otimes (H/I))^{\tau\otimes\sigma}\big)_{\bar{f}},$$ where
$\bar{f}$ denotes the image of $f$ in ${\mathcal A}_e\otimes
(H/I)$. In view of [\cite{BK}, (6.3)] we then have
\begin{eqnarray}\label{gk5}
{\mathrm Dim}(U(\g)_f/\widetilde{I}_f)\,=\,{\mathrm
Dim}\big(({\mathcal A}_e'\otimes (H/I)\big)^{\tau\otimes\sigma};
\end{eqnarray}
see also [\cite{Ja}, (11A.2)]. The Kazhdan filtration
$\{H^k\,|\,\,k\in\Z_+\}$ gives rise to the natural filtration
$\{(H/I)^k=H^k/(H^k\cap I)\,|\,\,k\in\Z_+\}$ of the algebra $H/I$.
Thanks to [\cite{P02}, Theorem~4.6(iii)], $\text{gr}\,H$ is a
Noetherian commutative $\k$-algebra. Hence so is the corresponding
graded algebra $\text{gr}(H/I)\,\cong\,\text{gr}\,H/\text{gr}\,I$.
Since $\sigma$ preserves both $I$ and the Kazhdan filtration of
$H$, it induces an automorphism of the graded algebra $\gr(H/I)$.

Next we observe that the subspaces $$\big\{\big({\mathcal
A}_e'\otimes (H/I)\big)_k\,=\,{\textstyle\sum_{i+j\le
k}}\,\,{\mathcal A}_{e,i}'\otimes (H/I)^j\,|\,\,\,k\in\Z_+\big\}$$
form an increasing filtration of the algebra ${\mathcal
A}_e'\otimes (H/I)$ such that
$$\text{gr}\big({\mathcal A}_e'\otimes (H/I)\big)\cong\,\text{gr}({\mathcal
A}_e')\otimes \text{gr}(H/I)\cong\,\k[X_1,\ldots,
X_{2s+2}]\otimes\text{gr}(H/I).$$ By construction, the involution
$\tau\otimes\sigma$ acts on the graded algebra
$\text{gr}({\mathcal A}_e'\otimes (H/I))$ and
$$\big(\text{gr}({\mathcal A}_e'\otimes
(H/I))\big)^{\tau\otimes\sigma}\cong\,\,\text{gr}\big(({\mathcal
A}_e'\otimes (H/I))^{\tau\otimes\sigma}\big)$$ as graded algebras.
Since the morphism
$$\text{Spec}\,\text{gr}\big({\mathcal A}_e'\otimes
(H/I)\big)\,\longrightarrow \,\text{Spec}\big(\text{gr}(({\mathcal
A}_e'\otimes (H/I))\big)^{\tau\otimes\sigma}$$ induced by
inclusion $\big(\text{gr}({\mathcal A}_e'\otimes
(H/I))\big)^{\tau\otimes\sigma}\hookrightarrow\,\text{gr}\big({\mathcal
A}_e'\otimes (H/I)\big)$ is finite, the Noetherian $\k$-algebras
$\big(\text{gr}({\mathcal A}_e'\otimes
(H/I))\big)^{\tau\otimes\sigma}$ and $\text{gr}({\mathcal
A}_e')\otimes \text{gr}(H/I)$ have the same Krull dimension. Since
the Krull dimensions of the graded algebras
$\big(\text{gr}({\mathcal A}_e'\otimes
(H/I))\big)^{\tau\otimes\sigma}$ and $\text{gr}(H/I)$ coincide
with the degrees of their respective Hilbert polynomials, we
derive that
$$\text{Dim}\,\text{gr}\big(({\mathcal A}_e'\otimes
(H/I))^{\tau\otimes\sigma}\big)\,=\,2(s+1)+\text{Dim}\,\text{gr}(H/I).$$
On the other hand, it follows from [\cite{McR}, Proposition~8.6.6]
that $\text{Dim}\,\text{gr}(H/I)\,=\,\text{Dim}(H/I)$ and
$\text{Dim}\,\text{gr}\big(({\mathcal A}_e'\otimes
(H/I))^{\tau\otimes\sigma}\big)\,=\,\,\text{Dim}\,\big(({\mathcal
A}_e'\otimes (H/I))^{\tau\otimes\sigma}\big)$. Combining this with
(\ref{gk5}) we get (i).

\smallskip

\noindent (c) Let $I\in\text{Spec}\,H$ and suppose ${\cal
J}_1{\cal J}_2\subseteq \widetilde{I}_f$ for some two-sided ideals
${\cal J}_1$ and ${\cal J}_2$ of $U(\g)_f=\big({\mathcal
A}_e\otimes H\big)^{\tau\otimes\sigma}$. By part~(a), there exist
two-sided ideals $J_1$ and $J_2$ of $H$ such that ${\cal
J}_i=\,{\mathcal A}_e^\tau\otimes J_{i,+}\,\textstyle{\bigoplus}\,
\,{\mathcal A}_e^\tau \Delta\otimes J_{i,-}$, $\,i=1,2$. Since
$\widehat{\Delta}^2$ is invertible, it is easy to see that
$J_1J_2\subseteq I$. Then either $J_1\subseteq I$ or $J_2\subseteq
I$, and hence either ${\cal J}_1\subseteq \widetilde{I}_f$ or
${\cal J}_2\subseteq \widetilde{I}_f$. As a result,
$\widetilde{I}_f\in\text{Spec}\,U(\g)_f$ and
$\widetilde{I}=\widetilde{I}_f\cap U(\g)\in\text{Spec}\,U(\g)$.

Conversely, let ${\cal I}$ be a prime ideal of infinite
codimension in $U(\g)$. The discussion at the beginning of this
subsection shows that $S_f^{-1}{\cal I}\in \text{Spec}\,U(\g)_f$.
By part~(a) of this proof, there is a two-sided ideal $I$ of $H$
such that $S_f^{-1}{\cal I}\,=\,{\mathcal A}_e^\tau\otimes
I_{+}\,\textstyle{\bigoplus}\, \,{\mathcal A}_e^\tau \Delta\otimes
I_{-}$. Clearly, ${\cal I}=\widetilde{I}_f\cap
U(\g)=\widetilde{I}$. If $PQ\subseteq I$ for some two-sided ideals
$P,Q$ of $H$, then $\widetilde{P}_f\widetilde{Q}_f\subseteq
S_f^{-1}\cal I$, forcing either $\widetilde{P}_f\subseteq
S_f^{-1}\cal I$ or $\widetilde{Q}_f\subseteq S_f^{-1}\cal I$. As
$S_f^{-1}{\cal I}=\widetilde{I}_f$, we obtain that either
$P\subseteq I$ or $Q\subseteq I$. Therefore, $I\in\text{Spec}\,H$.

Now let $\bar{f}$ be the image of $f$ in $U(\g)/\widetilde{I}$
where $I\in\text{Spec}\,H$. Since $\,U(\g)_f/\widetilde{I}_f\cong
(U(\g)/\widetilde{I})_{\bar{f}}$ by the exactness of localisation,
it follows from [\cite{BK}, (6.3)] and (i) that
$\text{Dim}(U(\g)/\widetilde{I})=\,{\mathrm Dim}(H/I)+\dim {\cal
O}_{\mathrm min}.$
\end{pf}
\subsection{} Given a ring $R$ we let $\text{Prim}\,R$ denote the primitive spectrum of $R$,
 the set of all primitive ideals of $R$ taken with the Jacobson topology. Set
 ${\mathcal X}:=\text{Prim}\,U(\g)$ and denote by ${\mathcal X}_{\rm fin}$ the set of all primitive ideals
 of finite codimension in $U(\g)$. Using the highest weight theory and [\cite{Dix}, (2.5.6),
 (3.2.3)] it is easy to observe that ${\mathcal X}_{\rm fin}$ is a
 countable dense subset of $\mathcal X$. The topology of $\mathcal X$ induces that on the complement
 ${\mathcal X}_{\rm inf}:={\mathcal X}\setminus {\mathcal X}_{\rm
 fin}$.

Recall that $\text{Prim}\,H$ is a subset of $\text{Spec}\,H$.
 By Theorem~\ref{main1}(ii), the map $I\mapsto\widetilde{I}$ given by
 (\ref{eqS5}) sets up a bijection between
 $\text{Spec\,}H$ and the set of all prime ideals of infinite
 codimension in $U(\g)$. Identify $Z(\g)$ with $Z(H)$
 according to
 Corollary~\ref{C5}(vi).

\begin{theorem}\label{main2}
The following are true:

\smallskip

\begin{enumerate}
\item The map $I\mapsto\widetilde{I}$ takes ${\mathrm Prim}\,H$
onto ${\mathcal X}_{\mathrm inf}$ and induces a homeomorphism of
topological spaces  ${\mathrm
Prim}\,H\stackrel{\sim}{\longrightarrow}\, {\mathcal X}_{\mathrm
inf}$.

\smallskip

\item Let $V$ be a finite dimensional irreducible $H$-module and
$I={\mathrm Ann}_H\, V$. Then $\widetilde{I}={\mathrm
Ann}_{U(\g)}\big(Q_\chi\otimes_H V\big)$ and $\,\,{\mathrm
rk}(U(\g)/\widetilde{I})\,=\,\dim{\mathrm Wh}(Q_\chi\otimes_H
V)\,=\,\dim V.$

\smallskip

\item Let $\cal I$ be a primitive ideal of $U(\g)$ with ${\mathcal
VA}({\cal I})\,=\,\overline{{\cal O}}_{\mathrm min}$. Then there
is a finite dimensional irreducible $H$-module $E$ such that
${\cal I}=\,{\rm Ann}_{U(\g)}\big(Q_\chi\otimes_H E\big)$.

\smallskip

\item Let $V_1$ and $V_2$ be two finite dimensional irreducible
$H$-modules. Then $V_1\cong V_2$ as $H$-modules if and only if
$\,{\mathrm Ann}_{U(\g)}\big(Q_\chi\otimes_H V_1\big)\,=\,{\mathrm
Ann}_{U(\g)}\big(Q_\chi\otimes_H V_2\big).$

\smallskip

\item For any algebra homomorphism $\eta\colon\,Z(\g)\rightarrow\,
\k$ there is a bijection between the isoclasses of finite
dimensional irreducible $H$-modules with central character $\eta$
and the primitive ideals $\cal I$ of $U(\g)$ with ${\cal I} \cap
Z(\g)\,=\,{\mathrm Ker}\,\eta$ and ${\mathcal VA}({\cal
I})\,=\,\overline{{\cal O}}_{\mathrm min}$.

\smallskip

\item A prime ideal $I$ of $H$ is primitive if and only if $I\cap
Z(H)$ is a maximal ideal of $Z(H)$.

\end{enumerate}
\end{theorem}
\begin{pf}
(a) Let $\{J^\alpha\,|\,\, \alpha\in A\}$ be a set of two-sided
ideals of $H$, and $J=\cap_{\,\alpha\in A}\,J^\alpha$. Since any
two-sided ideal of $H$ is $\sigma$-stable, we have
$$J_{\pm}=\{x\pm\sigma(x)\,|\,x\in J\}=\,\cap_{\,\alpha\in
A}\,\{x\pm\sigma(x)\,|\,x\in J^\alpha\}=\,\cap_{\,\alpha\in
A}\,J^\alpha_\pm.$$ Using (\ref{eqS5}) it is now easy to deduce
that $\widetilde{J}_f=\,\cap_{\,\alpha\in
A}\,\widetilde{J}^{\alpha}_f$ and
$\widetilde{J}=\,\cap_{\,\alpha\in A}\,\widetilde{J}^{\alpha}$.
Arguing similarly and using Theorem~\ref{main1} one also observes
that given $I,J\in\text{Spec}\,H$ one has $I\subsetneq J$ if and
only if $\widetilde{I}\subsetneq \widetilde{J}$.

\smallskip

\noindent (b) Let $I\in\text{Prim}\,H$. Then $I=\text{Ann}_H M_0$
for some irreducible $H$-module $M_0$. Let $M=Q\otimes_H M_0$ and
identify $M$ with $\k[h,z_1,\ldots,z_s]\otimes M_0$; see (5.1) for
detail. Theorem~\ref{main1} shows that $\text{Ann}_{U(\g)}
M=\widetilde{J}$ for some $J\in\text{Spec}\,H$, while
Theorem~\ref{S5} yields $\widetilde{J}_f=\widetilde{I}_f$. Then
$I=J$ in view of Theorem~\ref{main1}(i), forcing
$\widetilde{I}=\text{Ann}_{U(\g)} M$. As a consequence, the map
$I\mapsto\widetilde{I}$ takes $\text{Prim}\, H$ into ${\mathcal
X}_{\mathrm inf}$.

\smallskip

\noindent (c) Now suppose that $\widetilde{I}\in{\mathcal X}_{\rm
inf}$ for some two-sided ideal $I$ of $H$. By
Theorem~\ref{main1}(ii), $I\in\text{Spec}\,H$. Let ${\frak I}=
\{J\in{\mathrm Spec}\,H\,|\,\,J\supsetneq I\}$, and
$J_0=\cap_{J\in\frak I}\, J$. Our discussion in part~(a) in
conjunction with Theorem~\ref{main1}(ii) implies that
$\widetilde{\frak I}:=\{\widetilde{J}\,|\,\,J\in\frak I\}$
coincides with the set of all prime ideals of $U(\g)$ containing
$\widetilde{I}$ properly, and $\widetilde{J}_0=\,\cap_{\,{\cal
I}\in\widetilde{\frak I}}\,\,\cal I$. Since $\widetilde{I}$ is a
primitive ideal, [\cite{Dix}, (8.5.7)] applies yielding
$\widetilde{J}_0\supsetneq \widetilde{I}$. But then $J_0\supsetneq
I$ by our concluding remark in part~(a).

\smallskip

\noindent (d) We claim that the prime ideal $I$ from part~(c) is
the intersection of some primitive ideals of $H$. To see this one
can mimic the proof of Proposition~3.1.15 in [\cite{Dix}] which
deals with enveloping algebras but applies to a larger class of
filtered rings. For the reader's convenience we include the
argument which goes back to Duflo. As in {\it op.\,cit.} we put
$B:=H/I$, let $X$ be a variable, and set $C:=B\otimes \k[X]$. The
Kazhdan filtration of $H$ induces a filtration of $B$, which in
turn gives rise to a filtration of $C$. Since $\gr\,{H}$ is
finitely generated and commutative, by [\cite{P02}
Theorem~4.6(iii)], so are $\gr\,B$ and $\gr\,C$. Let $a\in J(B)$
where $J(B)$ is the Jacobson radical of $B$.

Suppose $C(1-aX)\ne C$. By Zorn's lemma, there exists a maximal
left ideal of $C$ containing $C(1-aX)$, say $L$. Put $M:=C/L$, a
simple $C$-module, and let $m_0$ denote the image of $1\in C$ in
$M$. Then $m_0\ne 0$ and $(1-aX)m_0=0$. Let $x$ and $a_M$ denote
the images of $X$ and $a$ in $\End M$. Since $X\in Z(C)$, we have
$x\in\End_C M$. Since $\gr\, C$ is finitely generated and
commutative, $x$ is invertible in $\End_C M$ and algebraic over
$\k$, by Quillen's lemma. Put $y:=x^{-1}$. Then $x=p(y)$ for some
$p\in \k[X]$, and $a_M(m_0)=y(m_0)$. Therefore,
$(1-ap(a))m_0=(1-yp(y))(m_0)=0$. On the other hand, $ap(a)\in
J(B)$, hence $1-ap(a)$ is invertible in $C$; see [\cite{Dix},
(3.1.12)] for instance. This contradiction shows that $C(1-aX)=C$.

As a consequence, $(a_0+a_1X+\cdots a_nX^n)(1-aX)=1$ for some
$a_i\in B$. Easy induction on $i$ gives $a_i=a^i$ for $0\le i\le
n$. Then $a^{n+1}=0$, showing that all elements in $J(B)$ are
nilpotent. As $(0)$ is a prime ideal of $B=H/I$, it follows from
[\cite{Dix}, (3.1.14)] that
$\cap_{\,J\in\text{Prim}\,B}\,J=J(B)=0$. But then $I$ is the
intersection of some primitive ideals of $H$, as claimed.

If $I\not\in\text{Prim}\,H$, then all primitive ideals of $H$
containing $I$ lie in $\frak I$. However, this is impossible as
$J_0\supsetneq I$; see part~(c). Thus $I$ must be a primitive
ideal of $H$ showing that the map $I\mapsto\widetilde{I}$ induces
a bijection between $\text{Prim}\, H$ and ${\mathcal X}_{\rm
inf}$, call it $\varkappa$.

Let $\cal Y$ be a closed set in $\text{Prim}\,H$. Then there is a
two-sided ideal $I_{\cal Y}$ of $H$ such that ${\cal
Y}=\{J\in\text{Prim}\,H\,|\,J\supseteq I_{\cal Y}\}$; see
[\cite{Dix}, (3.2.3)]. Our earlier remarks in the proof show that
${\varkappa}({\cal Y})=\,\{{\cal I}\in{\mathcal X}_{\rm inf}\,|\,
\,{\cal I}\supseteq\widetilde{I}_{\cal Y}\}$. Therefore,
$\varkappa\colon\,\text{Prim}\,H\rightarrow\,{\mathcal X}_{\rm
inf}$ is a closed map.

 Recall that the topology on ${\mathcal
X}_{\rm inf}$ is induced by the Jacobson topology on
$\text{Prim}\,U(\g)$. From [\cite{Dix}, (3.1.10)] it follows that
$\text{Prim}\,U(\g)$ is a Zariski space, that is any closed set in
$\text{Prim}\,U(\g)$ is a finite union of irreducible closed sets.
But then ${\mathcal X}_{\rm inf}$ is a Zariski space as well. Let
$\widetilde{\cal Y}$ be an irreducible closed set in ${\mathcal
X}_{\rm inf}$. Then there is an ${\cal I}\in\text{Spec}\,U(\g)$
such that $\widetilde{\cal Y}=\,\{{\cal J}\in {\mathcal X}_{\rm
inf}\,|\,\,{\cal J}\supseteq\cal I\}$; see [\cite{Dix}, (3.2.5)].
By Theorem~\ref{main1}(ii), ${\cal I}=\widetilde{I}$ for some
$I\in\text{Spec}\,H$. Furthermore, $\varkappa^{-1}(\widetilde{\cal
Y})=\,\{J\in\text{Prim}\,H\,|\,J\supseteq I\}$ by our remarks
earlier in  the proof. From this it is immediate that
$\varkappa^{-1}\colon\,{\mathcal X}_{\rm
inf}\rightarrow\,\text{Prim}\,H$ is a closed map too, proving (1).

\smallskip

\noindent (e) Let $V$ be a finite dimensional irreducible
$H$-module and $I=\text{Ann}_H\, V$, a primitive ideal of finite
codimension in $H$. By part~(b) of this proof,
$\widetilde{I}=\text{Ann}_{U(\g)}\,\widetilde{V}$ where
$\widetilde{V}=Q_\chi\otimes_H V$. As before, we identify
$\widetilde{V}$ with $\k[h,z_1,\ldots,z_s]\otimes V$; see (5.1).
Then $\widetilde{\mathcal B}_0:=U(\g)_f/\widetilde{I}_f$
identifies with a subalgebra of $\widetilde{\mathcal B}:={\mathcal
A}_e\otimes(H/I)$. More precisely, from part~(b) of the proof of
Theorem~\ref{main1} we know that $\tau\otimes\sigma$ acts on
$\widetilde{\mathcal B}$, and $\widetilde{\mathcal
B}_0=\widetilde{\mathcal B}^{\tau\otimes\sigma}$.

Since $U(\g)$ is a completely reducible $\ad\,\g$-module,
$(U(\g)/\widetilde{I})^{\n_\chi}\cong\,
U(\g)^{\n_\chi}/\widetilde{I}^{\n_\chi}$; see [\cite{Ja}, (3.2)]
for example. Since $f$ is central in $U(\g)^{\n_\chi}$, we then
have
$${\mathcal B}_0:=\,\widetilde{\mathcal B}_0^{\n_\chi}\cong\,
\big((U(\g)/\widetilde{I})_{\bar{f}}\big)^{\n_\chi}\cong\,
(U(\g)^{\n_\chi})_f/(U(\g)^{\n_\chi})_f\cap\widetilde{I}_f,$$
where $\bar{f}$ is the image of $f$ in $U(\g)/\widetilde{I}$. Put
$\bar{H}=H/I$ and $\mathcal{B}={\widetilde{\mathcal
B}}^{\n_\chi}$. Corollary~\ref{S5}(iii) implies that
$\widetilde{\mathcal B}=\widetilde{\mathcal
B}_0\oplus\widetilde{\mathcal B}_0\,t$ where $t$ stands for the
image of $\widehat{\Delta}$ in $\widetilde{\mathcal B}$. Since
$\widehat{\Delta}$ commutes with $\n_\chi$, it must be that
${\mathcal B}={\mathcal B}_0\oplus {\mathcal B}_0\,t$. In
conjunction with Corollary~\ref{S5}(v) this shows that the natural
map $({\mathcal A}_e\otimes H)^{\n_\chi}\rightarrow {\mathcal B}$
is surjective and ${\mathcal B}=\, \k[t,t^{-1}]\otimes \bar{H}$ as
algebras, where $\k[t,t^{-1}]$ is the Laurent polynomial ring in
$t$ over $\k$. As $\bar{H}\cong \,\End_\k V$ is a prime ring, so
too is ${\mathcal B}=\,\bar{H}[t,t^{-1}]$.

Since $\widetilde{I}$ is a prime ideal of $U(\g)$, the ring
$(U(\g)/\widetilde{I})^{\n_\chi}$ is prime with
$\text{rk}(U(\g)/I)^{\n_\chi}=\,\text{rk}(U(\g)/\widetilde{I})$;
see [\cite{Ja}, (13.10)]. Applying [\cite{Dix}, (3.6.15)] we
derive that the ring ${\mathcal
B}_0\,=\,\big((U(\g)/I)^{\n_\chi}\big)_{\bar{f}}$ is prime  with
$\text{rk}({\mathcal B}_0)\,=\,\text{rk}(U(\g)/\widetilde{I})$. On
the other hand, it follows from [\cite{Sh}] that
$\text{rk}\big(\k[t]\otimes\bar{H}\big)\,=\,\text{rk}(\bar{H})$.
As $\k[t,t^{-1}]\otimes\bar{H}$ is a localisation of
$\k[t]\otimes\bar{H}$ with respect to the Ore set
$\{t^i\,|\,i\in\Z_+\}$, we have that
$$\dim V\,=\,\text{rk}(\bar{H})\,=\,\text{rk}\big(\k[t,t^{-1}]\otimes\bar{H}\big)\,=
\text{rk}({\mathcal B});$$ see [\cite{McR}, Lemma~2.2.12] for
example. We thus need to show that $\text{rk}({\mathcal
B})\,=\,\text{rk}({\mathcal B}_0)$.

It should be mentioned that our present setting (involving a
quadratic extension of rings) resembles that of [\cite{J-pt},
(6.5)] where the above equality has been claimed in a more general
situation. However, the proof of Lemma~6.5 in [\cite{J-pt}] is
based on a faulty argument: in the notation of {\it op.\,cit.}
$J'$ is not a left ideal of $M_n\otimes L'$.\footnote{One hopes
that the proof of Lemma~6.5 in [\cite{J-pt}] can be corrected; see
[\cite{J-wp}]. When this is achieved, one would be able to
establish the quality
$\text{rk}(U(\g)/\widetilde{I})=\text{rk}(H/I)$ for {\it all}
primitive ideals of $H$.}  We are lucky here because in the
present case a different argument can be used to establish the
required equality of Goldie ranks.

Put $K:=\k(t),\,$ $K_0:=\k(t^2)$, and $S:=\k[t^2]\setminus\{0\}$,
a central Ore set in ${\mathcal B}_0$ and ${\mathcal B}$. Put
$B:=S^{-1}{\mathcal B}$ and $B_0:=S^{-1}{\mathcal B}_0$. Since
${\mathcal B}=\,\bar{H}[t,t^{-1}]$ and $S^{-1} \k[t,t^{-1}]=K$, we
have that $B\cong\,K\otimes_\k \bar{H}\,\cong\,\text{Mat}_n(K)$
where $n=\dim V$. Since $B$ is a simple Artinian ring, all regular
elements of $B$ are invertible. Since $S$ consists of regular
elements of $\mathcal B$, the universality property of quotient
rings yields $B\cong\,{\mathcal Q}({\mathcal B})$. In particular,
$K$ is the Goldie field of $\mathcal B$ (this argument  provides
another proof for the equality ${\mathrm rk}({\mathcal B})=\,\dim
V$.) As in [\cite{J-pt}, (6.5)], we regard $B$ as a Galois
extension of $B_0$. The involution $\tau\otimes\sigma$ induces a
$K_0$-automorphism of $B$, call it $\iota$. It is easy to see that
$\iota(t)=-t,\,$ $B=B_0\oplus B_0\, t,\,$ and $B_0=B^\iota$. Thus
$\iota$ can be viewed as the generator of the Galois group
$\text{Gal}(K/K_0)$.

Note that $Z(B_0)=Z(B)^\iota=K^\iota=K_0$ and $B\cong
B_0\otimes_{K_0}K$ as $K$-algebras. Since ${\mathcal B}_0$ is a
prime ring, so too is $B_0$; see [\cite{Dix}, (3.6.15)]. The
preceding remark then shows that $B_0$ is a simple algebra finite
dimensional over its centre $K_0$. Since $\k$ is algebraically
closed, $K_0\cong\k(t)$ is a $C_1$-field, by Tsen's theorem.
Therefore, $B_0\cong\,\text{Mat}_m(K_0)$ as $K_0$-algebras. As in
the previous paragraph one observes that $B_0\cong\,{\mathcal
Q}({\mathcal B}_0)$ and $K_0$ is the Goldie field of ${\mathcal
B}_0$. Since $B\cong\,B_0\otimes_{K_0}K$, one has $m=n$, proving
(2).

\smallskip

\noindent (f) Let $\cal I$ be a primitive ideal of $U(\g)$ with
${\mathcal VA}({\cal I})=\overline{{\cal O}}_{\rm min}$. Then
${\cal I}=\widetilde{I}$ for some $I\in\text{Prim}\,H$, by
part~(1) of this theorem. Thanks to Theorem~\ref{main1}(ii),
$\text{Dim}(H/I)=0$. Hence $H/I$ is finite dimensional over $k$;
see [\cite{McR}, (8.1.17)] for example. Since $H/I$ is a prime
ring, it must be that $H/I\cong\,\End(E)$ for some finite
dimensional $H$-module $E$. As $H/I$ is simple, $I=\,{\rm Ann}_H
E$. But then ${\cal I}=\,{\rm Ann}_{U(\g)}\big(Q_\chi\otimes_H
E\big)$ by part~(b) of this proof, as stated in (3).

Now let $V_1$ and $V_2$ be two finite dimensional irreducible
$H$-modules, and set $I_i:={\rm Ann}_H V_i$, $\,i=1,2$. If
$V_1\cong V_2$ as $H$-modules then, of course, $I_1=I_2$. Hence
$\widetilde{I}_1=\widetilde{I}_2$ yielding ${\rm
Ann}_{U(\g)}\big(Q_\chi\otimes_H V_1\big)\,=\,{\rm
Ann}_{U(\g)}\big(Q_\chi\otimes_H V_2\big)$, again by part~(b).
Conversely, if ${\rm Ann}_{U(\g)}\big(Q_\chi\otimes_H
V_1\big)\,=\,{\rm Ann}_{U(\g)}\big(Q_\chi\otimes_H V_2\big)$, then
$I_1=I_2$ in view of Theorem~\ref{main1}(ii) and part~(b). So
$H/I_1=\,H/I_2\,\cong\, \text{Mat}_m(\k)$ for some $m$. It is now
straightforward to see that $V_1\cong V_2$ as $H$-modules, giving
(4).

Fix an algebra homomorphism $\eta\colon\,Z(\g)\rightarrow\,\k$ and
identify $Z(\g)$ with $Z(H)$; see Corollary~\ref{S5}(vi). If $V$
is a finite dimensional $H$-module with central character $\eta$
and $I={\rm Ann}_H V$, then $I\cap Z(H)=\,\text{Ker}\,\eta$ is a
maximal ideal of $Z(H)$, by Schur's lemma. Hence
$\text{Ker}\,\eta\,=\,{\rm Ann}_{U(\g)}\big(Q_\chi\otimes_H
V\big)\cap Z(\g)=\,\widetilde{I}\cap Z(\g)$. Thanks to
Theorem~\ref{main1}(ii) we also have ${\mathcal
VA}(\widetilde{I})=\,\overline{\cal O}_{\rm min}$. Now let ${\cal
I}\in {\mathcal X}$ be such that $Z(\g)\cap {\cal
I}=\,\text{Ker}\,\eta$ and ${\mathcal VA}({\cal
I})=\,\overline{\cal O}_{\rm min}$. By parts~(3) and (4) of this
theorem, ${\cal I}=\,{\rm Ann}_{U(\g)}\big(Q_\chi\otimes_H E\big)$
for some finite dimensional irreducible $H$-module $E$, which is
uniquely determined up to isomorphism. Since
$\text{Ker}\,\eta\subset {\cal I}$ and $Z(\g)=Z(H)$, the
$H$-module $E$ has central character $\eta$. We obtain (5).

\smallskip

\noindent (g) Let $I\in\text{Spec}\,H$ and suppose that $I\cap
Z(H)$ is a maximal ideal of $Z(H)$. By Theorem~\ref{main1}(ii),
$\widetilde{I}\in\text{Spec}\,U(\g)$. As explained in the proof of
Corollary~\ref{S5}(vi),
$$Z(H)= Z({\mathcal A}_e\otimes H)=\,\k\otimes\big(H^\sigma\cap Z(H)\big)=
\big(Z({\mathcal A}_e\otimes H)\big)^{\tau\otimes\sigma}\subseteq
Z(U(\g))_f=\,Z(\g).$$ It follows that $I\cap Z(H)\subseteq\,
(\k\otimes I_+)\cap Z(\g)\subseteq \widetilde{I}_f\cap
Z(\g)=\widetilde{I}\cap Z(\g).$ As $Z(\g)=Z(H)$, we deduce that
$\widetilde{I}\cap Z(\g)$ is a maximal ideal of $Z(\g)$. But then
$\widetilde{I}\in{\mathcal X}$; see [\cite{Dix}, (8.5.7)]. Since
$\widetilde{I}$ has infinite codimension in $U(\g)$, part~(1) of
this theorem implies that
$I=\varkappa^{-1}(\widetilde{I})\in\text{Prim}\,H$. Finally,
suppose $I\in\text{Prim}\,H$ and let $M$ be an irreducible
$H$-module such that $I=\text{Ann}_H M$. Since $\gr\,H$ is
finitely generated and commutative, Quillen's lemma shows that
$Z(H)$ acts on $M$ as scalar operators. Consequently, $I\cap Z(H)$
is a maximal ideal of $Z(H)$. The proof of the theorem is now
complete.
\end{pf}

\subsection{} Before finishing this section off we wish to discuss a possible
extension of the above results to the case of a general algebra
$H_\chi$. Let $e$ be any nilpotent element in $\g$ and
$\chi=\chi_e\in\g^*$. Put ${\cal O}={\cal O}_\chi$ and denote by
${\mathcal X}_{\cal O}$ the set of all primitive ideals $\cal I$
of $U(\g)$ with ${\mathcal VA}({\cal I})\supset \cal O$. Take
$\text{Prim}\,H_\chi$ with the Jacobson topology and ${\mathcal
X}_{\cal O}$ with the topology induced by that of $\mathcal X$.
\begin{question} Are the following true?
\begin{itemize}
\item[1.] The centre of $H_\chi$ coincides with the image of
$Z(\g)$ in $H_\chi$.\footnote{At the Oberwolfach meeting on
enveloping algebras in March 2005 Victor Ginzburg has explained to
me that this is a consequence of the finiteness of the number of
symplectic leaves of ${\mathcal S}_e$ contained in the fibres of
the morphism $f\colon \,{\mathcal S}_e\rightarrow \g/\!\!/G$
iduced by the adjoint quotient map of $\g$.  Each homogeneous
element $z\in\text{gr}\,Z(H_\chi)$ lies in the Poisson centre of
$\text{gr}\,H_\chi=\,\k[{\mathcal S}_e]$, hence reduces to scalars
on all symplectic leaves of ${\mathcal S}_e$. The Poisson
structure on ${\mathcal S}_e$ induced by multiplication in
$H_\chi$ is determined in [\cite{GG}, (3.2)]. By [\cite{P02},
(5.4), (6.3)], all scheme-theoretic fibres of $f$ are reduced and
irreducible, and $\text{gr}\,H_\chi$ is a flat module over
$\text{gr}\,Z(\g)$. These results are needed to carry out
Ginzburg's argument: Since each fibre of $f$ contains a Zariski
dense symplectic leaf of ${\mathcal S}_e$, the regular function
$z$ is constant on each fibre of $f$. The flatness of the
$\text{gr}\,Z(\g)$-module $\text{gr}\,H_\chi$ along with the fact
that all scheme-theoretic fibres of $f$ are reduced then yields
$z\in\text{gr}\,Z(\g)$ implying $Z(H_\chi)=Z(\g)$.}

\medskip

\item[2.] There exists a homeomorphism $\varkappa\colon\,
\text{Prim}\,H_\chi\,\longrightarrow\,{\mathcal X}_{\cal O}$ such
that:

\medskip

\begin{enumerate}
\item[(a)] ${\rm Dim}(U(\g)/\varkappa(I))\,=\,{\rm Dim}(H/I)+\dim
{\cal O}\qquad\big(\forall\,I\in{\mathrm Prim}\,H_\chi\big);$

\smallskip

\item[(b)] ${\rm rk}(U(\g)/\varkappa(I))=\sqrt{\dim_\k H/I}\,\,$
for all $I\in{\rm Prim}\,H_\chi$ with $\codim_{H_\chi}\!
I<\infty$.
\end{enumerate}
\medskip

\item[3.] For every character $\eta$ of $Z(H_\chi)=Z(\g)$ the map
$\varkappa$ induces a bijection between the isoclasses of finite
dimensional $H_\chi$-modules with central character $\eta$ and the
primitive ideals of $U(\g)$ with ${\cal I}\cap
Z(\g)=\,\text{Ker}\,\eta$ and ${\mathcal VA}({\cal
I})=\overline{{\cal O}}$.
\end{itemize}
\end{question}

\section{\bf The Joseph ideal and a presentation of $H$}
\subsection{} In his seminal work [\cite{J}] Joseph has
discovered that outside type $\rm A$ the enveloping algebra
$U(\g)$ has a unique completely prime primitive ideal whose
associated variety is $\overline{\cal O}_{\rm min}$. This ideal is
denoted ${\cal J}_0$ and referred to as the {\it Joseph ideal} of
$U(\g)$.  For $\g$ of type $\rm A$ the completely prime primitive
ideals of $U(\g)$ with ${\mathcal VA}({\cal I})=\overline{{\cal
O}}_{\rm min}$ form a single family parametrised by the elements
of $\k$ (this will be explained in more detail in the course of
proving Theorem~\ref{main3}).

The Joseph ideal is prominent in several areas of representation
theory, especially in the theory of minimal representations of
$p$-adic groups. Different realisations of ${\cal J}_0$ can be
found in the literature for various types of $\g$ but most of them
are ad hoc. This seems almost inevitable as outside type $\rm A$
the orbit ${\cal O}_{\rm min}$ is ${\it rigid}$, that is forms a
single sheet in $\g^*$. Hence ${\cal J}_0$ cannot be obtained by
parabolic induction from a primitive ideal of a proper Levi
subalgebra of $\g$, the only `regular' way so far to obtain
primitive ideals.

It was noticed by Savin (in a letter to Vogan) that Joseph's
original proof of the uniqueness of ${\cal J}_0$ was incomplete.
This was recently fixed by W.T.~Gan and Savin with the assistance
of Wallach; see [\cite{GS}]. The argument in [\cite{GS}] relies on
some invariant theory and earlier results of Garfinkle. We shall
see in a moment that the existence and uniqueness of ${\cal J}_0$
follow readily from our results; see also Remark~\ref{joseph}.

\subsection{} Retain the
assumptions and conventions of Sections~4 and 5. Set
$\k_0:=H/H^+$. Since $\k_0$ is an irreducible $H$-module, so is
the $\g$-module $Q_{\chi,0}:=\,Q_\chi\otimes_{H}\k_0$. So
$J_0:=\text{Ann}_{U(\g)}\,Q_{\chi,0}$ is a primitive ideal of
$U(\g)$.
\begin{prop}\label{P3} The ideal
$J_0$ is completely prime and ${\mathcal
VA}\big(J_0)=\overline{\cal O}_{\mathrm min}$. If $\g$ is not of
type $\mathrm A$ then $J_0$ is the only primitive ideal of $U(\g)$
with these properties, and hence $J_0$ is the Joseph ideal in this
case.
\end{prop}
\begin{pf}
Theorem~\ref{main2}(2) shows that
$\,\text{rk}(U(\g)/J_0)=\,\dim_\k \k_0=1$. Hence ${\mathcal
Q}\big(U(\g)/J_0\big)$ is a division ring. Then $U(\g)/J_0$ is a
domain, that is $J_0$ is completely prime. Theorem~\ref{S3} gives
${\mathcal VA}(J_0)=\overline{\cal O}_{\rm min}$. Now suppose $\g$
is not of type $\rm A$. Then Corollary~\ref{C1} implies that $H$
has a unique one-dimensional representation. In view of
Theorem~\ref{main2} this means that $U(\g)$ has a {\it unique}
completely prime primitive ideal whose associated variety is
$\overline{\cal O}_{\rm min}$. So $J_0={\cal J}_0$ in this case.
\end{pf}

\begin{rem}
The existence part of our proof is hardly shorter than Joseph's
original proof in [\cite{J}] as it relies on the brute force
computations of Section~4. However, there is a slightly different
proof of the uniqueness of ${\cal J}_0$ which eludes Section~4
completely. We just sketch the argument leaving the details to the
interested reader: If $I$ is an ideal of codimension $1$ in $H$
then $[H,H]\subset I$. Since outside type $\rm A$ the Lie algebra
$\z_\chi(0)$ is semisimple, we have $\Theta(\z_\chi(i))\subset I$
for $i=1,2$. Also, $C-\mu\in I$ for some $\mu\in \k$. Using
Proposition~\ref{P1} it is not hard to observe that $\mu$ is
independent of $I$. Therefore, $H$ cannot afford more than one
ideal of codimension $1$. The rest of the proof is unchanged.
\end{rem}
\subsection{} In [\cite{J}], Joseph has also computed the
{\it infinitesimal character} of ${\cal J}_0$, that is the algebra
homomorphism $Z(\g)\rightarrow\k$ through which the centre $Z(\g)$
acts an the primitive quotient $U(\g)/{\cal J}_0$. We are going to
use his result to determine the elusive constant $c_0$ from
Proposition~\ref{P2}.
\begin{theorem}\label{main3}
In the notation of Section~4, the algebra $H$ is generated by the
Casimir element $C$ and the subspaces $\Theta(\z_\chi(i))$ for
$i=0,1$, subject to the following relations:

\smallskip

\begin{enumerate}
\item[(i)\,] $[\Theta_x,\Theta_y]=\Theta_{[x,y]}$ for all
$x,y\in\z_\chi(0)$;

\medskip

\item[(ii)\,] $[\Theta_x,\Theta_u]=\Theta_{[x,u]}$ for all
$x\in\z_\chi(0)$ and $u\in\z_\chi(1)$;

\medskip

\item[(iii)\,] $C$ is central in $H$;

\medskip

\item[(iv)\,]
$[\Theta_u,\Theta_v]\,=\,\frac{1}{2}(f,[u,v])\big(C-\Theta_{\rm
Cas}-c_0\big)+
\frac{1}{2}\sum_{i=1}^{2s}\,\big(\Theta_{[u,z_i]^\sharp}\,\Theta_{[v,z_i^*]^\sharp}+
\Theta_{[v,z_i^*]^\sharp}\,\Theta_{[u,z_i]^\sharp}\big)$

\medskip

\noindent for all $u,v\in\z_\chi(1)$, where $\Theta_{\rm
Cas}=\sum_i\Theta_{a_i}\Theta_{b_i}$ is a Casimir element
 of the Lie algebra
$\Theta(\z_\chi(0))$ and the constant $-c_0$ is given in the table
below:
\end{enumerate}

$$\begin{tabular}{|c|c|c|c|c|c|c|c|c|c|}
\hline ${\mathrm Type}$ & ${\mathrm A}_n$ & ${\mathrm B}_n$ &
${\mathrm C}_n$ & ${\mathrm D}_n$ & ${\mathrm E}_6$ & ${\mathrm
E}_7$ & ${\mathrm E}_8$ & ${\mathrm F}_4$ & ${\mathrm
G}_2$\\
\hline  $-c_0$ & $\frac{n(n+1)}{4}$ & $\frac{(2n+1)(2n-3)}{4}$ &
$\frac{n(2n+1)}{8}$ & $n(n-2)$ & $36$ &
$84$ & $240$ & $\frac{39}{2}$ & $\frac{28}{9}$\\
\hline
\end{tabular}$$

\bigskip

\smallskip

\noindent If $\g$ is not of type $\rm A$ then $c_0$ is the
eigenvalue of $C$ on the primitive quotient $U(\g)/{\cal J}_0$. If
$\g$ is of type ${\rm A}_1$ then $\z_\chi(0)=\,\z_\chi(1)= 0$ and
$H=\k[C]$.
\end{theorem}
\begin{pf} First we determine $c_0$. Recall from (4.1) that
$(\gamma,\gamma)=2$. Therefore, if $\g$ is not of type ${\mathrm
C}_n$ or ${\rm G}_2$, then the scalar product
$(\,\cdot\,,\,\cdot\,)$ on the ${\mathbb Q}$-span of $P$ in $\h^*$
coincides with the scalar product $(\,\cdot\,|\,\cdot\,)$ from
Bourbaki's tables in [\cite{Bou}]. In the remaining two cases,
$(\,\cdot\,,\,\cdot\,)=\frac{1}{2} (\,\cdot\,|\,\cdot\,)$ for $\g$
of type ${\rm C}_n$, and $(\,\cdot\,,\,\cdot\,)=\frac{1}{3}
(\,\cdot\,|\,\cdot\,)$ for $\g$ of type ${\rm G}_2$. Recall that
for any $\lambda\in \h^*$ the eigenvalue of the Casimir element
$C$ on the irreducible highest weight module $L(\lambda)$ equals
$(\lambda,\lambda+2\rho).$

\smallskip

\noindent (a) Suppose $\g$ is not of type $\rm A$. In [\cite{J},
p.~15], Joseph has found an irreducible highest weight module
$L(\lambda_0)$ with ${\rm Ann}_{U(\g)}\, L(\lambda_0)={\cal J}_0$.
It is immediate from the definition of $H^+$ that $C$ acts on
$\k_0=H/H^+$ as scalar $c_0$. But then $C_{\vert
Q_0}=\,c_0\,\text{id}$. Proposition~\ref{P3} now shows that $C$
acts as $c_0\,\text{id}$ on the primitive quotient $U(\g)/{\cal
J}_0$. In view of our remarks above this yields
$c_0=(\lambda_0,\lambda_0+2\rho)$.

If $\g$ is of type $\rm E$, then $\lambda_0=-\varpi_4$; see
[\cite{J}]. Using parts (VI) and (VII) of Tables~V--VII in
[\cite{Bou}] one finds out that $c_0=-240$ for $\g$ of type ${\rm
E}_8$, $c_0=-84$ for $\g$ of type ${\rm E}_7$, and $c_0=-36$ for
$\g$ of type ${\rm E}_6$. If $\g$ is of type ${\rm D}_n$, $n\ge
4$, then $\lambda_0=-\varpi_{n-2}$. Parts~(VI) and (VII) of
Table~V in [\cite{Bou}] yield $(\varpi_{n-2},\varpi_{n-2})=n-2$
and $(\varpi_{n-2},2\rho)=n^2-n-2$. Therefore, $c_0=-n(n-2)$ in
this case.

If $\g$ is of type ${\rm B}_n$, $n\ge 3$, then
$\lambda_0=-\frac{1}{2}(\varpi_{n-2}+\varpi_{n-1})$. Form Table~II
in [\cite{Bou}] we get $(\lambda_0,\lambda_0)=\frac{4n-7}{4}$ and
$(\lambda_0,2\rho)=-n^2+\frac{5}{2}$. Then
$c_0=-(n^2-n-\frac{3}{4})=-\frac{(2n+1)(2n-3)}{4}.$ If $\g$ is of
type ${\rm C}_n$, $n\ge 2$, then $\lambda_0=-\frac{1}{2}\varpi_n$.
Table~III in [\cite{Bou}] yields
$(\lambda_0|\lambda_0)=\frac{n}{4}$ and
$(\lambda_0|\,\!2\rho)=-\frac{n(n+1)}{2}$. Consequently,
$(\lambda_0|\lambda_0+2\rho)=-\frac{n(2n+1)}{4}$ and
$c_0=\,\frac{1}{2}(\lambda_0|\lambda_0+2\rho)=-\frac{n(2n+1)}{8}.$
If $\g$ is of type ${\rm G}_2$, then
$\lambda_0=-\frac{2}{3}\varpi_2=-\frac{2}{3}\,\widetilde{\alpha}$.
Hence $(\lambda_0|\lambda_0)=\frac{8}{3}$ and
$$(\lambda_0|2\rho)=(-\frac{2}{3}\varpi_2|10\alpha_1+6\alpha_2)=-4(\varpi_2|\alpha_2)=
-2(\alpha_2|\alpha_2)\la\varpi_2,\alpha_2\ra=-12;$$ see
[\cite{Bou}, Table~IX]. Therefore,
$(\lambda_0|\mu_0+2\rho)=-\frac{28}{3}$ and
$c_0=\frac{1}{3}(\lambda_0|\lambda_0+2\rho)=-\frac{28}{9}.$ If
$\g$ is of type ${\rm F}_4$, then
$\lambda_0=-\frac{1}{2}(\varpi_1+\varpi_2)$. Using Table~VIII in
[\cite{Bou}] we get $(\lambda_0|\lambda_0)=\frac{7}{2}$ and
$(\lambda_0|2\rho)=-23.$ Hence
$c_0=(\lambda_0,\lambda_0+2\rho)=(\lambda_0|\lambda_0+2\rho)=-\frac{39}{2}.$

\smallskip

\noindent (b) Now suppose $\g$ is of type ${\rm A}_n$, $n\ge 2$.
This case is more subtle because here we have an infinite family
of completely prime ideals in $\mathcal X$ sharing the same
associated variety $\overline{\cal O}_{\rm min}$. In order to
determine $c_0$ in the present case we shall have to locate a
special member of this family. Theorem~\ref{main2} will play a
crucial r{\^o}le here.

Let $\p_1$ be the standard parabolic subalgebra of $\g$ whose Levi
subalgebra is generated by $\h$ and all $e_{\pm \alpha_i}$ with
$2\le i\le n$. For any $t\in\k$ the linear function $t\varpi_1$
vanishes on all $h_{\alpha_i}$ with $i\ge 2$, hence extends
uniquely to a one-dimensional representation of $\p_1$. Let $\k_t$
denote the corresponding one-dimensional $\p_1$-module, and put
$I(\p_1,t):={\rm Ann}_{U(\g)}\big(U(\g)
\otimes_{U(\p_1)}\k_t\big)$. Although some of the induced
$\g$-modules $U(\g)\otimes_{U(\p_1)}\k_t$ are reducible, it
follows from [\cite{Dix}, (8.5.7)] and Conze's theorem [\cite{C}]
that all two-sided ideals $I(\p_1,t)$ are primitive and completely
prime. It is not hard to check that $t_1\varpi_1$ and
$t_2\varpi_1$ are conjugate under the dot action of the Weyl group
$W\cong{\frak S}_{n+1}$ if and only if $t_1=t_2$. Consequently,
all members of the family ${\frak
I}_\k:=\{I(\p_1,t)\,|\,\,t\in\k\}$ have pairwise distinct
infinitesimal characters. In view of [\cite{Ja}, (17.17)] they
share the same associated variety $\overline{\cal O}_{\rm min}$.

Let ${\cal I}\in{\mathcal X}$ be a completely prime ideal with
${\mathcal VA}({\cal I})=\overline{\cal O}_{\rm min}$. It follows
from the main result of M{\oe}glin in [\cite{Mo}] that there exist
a standard parabolic subalgebra $\p$ of $\g$ and a one-dimensional
representation $f\colon\,\p\rightarrow \k$ such that ${\cal
I}={\rm Ann}_{U(\g)}\big( U(\g)\otimes_{U(\p)} \k_f\big)$. In
conjunction with [\cite{Ja}, (17.17), (15.27)] this yields ${\cal
I}\in{\frak I}_\k$.

Let $I_0$ be the two-sided ideal of $H$ generated by $[H,H]$. In
order to describe the one-dimensional representations of $H$ we
have to take a close look at the commutative $\k$-algebra $H^{\rm
ab}:=H/I_0$. Given $x\in H$ we denote by $\bar{x}$ the image of
$x$ in $H^{\rm ab}$. We may assume that $\g={\mathfrak
sl}_{n+1}(\k)$ and $\h$ is the subalgebra of all diagonal matrices
in $\g$. Let $\{e_{ij}\,|\,1\le i,j\le n+1\}$ be the matrix units
in ${\mathfrak gl}_{n+1}(\k)$. We may also assume that
$\alpha_i=\varepsilon_i-\varepsilon_{i+1}$ and
$e_{\alpha_i}=e_{i,i+1}$ for $1\le i\le n$; see Table~I in
[\cite{Bou}]. Then $e=e_{n,n+1}$, $h=e_{nn}-e_{n+1,n+1},$ and
$f=e_{n+1,n}$ by our conventions in (4.1). No generality will be
lost by assuming that $z_i=e_{n+1,i}$ and $z_i^*=-e_{in}$ for
$1\le i\le n-1$ (notice that $s=n-1$ in the present case). It is
straightforward to see that the centre of the subalgebra
$\g(0)^\sharp=\z_\chi(0)$ is spanned by the element
$z:=e_{nn}+e_{n+1,n+1}-\frac{2}{n+1}I_{n+1}$. Therefore,
$\z_\chi(0)=\k z\oplus \z_\chi(0)'$ where
$\z_\chi(0)'=[\z_\chi(0),\z_\chi(0)]$. Since $\z_\chi(1)$ has no
zero weight relative to $\h_e=\,\h\cap\z_\chi(0)$, this implies
that the $\k$-algebra $H^{\rm ab}$ is generated by $\bar{C}$ and
$\bar{\Theta}_z$.

Set $u=e_{1,n+1}$ and $v=-e_{n,1}$. We have $[u,z_i]^\sharp\in
[\z_\chi(0),\z_\chi(0)]$ for $2\le i\le s$ and $[u,z_i^*]=0$ for
$1\le i\le s$. Also,
$[u,z_1]^\sharp=(e_{11}-e_{n+1,n+1})-\frac{1}{2}h=e_{11}-
\frac{1}{2}(e_{nn}+e_{n+1,n+1})$. Likewise,
$[v,z_1^*]^\sharp=(e_{nn}-e_{11})-\frac{1}{2}h=-e_{11}+
\frac{1}{2}(e_{nn}+e_{n+1,n+1})$. Since
$(z,z)=2-\frac{8}{n+1}+\frac{4}{n+1}=\frac{2(n-1)}{n+1},$ we can
take $a_1=z$, $b_1=\frac{n+1}{2(n-1)}\,z,$ and $a_i,
b_i\in\z_\chi(0)'$ for $i>1$. Next observe that $(z,e_{11}-
\frac{1}{2}(e_{nn}+e_{n+1,n+1}))=-1$. As $z\perp \z_\chi(0)'$, it
follows that $e_{11}- \frac{1}{2}(e_{nn}+e_{n+1,n+1})$ is
congruent to $-\frac{1}{(z,z)}\,z$ modulo $\z_\chi(0)'$. As
$(f,[u,v])=1$, Proposition~\ref{P2} now yields
$$\frac{1}{2}\Big(\bar{C}-c_0-\frac{n+1}{2(n-1)}\,
\bar{\Theta}_z^2\Big)+\frac{1}{2}\cdot
(-2)\cdot\frac{(n+1)^2}{4(n-1)^2}\,\bar{\Theta}_z^2\,=\,0.$$ As a
consequence, we obtain that the following relation holds in
$H^{\rm ab}$:
\begin{eqnarray}\label{eqS6}
\bar{\Theta}_z^2\,=\,\frac{(n-1)^2}{n(n+1)}\,(\bar{C}-c_0).
\end{eqnarray} Since $n\ge 2$, this shows that $H^{\rm ab}$ is a
homomorphic image of the polynomial algebra $\k[X]$. On the other
hand, Theorem~\ref{main2} in conjunction with our earlier remarks
entails that the one-dimensional representations of $H^{\rm ab}$
are in $1$-$1$ correspondence with the elements in ${\mathfrak
I}_\k$. But then the maximal spectrum ${\rm Max}(H^{\rm ab})$ of
$H^{\rm ab}$ is infinite, forcing $H^{\rm ab}\cong\,
\k[X,Y]/(X^2-\frac{(n-1)^2}{n(n+1)}\,Y)$ (under the algebra map
$X\mapsto \bar{\Theta}_z$, $\,Y\mapsto\bar{C}-c_0$).

For $c\in\k$ we let ${\rm Max}_c(H^{\rm ab})$ (resp. ${\frak
I}_{\k,\,c}$) denote the set of all $I$ in ${\rm Max}(H^{\rm ab})$
(resp. $I(\p_1,t)$ in ${\frak I}_\k$) containing $\bar{C}-c$
(resp. $C-c$). It is immediate from (\ref{eqS6}) that
\[
|{\mathrm Max}_c(H^{\rm ab})|\,=\, \left\{
\begin{array}{ll}
2 & \mbox{when $\,c\ne c_0$},\\
1 & \mbox{when $\,c=c_0$}.
\end{array}
\right.
\]
Theorem~\ref{main2} implies that for any $c\in\k$ the map
$\varkappa\colon\,{\rm Prim}\,H\rightarrow {\mathcal X}_{\rm
inf},\ I\mapsto\widetilde{I},$ induces a bijection between ${\rm
Max}_c(H^{\rm ab})$ and ${\frak I}_{\k,\,c}$. It is well-known
that $C$ acts on the induced module $U(\g)\otimes_{U(\p_1)}\k_t$
as $(t\varpi_1,t\varpi_1+2\rho)\,\text{id}$. So $I(\p_1,t)$
contains $C-(t\varpi_1,t\varpi_1+2\rho)$. It is immediate from
[\cite{Bou}, Table~I] that
$(t\varpi_1,t\varpi_1+2\rho)=\frac{n}{n+1}\,t^2+nt.$ The equation
$\frac{n}{n+1}\,t^2+nt-c=0$ has two distinct roots if and only if
$n^2+\frac{4nc}{n+1}\ne 0$. Therefore,
\[
|{\frak I}_{\k,c}|\,=\, \left\{
\begin{array}{ll}
2 & \mbox{when $\,c\ne -\frac{n(n+1)}{4}$},\\
1 & \mbox{when $\,c= -\frac{n(n+1)}{4}$}.
\end{array}
\right.
\]
But then ${\rm Max}_{c_0}(H^{\rm ab})$ must be mapped onto ${\frak
I}_{\k,\,-\frac{n(n+1)}{4}}$, forcing $c_0=-\frac{n(n+1)}{4}$.

\noindent (c) Now let $\widehat{H}$ be the associative
$\k$-algebra generated by an element $\widehat{C}$ and isomorphic
copies $\widehat{\Theta}(\z_\chi(i))$ of the subspaces
$\z_\chi(i)$ with $i=1,2$, subject to the relations (i)\! -\! (iv)
from the formulation of this theorem. Define an increasing
filtration in $\widehat{H}$ by giving $\widehat{C}$ filtration
degree $4$, by assigning to all nonzero elements of
$\widehat{\Theta}(\z_\chi(i))$ filtration degree $i+2$, and by
extending to $\widehat{H}$ algebraically.

Choose bases $x_1,\ldots,x_q$ and $y_1,\ldots,y_{2s}$ in
$\z_\chi(0)$ and $\z_\chi(1)$, respectively, and set $X_i=
\widehat{\Theta}_{x_i}$ for $1\le i\le q$ and
$Y_i=\widehat{\Theta}_{y_i}$ for $1\le i\le 2s$. Let
$\widehat{H}'$ be the $\k$-span of all monomials $m({\bf a},{\bf
b},l):=X_1^{a_1}\cdots X_q^{a_q}\cdot Y_1^{b_1}\cdots
Y_{2s}^{b_{2s}}\cdot \widehat{C}^l$ with $a_i,b_j,l\in\Z_+$. Note
that $m({\bf a},{\bf b},l)$ has filtration degree $2\sum a_i+3\sum
b_j+4l$. Using the relations (i)\! -\! (iv) and induction on the
filtration degree of $m({\bf a},{\bf b},l)$ it is straightforward
to see that $\widehat{H}'$ is a two-sided ideal of $\widehat{H}$.
Since $1\in\widehat{H}'$ it must be that
$\widehat{H}=\widehat{H}'$.

It follows from Proposition~\ref{P2} and Lemmas~\ref{L4} and
~\ref{L5} that there is a surjective algebra homomorphism
$f\colon\,\widehat{H}\to H$ such that$f(\widehat{C})=C$,
$f(X_i)=\Theta_{x_i}$ for $i\le q$, and $f(Y_i)=\Theta_{y_i}$ for
$i\le 2s$. Since the vectors $f(m({\bf a},{\bf b},l))$ are
linearly independent in $H$, by [\cite{P02}, Theorem~4.6(ii)], the
equality $\widehat{H}=\widehat{H}'$ shows that $f$ is injective.
But then $\widehat{H}\cong H$, and our proof is complete.
\end{pf}
\begin{rem}
We have originally computed the infinitesimal character of $J_0$
by using a direct approach in the spirit of Section~4; this was
done before we established a link between $H$ and ${\cal J}_0$.
Having established that link we discovered that outside type ${\rm
A}$ our result was consistent with [\cite{J}, p.~15]. In type
${\rm A}$ we have found two different proofs yielding the same
result. This eventually convinced us that the quadratic relation
of Theorem~\ref{main3} was correct. Our computations are rather
lengthy, especially in type ${\rm C}$, and will not be presented
here.
\end{rem}

\begin{rem}\label{Osc} Suppose $\g$ is of type ${\rm C}_2$. Then $\z_\chi(0)\cong {\frak sl}(2)$
and $\z_\chi(1)$ is an irreducible $2$-dimensional
$\z_\chi(0)$-module. In this case Theorem~\ref{main3} shows that
the algebra $H$ is generated by six elements $e,h,f,u,v,c,\,$
subject to the following relations:
\begin{itemize}

\item[1.] $(e,h,f)$ is an ${\frak sl}(2)$-triple relative to the
commutator product in $H$;

\item[2.] $[e,u]=0=[f,v],\ [e,v]=u,\ [f,u]=v,\ [h,u]=u,\
[h,v]=-v$;

\item[3.] $[u,v]=ef+fe+\frac{1}{2}h^2-\frac{1}{2}\,c-\frac{5}{8}$;

\item[4.] $c$ is central in $H$.
\end{itemize}
For any $t\in\k$ the factor algebra $H_t:=H/(c-t)$ is isomorphic
to one of the deformed oscillator algebras $H_f$ studied by Khare
in [\cite{Kh}] (in type ${\rm C}_n,$ $n\ge 3$, the defining
relations of the algebras $H/(C-t)$ differ from those of $H_f$).
Arguing as in the proof of Theorem~6.4 in [\cite{P02}] one can
observe that image of $Z(\g)$ in $H_t$ is isomorphic to a
polynomial algebra in one variable. It is likely that the centre
of $H_t$ is generated by that image. It would be very interesting
to describe the Goldie field of $H$ in the present case. In view
of Corollary~\ref{C5}(iv) this might help to resolve the
Gelfand--Kirillov conjecture for $\g={\frak sp}_{4}(\k)$.
\end{rem}

\begin{rem}\label{joseph}
(A. Joseph) Assume $\g$ is not of type $\mathrm A$. The argument
below gives a short proof of the uniqueness of ${\cal J}_0$
relying only on the information available at the time when
[\cite{J}] was written. Let $(e,h,f)$ be an $\sl_2$-triple in $\g$
with $e$ being a highest root vector. Let ${\mathfrak
d}=\text{Ker}\,(\ad h-\text{id})\oplus \text{Ker}\,(\ad
h-2\,\text{id})$, a Heisenberg Lie subalgebra of $\g$, and
${\mathfrak r}=\k h\oplus{\mathfrak d}.$ Let $J$ be a completely
prime primitive ideal of $U(\g)$ such that
$\text{Dim}(U(\g)/J)=\dim {\cal O}_{\rm min}$. Since $\ad e$ is
nilpotent, $U(\g)/J$ embeds into its localisation $U$ at $e$ which
contains the localisation $A$ of $U({\mathfrak r})$ at $e$. Let
$Z$ denote the centraliser of $A$ in $U$. It follows from
[\cite{J}, Lemma~4.1] that $A$ is a localised Weyl algebra with
$\text{Dim}(A)=\dim{\cal O}_{\rm min}$. Clearly, $Z$ inherits a
filtration from $U(\g)$ such that $\text{gr}\,Z$ is commutative.
Since $\mathfrak d$ is the nilradical of a parabolic subalgebra of
$\g$, Hadziev's theorem yields that the algebra $\text{gr}\,Z$ is
finitely generated (one also needs the fact that $\ad h$ is
semisimple). Since $A$ is central simple, the multiplication map
$Z\otimes A\rightarrow U$ is injective. Since both $\text{gr}\,Z$
and $\text{gr}\,A$ are commutative and finitely generated, we have
that $\text{Dim}(Z\otimes A)=\text{Dim}(Z)+\text{Dim}(A)$. Since
$\text{Dim}(U)=\text{Dim}(A)$, we get $\text{dim}(Z)=0$. Hence $Z$
is algebraic over $\k$. Since $U$ is a domain, we now obtain
$Z=\k$. Since $\mathfrak d$ consists of nilpotent elements of
$\g$, Taylor's lemma proved in [\cite{J'}] implies that for any
$h$-weight vector $u\in U$ there is an $h$-weight vector $a\in A$
of the same weight as $u$ such that $u-a$ commutes with the image
of $\mathfrak d$ in $U$ (a preprint version of [\cite{J'}] was
available since 1973 and is quoted in [\cite{J}]). Taking $u$ to
be the image of $f$ in $U$ we get $e(f-a)\in Z=\,\k$. But then
$f\in A$ and so $U=A$, by the simplicity of $\g$. Now apply
[\cite{J}, Theorem~4.3] to deduce the equality $J={\cal J}_0$.
\end{rem}

\begin{lemma}\label{L6}
If $\g$ is of type ${\rm G}_2$ then the algebra $H$ admits a
$2$-dimensional irreducible representation $\rho$ such that
$\rho(C)=-\frac{16}{9}\,{\rm id}$ and $\rho(\Theta_u)=0$ for all
$u\in\z_\chi(1)$.
\end{lemma}
\begin{pf}
 For any
$\alpha=m\alpha_1+n\alpha_2\in\Phi^+$ we set $e_{m,n}=e_\alpha,\,$
$h_{m,n}=h_\alpha,$ and $f_{m,n}=e_{-\alpha}$. Recall from (4.1)
that $\beta=\alpha_2$. It is easy to see that in the present case
$\z_\chi(0)\cong\sl(2)$ and $\z_\chi(1)$ is an irreducible
$4$-dimensional $\z_\chi(0)$-module. Furthermore, $\z_\chi(0)=\k
e_{2,1}\oplus \,\k h_{2,1}\oplus \,\k f_{2,1}$. We can assume,
after an admissible sign change possibly, that $z_1=f_{3,2},\,$
$z_1^*=e_{3,1},\,$ $z_2=f_{1,1}$, and $z_2^*=\frac{1}{3}e_{1,0}$.
Then $[e_{3,1},f_{3,2}]=\frac{1}{3}[e_{1,0},f_{1,1}]=f_{0,1}=f$.
Put $u_1^*=e_{3,2},\,$ $u_1=-f_{3,1},\,$ $u_2^*=e_{1,1},$ and
$u_2=af_{1,0}$. where $a\in\k^\times$. Clearly, $\z_\chi(1)$ is
spanned by the $u_1, u_2,u_1^*, u_2^*$. Since
$$([e_{3,2},-f_{3,1}],[f_{3,2},e_{3,1}])=
(e_{3,2},[f_{3,2},h_{3,1}])=-(e_{3,2},[h_{3,1},f_{3,2}])=(e_{3,2},f_{3,2})=1,$$
there is $a\in\k^\times$ such that $[u_i,u_j]=[u_i^*,u_j^*]=0$ and
$[u_i^*,u_j]=-\delta_{ij}\,e$ for $1\le i,j\le 2$.

Let $\{E,H,F\}$ be the standard basis of $\sl_2(\k)$ and
${\mathcal E}= \Theta(\z_\chi(0))\cup\Theta(\z_\chi(1))\cup\{C\}$.
Let $\rho\colon {\mathcal E}\rightarrow \Mat_2(\k)$ be such that
$\rho(\Theta_{x e_{2,1}+yh_{2,1}+zf_{2,1}})=xE+yH+zF,\,$
$\rho(\Theta_u)=0,$ and $\rho(C)=-\frac{16}{9}\,I_2$ for all
$x,y,z\in\k$ and $u\in\z_\chi(1)$. We claim that the elements from
$\rho(\mathcal E)$ satisfy the relations (i)\! -\! (iv) of
Theorem~\ref{main3}. Since the relations (i)\,\! -\! (iii) are
satisfied for obvious reasons we just need to check the quadratic
relation (iv).

Since $2\alpha_1+\alpha_2$ is a short root, we have
$(h_{2,1},h_{2,1})=6$ and $(e_{2,1},f_{2,1})=3$. Therefore,
$C_0=\frac{1}{3}(e_{2,1}f_{2,1}+f_{2,1}e_{2,1}+\frac{1}{2}\,h_{2,1}^2)$.
Note that $EF+FE+\frac{1}{2}\,H^2=\frac{3}{2}\,I_2$. Since
$c_0=-\frac{28}{9}$ in the present case, we have to show that
\begin{eqnarray}\label{relation}
\sum_{i=1}^{4}\,\big(\rho(\Theta_{[u,z_i]^\sharp})\,\rho(\Theta_{[v,z_i^*])^\sharp})+
\rho(\Theta_{[v,z_i^*]^\sharp})\,\rho(\Theta_{[u,z_i]^\sharp})\big)\,=\,-\frac{5}{6}(f,[u,v])I_2
\end{eqnarray}
for all $u,v\in\z_\chi(1)$. When $u$ and $v$ run through the set
$\{u_1,u_2,u_1^*,u_2^*\}$, the LHS of (\ref{relation}) is always a
linear combination of matrices  $XY+YX$ with $X,Y\in\sl_2(\k)$.
Since all such matrices are multiples of $I_2$, the LHS of
(\ref{relation}) equals $g(u,v)I_2$ for some skew-symmetric
bilinear form $g$ on $\z_\chi(1)$. Using the relations (i) and
(ii) of Theorem~\ref{main3} it is easy to observe that this form
in $\z_\chi(0)$-invariant. As $\z_\chi(1)$ is an irreducible
$\z_\chi(0)$-module, there is a scalar $c\in\k$ such that
$g(u,v)=c(f,[u,v])$ for all $u,v\in\z_\chi(1)$. Thus we need to
check that $c=-\frac{5}{6}$.

Note that $[u_1^*,z_i^*]=[u_1,z_i]=0$ for $i=1,2$. Also,
$[u_1^*,z_2]^\sharp=[u_1^*,z_2]$ and
$[u_1,z_2^*]^\sharp=[u_1,z_2^*]$. It follows that for $i=1,2$,
$$([u_1^*,z_i],[u_1,z_i^*])=(u_1^*,[u_1,[z_i,z_i^*]])=-(u_1^*,[u_1,f])=
-([u_1^*,u_1],f)=(e,f)=1.$$ As $(u_1^*,z_1)=-(u_1,z_1^*)=1$, we
have $[u_1^*,z_1]^\sharp=[u_1^*,z_1]-\frac{1}{2}h$ and
$[u_1,z_1^*]^\sharp=[u_1,z_1^*]+\frac{1}{2}h$. Consequently,
\begin{eqnarray*}
([u_1^*,z_1]^\sharp,[u_1,z_1^*]^\sharp)&=&([u_1^*,z_1],[u_1,z_1^*])+\frac{1}{2}([u_1^*,z_1],h)
-\frac{1}{2}(h,[u_1,z_1^*])-\frac{1}{4}(h,h)\\
&=&1+\frac{1}{2}(u_1^*,z_1)-\frac{1}{2}(u_1,z_1^*)-\frac{1}{2}=2-\frac{1}{2}=\frac{3}{2}.
\end{eqnarray*}
Since $[u_1^*,z_1]^\sharp,\,$ $[u_1^*,z_2],\,$
$[u_1,z_1^*]^\sharp,$ and $[u_1,z_2^*]$ are multiples of
$h_{2,1},\,$ $e_{2,1},\,$ $h_{2,1}$, and $f_{2,1}$, respectively,
and $(e_{2,1},f_{2,1})=\frac{1}{2}(h_{2,1},h_{2,1})=3$, the
preceding remarks show that
$$
\sum_{i=1}^{4}\,\big(\rho(\Theta_{[u_1^*,z_i]^\sharp})\,\rho(\Theta_{[u_1,z_i^*])^\sharp})+
\rho(\Theta_{[u_1,z_i^*]^\sharp})\,\rho(\Theta_{[u_1^*,z_i]^\sharp})\big)=
\frac{1}{3}(EF+FE)+\frac{1}{2}H^2=\frac{5}{6}I_2.
$$
As $(f,[u_1^*,u_1])=-1$, we deduce from (\ref{relation}) that
$c=-\frac{5}{6}$, as wanted. As Theorem~\ref{main3} gives a
presentation of $H$ by generators and relations, the result
follows.
\end{pf}
\begin{rem}
It is immediate from the proof of Lemma~\ref{L6} that for $\g$ of
type ${\rm G}_2$ the following relation holds in $H$:
\begin{eqnarray}\label{G2}
[\Theta_{u_1^*},\Theta_{u_1}]=\,-\frac{1}{2}\,C+
\frac{1}{3}\big(\Theta_{e_{2,1}}\Theta_{f_{2,1}}+\Theta_{f_{2,1}}\Theta_{e_{2,1}}+
\frac{1}{2}\Theta^2_{h_{2,1}}\big)+\frac{1}{6}\Theta^2_{h_{2,1}}-\frac{14}{9}.
\end{eqnarray}
The expressions for all $[\Theta_u,\Theta_v]$ with
$u,v\in\z_\chi(1)$ can be derived from (\ref{G2}) by using the
action of $\ad\,\Theta(\z_\chi(0))$ on $\Theta(\z_\chi(i)),\,$
$i=1,2$. For example, it can be deduced easily that
$[\Theta_{u_1^*},\Theta_{u_2^*}]$ is a nonzero scalar multiple of
$\Theta_{e_{2,1}}^2$. This implies that the span of all PBW
monomials in
$C,\Theta_{h_{2,1}},\Theta_{e_{2,1}},\Theta_{u_1^*},\Theta_{u_2^*}$
is a subalgebra of $H$. It can be regarded as a Borel subalgebra
of $H$.
\end{rem}
\subsection{}
As yet another application of Theorems~\ref{main2} and \ref{main3}
we are going to classify all irreducible finite dimensional
representations of $H$ in the case where $\g$ is of type ${\rm
C}_n$ or ${\rm G}_2$. Dimension formulae for these representations
will be given. We shall rely on Joseph's theory of Goldie-rank
polynomials. The reader will notice that our method is quite
general and can be applied to any simple Lie algebra $\g$.
However, various problems remain in the general case, especially
for Lie algebras of type ${\rm E}_7$ and ${\rm E}_8$. We hope to
return to this interesting subject in the future.

Given $\nu\in \h^*$ we denote by $I(\nu)$ the annihilator of the
irreducible highest weight module $L(\nu)$ in $U(\g)$.  Recall
from the proof of Theorem~\ref{main3} that
$\lambda_0=-\frac{1}{2}\varpi_n$ for $\g$ of type ${\rm C}_n$ and
$\lambda_0=-\frac{2}{3}\varpi_2$ for $\g$ of type ${\rm G}_2$. Let
$\Phi_0=\{\alpha\in\Phi\,|\,\la \lambda_0,\alpha\ra\in\Z\}$. It is
easy to see that $\Phi_0$ coincides with the set of all {\it
short} roots in $\Phi$. In particular, $\Phi_0$ is a root system
in $\h^*$ but not a closed subsystem of $\Phi$. The set
$\Pi_0=\{\alpha_1,\ldots,\alpha_{n-1},\alpha_{n-1}+\alpha_n\}$ is
the basis of simple roots in $\Phi_0$ contained in $\Phi^+$. This
implies that $\Phi_0$ is of type ${\rm D}_n$ and ${\rm A}_2$ when
$\g$ is of type ${\rm C}_n$ and ${\rm G}_2$, respectively (our
convention here is that ${\rm D}_2\cong{\rm A}_1\times {\rm A}_1$
and ${\rm D}_3\cong {\rm A}_3$). Note that
$\lambda_0=\,-\frac{d-1}{d}\,\varpi_n$ where $$
d=\frac{(\alpha_{n},\alpha_n)}{(\alpha_{n-1},\alpha_{n-1})}\,=\,
\left\{
\begin{array}{ll}
2 & \mbox{when $\g$ is of type ${\rm C}_n$},\\
3 & \mbox{when $\g$ is of type ${\rm G}_2$}.
\end{array}
\right.
$$

It is well-known that the subgroup $W_0:=\{w\in
W\,|\,w(\lambda_0)-\lambda_0\in \Z\Phi\}$ of $W$ is generated by
the reflections $s_\alpha$ with $\alpha\in\Phi_0$, hence
identifies with the Weyl group of $\Phi_0$. We note for further
references that
\begin{eqnarray*}
\la\lambda_0+\mu+\rho,\alpha_{n-1}+\alpha_n\ra&=&\la
\lambda_0+\mu+\rho,\alpha_{n-1}\ra+d\la
\lambda_0+\mu+\rho,\alpha_{n}\ra\\
&=&\la \mu,\alpha_{n-1}\ra+1+d\big(\la
\mu,\alpha_{n}\ra-\frac{d-1}{d}+1\big)\\
&=&\la \mu,\alpha_{n-1}\ra +d\la \mu,\alpha_{n}\ra+2.
\end{eqnarray*}
\begin{theorem}\label{CG}
Let $\Phi_0^+$ denote the set of all short roots in $\Phi^+$.
\begin{enumerate}

\item[1.] If $\g$ is of type ${\rm C}_n$, $n\ge 2$, then to every
$\mu\in P^+$ there corresponds a unique finite dimensional simple
$H$-module $V_H(\mu)$ such that
$$\dim V_H(\mu)\,=\prod_{\alpha\in\Phi_0^+}\frac{\la2\mu+2\rho-\varpi_n,\alpha\ra}
{\la2\rho-\varpi_n,\alpha\ra}.$$ Any finite dimensional simple
$H$-module is isomorphic to one of the modules $V_H(\mu),\,$
$\mu\in P^+$. The central characters of these modules are pairwise
distinct.

\smallskip

\item[2.] If $\g$ is of type ${\rm G}_2$, then to every
$\mu=a\varpi_1+b\varpi_2\in P^+$ there correspond two finite
dimensional simple $H$-modules $V_H^\pm(\mu)$ such
 that
\begin{eqnarray*}
\dim V^+_H(\mu)&=&\frac{(a+1)(a+3b+2)(2a+3b+3)}{6};\\
\dim V^-_H(\mu)&=&\frac{(a+1)(a+3b+3)(2a+3b+4)}{6}.
\end{eqnarray*}
Any finite dimensional simple $H$-module is isomorphic to one of
the modules $V^{\pm}_H(\mu),\,$ $\mu\in P^+$. The central
characters of these modules  are pairwise distinct.
\end{enumerate}
\end{theorem}
\begin{pf}
By Theorem~\ref{main2} the isoclasses of finite dimensional simple
$H$-modules are in $1$-$1$ correspondence with the primitive
ideals $\cal I$ of $U(\g)$ such that ${\mathcal VA}({\cal
I})=\overline{\cal O}_{\rm min}$. By Duflo's theorem, ${\cal
I}=I(\lambda)$ for some $\lambda\in\h^*$. Let
$\Phi_\lambda=\{\alpha\in\Phi\,|\,\la\lambda,\alpha\ra\in\Z\}$ and
let $\Pi_\lambda$ be the basis of simple roots of $\Phi_\lambda$
contained in $\Phi_\lambda\cap\Phi^+$. As explained in [\cite{J3},
p.~41],  the equality ${\mathcal VA}(I(\lambda))=\overline{\cal
O}_{\rm min}$ holds if and only if $\dim {\cal O}_{\rm
min}=|\Phi|-|\Phi_\lambda|$ and $\la\lambda+\rho,\alpha\ra>0$ for
all $\alpha\in\Pi_\lambda$ (the argument in {\it loc.\,cit.}
relies on the fact that ${\cal O}_{\rm min}$ is  not a special
orbit in the sense of Lusztig when $\g$ is not simply laced).
Since $\dim {\cal O}_{\rm min}=|\Phi|-|\Phi_0|$, we have
$|\Phi_0|=|\Phi_\lambda|$. Now $\Phi_\lambda^\vee$ is a closed
symmetric subsystem of the dual root system $\Phi^\vee\subset\h$.
The Borel--de Siebenthal algorithm implies that there is only one
such subsystem in $\Phi^\vee$ of size $|\Phi_0|$, namely
$\Phi_0^\vee$;\, see [\cite{Bou}, Ch. VI, Sect.~4, Exerc.~4]. This
shows that $\Phi_\lambda=\Phi_0$ and $\Pi_\lambda=\Pi_0$.

 Write $\lambda=\lambda_0+\sum_{i=1}^n l_i\varpi_i$ with
 $l_i\in\k$. Since $\alpha_1,\ldots,\alpha_{n-1}\in\Pi_0$, it must be  that
 $l_i\in\Z_+$ for $1\le i\le n-1$, while our earlier remarks show
 that $\la \lambda+\rho,\alpha_{n-1}+\alpha_n\ra=
 l_{n-1}+dl_n+2$
 is a positive integer. Hence $l_n\in\frac{1}{d}\Z$.
 Since $\Phi_\lambda\ne\Phi$ we have $\lambda\not\in P$, giving
 $l_n\not\in-\frac{1}{d}+\Z$. For $\g$ of type ${\rm C}_n$ this says
 $l_n\in\Z$, while for
$\g$ of type ${\rm G}_2$ we get that either $l_n\in\Z$ or
$l_n\in\frac{1}{3}+\Z$.

It is easy to see that $s_{\alpha_n}$ permutes the positive short
roots in $\Phi$. Therefore,
$I(\lambda)=I(s_{\alpha_n}\centerdot\lambda)$; see [\cite{Ja},
(5.16)]. As $\alpha_n=-d\varpi_{n-1}+2\varpi_n$ we have
\begin{eqnarray*}
s_{\alpha_n}\centerdot\lambda&=&s_{\alpha_n}\Big(\sum_{i=1}^{n-1}\,(l_i+1)\varpi_i+
(l_n+\frac{1}{d})\varpi_n\Big)-\rho
\,=\,\lambda-\frac{dl_n+1}{d}\,\alpha_n\\
&=&\lambda_0+\sum_{i=1}^{n-2}\,l_i\varpi_i+(l_{n-1}+dl_n+1)\varpi_{n-1}+\Big(
l_n-\frac{2}{d}(dl_n+1)\Big)\varpi_n\\
&=&\lambda_0+\sum_{i=1}^{n-2}\,l_i\varpi_i+(l_{n-1}+dl_{n}+1)
\varpi_{n-1}-\big(l_n+\frac{2}{d}\big)\varpi_n.
\end{eqnarray*}
Thus replacing $\lambda$ by $s_{\alpha_n}\centerdot\lambda$ if
necessary we may assume further that $l_n\ge 0$.

Suppose $l_n\in\Z$. Then the above discussion shows that
$\lambda-\lambda_0\in P^+$. It follows from Joseph's theory of
Goldie-rank polynomials that
$$\text{rk}\big(U(\g)/I(\lambda_0+\mu)\big)=\,
c\,{\textstyle\prod_{\alpha\in\Phi_0^+}}\,\la\lambda_0+\mu+\rho,\alpha\ra$$
for all $\mu\in P^+$, where $c$ is a constant independent of
$\mu$; see [\cite{J4}, p.~303]. Recall that Theorem~\ref{main2}
associates to each ${\cal I}\in{\mathcal X}$ with ${\mathcal
VA}({\cal I})=\overline{\cal O}_{\rm min}$ an irreducible finite
dimensional $H$-module (up to isomorphism). Abusing notation we
shall denote this module by ${\varkappa}^{-1}({\cal I})$. We have
already mentioned that ${\mathcal
VA}(I(\lambda_0+\mu))=\overline{\cal O}_{\rm min}$. Therefore, to
each $\mu\in P^+$ there corresponds an irreducible finite
dimensional $H$-module
$V_H(\mu):=\varkappa^{-1}(I(\lambda_0+\mu))$. By
Theorem~\ref{main2}(2),
$$\dim
V_H(\mu)=\,\text{rk}\big(U(\g)/I(\lambda_0+\mu)\big)\qquad(\forall\,\mu\in
P^+).$$ Since $I(\lambda_0)$ is the Joseph ideal,
Theorem~\ref{main2} together with Proposition~\ref{P3} gives $\dim
V_H(0)=\dim H/H^+=1.$ Therefore,
$c^{-1}=\prod_{\alpha\in\Phi_0^+\,}\la\lambda_0+\rho,\alpha\ra$
and \begin{eqnarray}\label{weyl}\dim
V_H(\mu)=\prod_{\alpha\in\Phi_0^+}\frac{\la\lambda_0+\mu+\rho,\alpha\ra}{\la\lambda_0+\rho,\alpha\ra}.
\end{eqnarray}
Since $\lambda_0+\rho+P^+$ is contained in the interior of the
dominant Weyl chamber, the modules in the set
$\{V_H(\mu)\,|\,\mu\in P_+\}$ have pairwise distinct central
characters. This settles the case where $\g$ is of type ${\rm
C}_n$.

Suppose $\g$ of type ${\rm G}_2$. For $\mu=a\varpi_1+b\varpi_2\in
P^+$ we put $V^+_H(\mu):=V_H(\mu)$. Since
$\Phi_0^+=\{\alpha_2,\alpha_1+\alpha_2,2\alpha_1+\alpha_2\}$ and
$\lambda_0=-\frac{2}{3}\varpi_2$, the dimension formula
(\ref{weyl}) reads
$$\dim V^+_H(\mu)=\,\frac{(a+1)(a+3b+2)(2a+3b+3)}{6}.$$

Now suppose $l_n\not\in\Z$. Our earlier remarks show that $\g$ is
of type ${\rm G}_2$ and $l_n\in\frac{1}{3}+\Z_+$. As a
consequence, $\lambda\in\frac{1}{2}\lambda_0+P^+$. For any $\mu\in
P^+$ we have $\Phi_{\frac{1}{2}\lambda_0+\mu}=\Phi_0$. As
$\frac{1}{2}\lambda_0+\mu+\rho$ lies in the interior of the
dominant Weyl chamber, the above argument applies yielding
${\mathcal VA}(I(\frac{1}{2}\lambda_0+\mu))=\overline{\cal O}_{\rm
min}$. Theorem~\ref{main2} shows that
$V^-_H(\mu):=\varkappa^{-1}(I(\frac{1}{2}\lambda_0+\mu))$ is an
irreducible finite dimensional $H$-module with
$$\dim
V_H^-(\mu)=\,\text{rk}\big(U(\g)/I({\textstyle\frac{1}{2}}\lambda_0+\mu)\big).$$
In conjunction with the discussion in [\cite{J4}, p.~303] this
entails that
\begin{eqnarray}\label{weyl'}
\dim
V_H^-(\mu)=\,\text{rk}\big(U(\g)/I({\textstyle\frac{1}{2}}\lambda_0+\mu)\big)=\,
c'\,{\textstyle\prod_{\alpha\in\Phi_0^+}}\,
\la{\textstyle\frac{1}{2}}\lambda_0+\mu+\rho,\alpha\ra
\end{eqnarray}
for all $\mu\in P^+$ where $c'$ is a constant independent of
$\mu$.

It also follows from Theorem~\ref{main2} that $C$ acts on
$V_H^\pm(\mu)$ as $f^{\pm}(\mu)\,\text{id}$ where
$f^+(\mu)=(\lambda_0+\mu,\lambda_0+\mu+2\rho)$ and
$f^-(\mu)=(\frac{1}{2}\lambda_0+\mu,\frac{1}{2}\lambda_0+\mu+2\rho)$.
We have noted in the proof of Theorem~\ref{main3} that
$(\,\cdot\,,\,\cdot\,)=\frac{1}{3}(\,\cdot\,|\,\cdot\,)$. By
[\cite{Bou}, Table~IX], $(\varpi_1|\varpi_1)=2,\,$
$(\varpi_1|\varpi_2)=3,$ and $(\varpi_2|\varpi_2)=6$. Using this
fact it is straightforward to see that
$f^-(\mu)=f^-(a\varpi_1+b\varpi_2)$ is a quadratic polynomial in
$a,\,b$ with all coefficients positive except for the constant
term
$f^{-}(0)=\,(-\frac{1}{3}\varpi_2,-\frac{1}{3}\varpi_2+2\rho)=-\frac{16}{9}.$
Furthermore,
$$f^+(a\varpi_1+b\varpi_2)=\,
\textstyle{\frac{2}{3}a^2+2b^2+2ab+2a+\frac{10}{3}b-{\textstyle\frac{28}{9}}}.$$
As a consequence, $f^+(\mu)>-\frac{16}{9}$ for all nonzero $\mu\in
P^+$.

Let $M$ be an $H$-module affording the representation $\rho$ from
Lemma~\ref{L6}. The above discussion shows that $M\cong
V^\pm(\nu)$ for some $\nu\in P^+$. Since $C$ acts on $M$ as
$-\frac{16}{9}\,\text{id}$, the preceding remark yields $M\cong
V_H^-(0)$. As $\dim M=2$, this allows us to determine the scale
factor $c'$. In view of (\ref{weyl'}) we then get
\begin{eqnarray*}
\dim V_H^-(\mu)&=& 2\prod_{\alpha\in\Phi_0^+}\frac{
\la-{\textstyle\frac{1}{3}}\varpi_2+\mu+\rho,\alpha\ra}{
\la{-\textstyle\frac{1}{3}}\varpi_2+\rho,\alpha\ra}\\
&=&\frac{(a+1)(a+3b+3)(2a+3b+4)}{6}
\end{eqnarray*}
for all $\mu\in P^+$. Since $(-\frac{2}{3}\varpi_2+\rho+P^+)\cap
(-\frac{1}{3}\varpi_2+\rho+P^+)=\emptyset$ and the union
$(-\frac{2}{3}\varpi_2+\rho+P^+)\cup
(-\frac{1}{3}\varpi_2+\rho+P^+)$ is contained in the interior of
the dominant Weyl chamber, the modules in the set
$\{V_H^\pm(\mu)\,|\,\mu\in P^+\}$ have pairwise distinct central
characters. This completes the proof.
\end{pf}
\begin{rem}
When $\g$ is of type ${\rm C}_2$ we have
$\Phi_0^+=\{\varepsilon_1-\varepsilon_2,\,\varepsilon_1+\varepsilon_2\},$
$\varpi_1=\varepsilon_1$, $\varpi_2=\varepsilon_1+\varepsilon_2$,
and $\lambda_0=-\frac{1}{2}(\varepsilon_1+\varepsilon_2)$; see
[\cite{Bou}, Table~III]. In this case our dimension formula reads
$$\dim V_H(\mu)=\,\frac{(a+1)(a+2b+2)}{2}=\,\frac{(r-s+1)(r+s+2)}{2},$$
where
$\mu=a\varpi_1+b\varpi_2=(a+b)\varepsilon_1+b\varepsilon_2=r\varepsilon_1+s\varepsilon_2$
and $r,s\in\Z_+$, $r\ge s$. The same dimension formula can be
found in [\cite{Kh}] where it was obtained by a completely
different method in the context of deformed oscillator algebras of
rank one; see Remark~\ref{Osc}.
\end{rem}
\section{\bf Highest weight modules for $H$}
\subsection{}

Let $\Phi_e$ denote the set of all $\alpha\in\Phi$ with
$\alpha(h)\in\{0,1\}$, and put $\Phi_e^\pm=\Phi_e\cap \Phi^\pm$,
$\Phi_{e,i}^{\pm}=\{\alpha\in\Phi_e^\pm\,|\,\alpha(h)=i\}$. Recall
that $\z_\chi$ is spanned by $\h_e$, by all $e_\alpha$ with
$\alpha\in\Phi_e$, and by $e$.  Let $h_1,\ldots, h_{l-1}$ be a
basis of $\h_e$, and let $\n^\pm(i)$ be the span of all $e_\alpha$
with $\alpha\in\Phi^\pm_{e,i}$. Clearly, $\n^+(0)$ and $\n^-(0)$
are maximal nilpotent subalgebras of $\g(0)^\sharp$. Let
$\{x_1,\dots ,x_t\}$ and $\{y_1,\dots,y_t\}$ be bases of $\n^+(0)$
and $\n^-(0)$ consisting of root vectors $e_\alpha$ with
$\alpha\in\Phi$. Recall that the $z_i$'s with $1\le i\le 2s$ are
root vectors for $\h$. For $1\le i\le s$, set $u_i=[e,z_i]$ and
$u_i^*=[e,z_i^*]$. It follows from our discussion in (4.1) that
$u_i$ (resp. $u_i^*$) is a root vector for $\h$
 corresponding to the root $\beta+\gamma_i\in\Phi_{e,1}^-$ (resp.
$\beta+\gamma_i^*\in\Phi_{e,1}^+$). Furthermore,  $\{u_1,\ldots,
u_s,u_1^*,\ldots,u_s^*\}$ is a $\k$-basis of $\z_\chi(1)$.

Given a linear function $\lambda$ on $\h_e$ and $c\in\k$ we denote
by $J_{\lambda,c}$ the linear span in $H$ of all PBW monomials of
the form
$$\prod_{i=1}^t\Theta_{y_i}^{l_i}\cdot
\prod_{i=1}^s\Theta_{u_i}^{m_i} \cdot\prod_{i=1}^{\ell-1}
\big(\Theta_{h_i}-\lambda(h_i)\big)^{n_i}\cdot (C-c)^{n_\ell}\cdot
\prod_{i=1}^s\Theta_{u_i^*}^{r_i}\cdot
\prod_{i=1}^t\Theta_{x_i}^{q_i},$$ where
$\sum_{i=1}^{\ell}n_i+\sum_{i=1}^t r_i+\sum_{i=1}^s q_i>0$.
\begin{lemma}\label{L7}
The subspace $J_{\lambda,c}$ is a left ideal of the algebra $H$.
\end{lemma}
\begin{pf}
For ${\bf a},{\bf b}\in\Z_+^t$, ${\bf c},{\bf d}\in\Z_+^s$, ${\bf
m}\in\Z_+^{\ell}$, set
$$\Theta({\bf a,b,c,d,m}):=\,
\textstyle{(\prod_{i=1}^t\Theta_{y_i}^{a_i})\,(\prod_{i=1}^s\Theta_{u_i}^{c_i})\,
(\prod_{i=1}^{\ell-1} \Theta_{h_i}^{m_i})\,
C^{m_l}\,(\prod_{i=1}^s\Theta_{u_i^*}^{d_i})\,
(\prod_{i=1}^t\Theta_{x_i}^{b_i}}).$$ By [\cite{P02},
Theorem~4.6(ii)], the PBW monomials $\Theta({\bf a,b,c,d,m})$ form
a $\k$-basis of $H$. Note that $\deg_e\big(\Theta({\bf
a,b,c,d,m})\big)=4m_\ell+3(|{\bf c}|+|{\bf d}|)+2({|\bf a}|+|{\bf
b}|)+2\sum_{i=1}^{\ell-1}m_i$.

Since $C-c$ is central in $H$ we have $\Theta({\bf a,b,c,d,m})
(C-c)\in J_{\lambda,c}$. Relations~(i) and (ii) of
Theorem~\ref{main3} imply that $\Theta({\bf a,b,c,d,m})
(\Theta_{h_i}-\lambda(h_i))\in J_{\lambda,c}$ for $1\le i\le
\ell-1$.  Since $\Theta(\n^+(0))$ is a Lie subalgebra of
$\Theta(\z_\chi(0))$, by Theorem~\ref{main3}, we also have that
$\Theta({\bf a,b,c,d,m})\cdot \Theta_{e_\alpha}\in J_{\lambda,c}$
for all $\alpha\in\Phi_{e,0}^+$.

It remains to show that $\Theta({\bf a,b,c,d,m})\cdot
\Theta_{u_i^*}\in J_{\lambda,c}$ for all $i\le s$. We shall use
induction on $\deg_e\big(\Theta({\bf a,b,c,d,m}) \big)$, so assume
from now that $\deg_e\big(\Theta({\bf a,b,c,d,m}) \big)=N$ and
$H^k\cdot \Theta_{u_i^*}\in J_{\lambda,c}$ for all $i\le s$ and
all $k<N$. First note that the span of $u_1^*,\ldots, u_s^*$
equals $\n^+(1)$, hence is stable under the adjoint action of
$\n^+(0)$. Since we have already established that $H\cdot
\Theta_{e_\alpha}\in J_{\lambda,c}$ for all $\alpha\in\Phi_{e,0}$,
relation~(ii) of Theorem~\ref{main3} yields $$\Theta({\bf
a,b,c,d,m})\cdot\Theta_{u_i^*}\in\, \Theta({\bf
a,b,c,0,m})\cdot\Theta(\n^+(1))+J_{\lambda,c}.$$ Thus we may
assume that ${\bf d}={\bf 0}$. If $b_j=0$ for all $j>i$, then
$\Theta({\bf a,b,c,0,m})\cdot\Theta_{u_i^*}=\Theta({\bf a,b+e_{\it
i},c,0,m})\in J_{\lambda,c}$. So suppose ${\bf b}=(b_1,\ldots,
b_k,0,\ldots,0)$ where $b_k>0$ and $k>i$. Then in view of
[\cite{P02}, Theorem~4.6(iv)] and our induction assumption we have
\begin{eqnarray*}
\Theta({\bf a,b,c,0,m})\cdot\Theta_{u_i^*}&\in&\Theta({\bf
a,b+e_{\it i},c,0,m})+\Theta({\bf a,b-e_{\rm
k},c,0,m})[\Theta_{u_k^*},\Theta_{u_i^*}]\\ &+& H^{N-2}\cdot
\Theta_{u_k^*}\,\subseteq\, \Theta({\bf a,b-e_{\rm
k},c,0,m})[\Theta_{u_k^*},\Theta_{u_i^*}]+J_{\lambda,c}.
\end{eqnarray*}
Since $(f,[u_k^*,u_i^*])=0$,  Theorem~\ref{main3} shows that $$
[\Theta_{u_k^*},\Theta_{u_i^*}]\,=\,
\frac{1}{2}\sum_{j=1}^{2s}\,\big(\Theta_{[u_{k}^*,z_j]^\sharp}\,\Theta_{[u_i^*,z_j^*]^\sharp}+
\Theta_{[u_i^*,z_j^*]^\sharp}\,\Theta_{[u_k^*,z_j]^\sharp}\big)\in
\textstyle{\sum}_{\alpha\in\Phi^+_{e,0}}\,H\cdot
\Theta_{e_\alpha}$$ (one should take into account that
$[u_i^*,z_j^*],\,\big[[u_i^*,z_j^*]^\sharp,[u_k^*,z_j]^\sharp\big]\in\bigcup_{{\alpha\in\Phi^+_{e,0}}}
\,\k e_\alpha$ for all $j\le s$). So $\Theta({\bf
a,b,c,0,m})\cdot\Theta_{u_i^*}\in J_{\lambda,c}$, and the result
follows by indiction on $N$.
\end{pf}

\subsection{}
Put $Z_H(\lambda,c):=H/J_{\lambda,c}$ and let $v_0$ denote the
image of $1$ in $Z_H(\lambda,c)$. Clearly, $Z_H(\lambda,c)$ is a
cyclic $H$-module generated by $v_0$. We call $Z_H(\lambda,c)$ the
{\it Verma module of level $c$ corresponding to} $\lambda$. By
Lemma~\ref{L7}, the vectors
$$\{\Theta_{y_1}^{l_1}\cdots\Theta_{y_t}^{l_t}\Theta_{u_1}^{m_1}\cdots\Theta_{u_s}^{m_s}(v_0)\,|
\,\,l_1,\ldots, l_t,m_1,\ldots, m_s\in\Z_+\}$$ form a $\k$-basis
of the Verma module $Z_H(\lambda,c)$. Let $Z_H^+(\lambda,c)$
denote the $\k$ span of all
$\{\Theta_{y_1}^{l_1}\cdots\Theta_{y_t}^{l_t}\Theta_{u_1}^{m_1}\cdots\Theta_{u_s}^{m_s}(v_0)$
with $\sum_i l_i+\sum_i m_i>0$. Let $Z_H^{\rm max}(\lambda,c)$
denote the sum of all $H$-submodules of $Z_H(\lambda,c)$ contained
in $Z_H^+(\lambda,c)$, and put
$$L_H(\lambda,c):=Z_H(\lambda,c)/Z_H^{\rm max}(\lambda,c).$$
\begin{prop}\label{P4} The following are true:
\begin{itemize}
\item[(i)] $Z_H^{\mathrm max}(\lambda,c)$ is a unique maximal
submodule of the Verma module $Z_H(\lambda,c)$ and hence
$L_H(\lambda,c)$ is a simple $H$-module.

\smallskip

\item[(ii)] The simple $H$-modules $L_H(\lambda,c)$ and
$L_H(\lambda',c')$ are isomorphic if and only if
$\lambda=\lambda'$ and $c=c'$.

\smallskip

\item[(iii)] Any finite dimensional simple $H$-module is
isomorphic to one of the modules $L_H(\lambda,c)$ with
$\lambda\in\h_e^*$ satisfying $\lambda(h_\alpha)\in\Z_+$ for all
$\alpha\in\Phi_{e,0}^+$. Furthermore, if $\g$ is not of type $\rm
A$ then $c$ is a rational number.
\end{itemize}
\end{prop}
\begin{pf}
(a) For a root $\alpha=\sum_{i=1}^\ell n_i\alpha_i$ in $\Phi$ we
put
$${\rm ht}_\beta(\alpha):=\sum_{\alpha_i\ne\beta}\,n_i.$$
Clearly, ${\rm ht}_\beta(\alpha)=0$ if and only if
$\alpha=\pm\beta$. As all derivations of $\g$ are inner, there is
a unique $h_0\in\h$ such that $[h_0,e_\alpha]={\rm
ht}_\beta(\alpha)e_\alpha$ for all $\alpha\in\Phi$. As
$[h_0,e_{\pm\beta}]=0$ we have that $h_0\in\h_e$.  Clearly,
$\Theta_{h_0}(v_0)=\lambda(h_0)v_0$ and $Z_H(\lambda,c)\,=\,\k
v_0\oplus Z_H^+(\lambda,c)$. Since all $y_i$ and $z_i$ are root
vectors for $\h$ corresponding to negative roots different from
$-\beta$, it follows from Theorem~\ref{main3} that the subspace
$Z_H^+(\lambda,c)$ decomposes into eigenspaces for $\Theta_{h_0}$
and the eigenvalues of $\Theta_{h_0}$ on $Z_H^+(\lambda,c)$ are of
the form $\lambda(h_0)-k$ where $k$ is a positive integer.

Let $V$ be a nonzero $H$-submodule of $Z_H(\lambda,c)$. If
$V\not\subseteq Z_H^+(\lambda,c)$, the above discussion shows that
$v_0\in V$. But then $V=Z_H(\lambda,c)$. Thus any proper submodule
of $Z_H(\lambda,c)$ is contained in $Z_H^+(\lambda,c)$. As a
consequence, $Z_H^{\rm max}(\lambda,c)$ is a unique maximal
submodule of $Z_H(\lambda,c)$, proving (i).

\smallskip

\noindent (b) It follows from part~(a) that each $H$-module
$L_H(\lambda,c)$ decomposes into eigenspaces for $\Theta_{h_0}$,
the eigenvalues of $\Theta_{h_0}$ on $L_H(\lambda,c)$ lie in the
set $\lambda(h_0)-\Z_+$, and the eigenspace
$L_H(\lambda,c)_{\lambda(h_0)}$ is spanned $v_0$. If
$L_H(\lambda,c)\cong L_H(\lambda',c')$ as $H$-modules then it must
be that $\lambda(h_0)\in\lambda'(h_0)-\Z_+$ and
$\lambda'(h_0)\in\lambda(h_0)-\Z_+$. This implies that
$\lambda(h_0)=\lambda'(h_0)$ and
$L_H(\lambda,c)_{\lambda(h_0)}\cong
L_H(\lambda',c')_{\lambda'(h_0)}$ as modules over the commutative
subalgebra $\Theta(\h_e)\oplus\k C$ of $H$. But then
$\lambda=\lambda'$ and $c=c'$, hence (ii).

\smallskip

\noindent (c) Let $M$ be a finite dimensional simple $H$-module.
Then $C\in Z(H)$ acts on $M$ as $c\,\text{id}$ for some $c\in \k$.
Since $\Theta(\h_e)$ is abelian, $M$ contains at least one weight
subspace for $\Theta(\h_e)$. From Theorem~\ref{main3} it follows
that the direct sum $\bigoplus_{\mu\in\h_e^*} M_\mu$ of all weight
subspaces of $M$ is an $H$-submodule of $M$.  Hence $M$ decomposes
into weight spaces relative to $\Theta(\h_e)$. Note that any
linear function on $\h$ vanishing on $\h_e$ is a scalar multiple
of $\beta$. Since $\beta$ is a simple root, any sum of roots from
$\Phi_e^+$ restricts to a nonzero function on $\h_e$. But then the
relation $$\phi\ge \psi\iff
\phi=\psi+(\textstyle{\sum}_{\gamma\in\Phi_e^+}\,
r_{\gamma}\gamma)_{\vert {\h_e}},\quad
r_\gamma\in\Z_+,\qquad\,(\forall\,\phi,\psi\in\h_e)$$ is a partial
ordering on $\h_e^*$. Since the set of $\Theta(\h_e)$-weights of
$M$ is finite, it contains at least one maximal element with
respect to this ordering, $\lambda$ say. Let $m$ be a nonzero
vector in $M_\lambda$. Then $\Theta_{x_i}(m)=\Theta_{u_i^*}(m)=0$
for all admissible $i$.  As a consequence, there exists a
homomorphism of $H$-modules $\xi\colon\,Z_H(\lambda,c)\rightarrow
M$ such that $\xi(v_0)=m$. The simplicity of $M$ implies that
$\xi$ is surjective, while part~(a) yields ${\rm
Ker}\,\xi=Z_H^{\rm max}(\lambda,c)$. Restricting $M$ to the
${\frak sl}_2$-triple
$(\Theta_{e_\alpha},\Theta_{h_\alpha},\Theta_{e_{-\alpha}})\subset
H$ with $\alpha\in\Phi_{e,0}$ it is easy to observe that
$\lambda(h_\alpha)\in\Z_+$ for all $\alpha\in\Phi_{e,0}^+$.

Finally, suppose $\g$ is not of type $\rm A$. Then $\z_\chi(0)$ is
a semisimple Lie algebra. By Weyl's theorem, $M$ is a completely
reducible $\Theta(\z_\chi(0))$-module. Let $\g_{\mathbb Q}$ be the
$\mathbb Q$-form in $\g$ spanned by the Chevalley system from
(4.1), and $\z_{\chi,\mathbb Q}(i)=\g_{\mathbb Q}\cap\z_\chi(i)$
where $i=0,1$.  Choose $u,v\in\z_{\chi,\mathbb Q}(1)$ with
$(f,[u,v])=2$. Then
$[u,z_i]^\sharp,[v,z_i^*]^\sharp\in\z_{\chi,\mathbb Q}(0)$ for all
$i$. The highest weight theory implies that there is a $\mathbb
Q$-form in $M$ stable under the action of $\Theta(\z_{\chi,\mathbb
Q}(0))$. It follows that ${\rm
tr}_M\big(\Theta_{[u,z_i]^\sharp}\Theta_{[v,z_i^*]^\sharp}\big)\in{\mathbb
Q}$ for $1\le i\le 2s$. Since ${\rm tr}_M
[\Theta_{u},\Theta_v]=0$, Theorem~\ref{main3} entails that
$(c-c_0)\dim M\in\mathbb Q.$ Since $c_0\in\mathbb Q$ by
Theorem~\ref{main3}, we obtain $c\in\mathbb Q$.
\end{pf}
\subsection{}
To determine the composition factors of the Verma modules
$Z_H(\lambda,c)$ with their multiplicities we are going to
establish a link between these $H$-modules and the $\g$-modules
obtained by parabolic induction from Whittaker modules for ${\frak
sl}(2,\k)$. The latter modules have been studied in much detail in
[\cite{Mc}, \cite{MS}, {\cite{Bac}], and it is known that their
composition multiplicities can be calculated by using the
Kazhdan-Lusztig algorithm. We are going to rely on Skryabin's
equivalence (3.1); the Kazhdan filtration of $H$ will play an
important r{\^o}le  too.

Let ${\frak s}_\beta$ denote the subalgebra of $\g$ spanned
$(e,h,f)= (e_\beta,h_\beta,f_\beta)$, and put
$$\p_\beta:=\,{\frak s}_\beta+\h+\textstyle{\sum}_{\alpha\in\Phi^+}\,\,
\k e_\alpha,\quad\,
\n_\beta:=\,\textstyle{\sum}_{\alpha\in\Phi^+\setminus\{\beta\}}\,\,\k
e_\alpha,\quad\, \widetilde{\frak s}_\beta:=\,\h_e\oplus
{\mathfrak s}_\beta.$$ Clearly, $\p_\beta=\widetilde{\mathfrak
s}_\beta\oplus\n_\beta$ is a parabolic subalgebra of $\g$ with
nilradical $\n_\beta$ and $\widetilde{\mathfrak s}_\beta$ is a
Levi subalgebra of $\p_\beta$. Let
$C_\beta=ef+fe+\frac{1}{2}h^2=2ef+\frac{1}{2}h^2-h$ be the Casimir
element of $U({\mathfrak s}_\beta)$. Given $\lambda\in\h_e^*$ and
$c\in\k$ we denote by $I_\beta(\lambda,c)$ the left ideal of
$U(\p_\beta)$ generated by $f-1,C_\beta-c$, all $h-\lambda(h)$
with $h\in\h_e$, and all $e_\gamma$ with
$\gamma\in\Phi^+\setminus\{\beta\}$.

Define $Y(\lambda,c):=\,U(\p_\beta)/I_\beta(\lambda,c)$, a
$\p_\beta$-module with the trivial action of $\n_\beta$, and let
$1_{\lambda,c}$ denote the image of $1$ in $Y(\lambda,c)$. Since
$f . 1_{\lambda,c}=1_{\lambda,c}$ we have that
$$e . 1_{\lambda,c}=\,\textstyle{\frac{1}{2}}(C_\beta-
\textstyle{\frac{1}{2}}h^2+h) .
1_{\lambda,c}=(\textstyle{-\frac{1}{4}h^2+\frac{1}{2}h+\frac{1}{2}c})
. 1_{\lambda,c}.$$ Together with the PBW theorem this shows that
the vectors $\{h^k\cdot 1_{\lambda,c}\,|\,\,k\in\Z_+\}$ form a
$\k$-basis of $Y(\lambda,c)$ (the independence of these vectors
follows from the fact that $Y(\lambda,c)$ is infinite
dimensional). We mention for completeness that $Y(\lambda,c)$ is
isomorphic to a Whittaker module for ${\mathfrak
s}_\beta\cong\sl(2,\k)$.

The above discussion shows that the vectors
$$m({\bf i, j, k},l):=\,z_1^{i_1}\cdots
z_s^{i_s}\cdot y_1^{j_1}\cdots y_t^{j_t}\cdot u_1^{k_1}\cdots
u_s^{k_s}\cdot h^l(1_{\lambda,c})$$ with ${\bf i,k}\in \Z_+^s,\,
{\bf j}\in \Z_+^t$, and $l\in \Z_+$ form a $\k$-basis of the
induced $\g$-module $$M(\lambda,c):=\,U(\g)\otimes_{U(\p_\beta)}
Y(\lambda,c).$$
\subsection{} Recall from (4.1) that each $z_i^*$ with $1\le i\le s$ is
a root vector for $\h$ corresponding to
$\gamma_i^*=-\beta-\gamma_i\in\Phi^+$. Put
$\delta=\frac{1}{2}(\gamma_1^*+\cdots+\gamma_s^*)$ and
$\rho_0=\rho-2\delta-(s+1)\beta=\sum_{\alpha\in\Phi_{e,0}^+}\,\alpha$.
Since the restriction of $(\,\cdot\,,\,\cdot\,)$ to $\h_e$ is
nondegenerate, for any $\eta\in\h_e^*$ there exists a unique
$t_\eta\in \h_e$ such that $\varphi=(t_\eta,\,\cdot\,)$. Hence
$(\,\cdot\,,\,\cdot\,)$ induces a nondegenerate bilinear form on
$\h_e^*$ via $(\mu,\nu):=(t_{\mu},t_{\nu})$ for all
$\mu,\nu\in\h_e^*$. Given a linear function $\varphi$ on $\h$ we
denote by $\bar{\varphi}$ the restriction of $\varphi$ to $\h_e$.
\begin{theorem}\label{main4}
Each $\g$-module $M(\lambda,c)$ is an object of the category
${\mathcal C}$. Furthermore, ${\mathrm Wh}(M(\lambda,c))\cong
Z_H(\lambda+\bar{\delta},c+(\lambda+2\bar{\rho},\lambda))$ as
$H$-modules.
\end{theorem}
\begin{pf}
Put $M:=M(\lambda,c)$, and let $M_0$ (resp. $M_1$) denote the
$\k$-span of all $m({\bf i, j, k},l)\in M$ with $|{\bf i}|+l=0$
(resp. $|{\bf i}|+l>0$). Clearly, $M=M_0\oplus M_1$ as vector
spaces. Let ${\rm pr}\colon M\,=\,M_0\oplus M_1\twoheadrightarrow
M_0$ denote the first projection.

If $|{\bf i}|+2|{\bf j}|+3|{\bf k}|+2l=k$, we say that $m({\bf i,
j, k},l)$ has {\it Kazhdan degree} $k$. Let $M^k$ denote the
$\k$-span in $M$ of all $m({\bf i, j, k},l)$ of Kazhdan degree
$\le k$. Then $\{M^k\,|\,\,k\in\Z_+\}$ is an increasing filtration
in $M$ and $M^0=\k 1_{\lambda,c}$. Taking $U(\g)$ with its Kazhdan
filtration (as defined in [\cite{GG}] for example) we can thus
regard $M$ as a filtered $U(\g)$-module.

Let $z=\lambda f+\sum_{i=1}^s \mu_i z_i^*\in\m_\chi$ where
$\lambda,\mu_i\in\k$. Since $z_i^*\in\n_\beta$  for $1\le i\le s$
and $f . 1_{\lambda,c}=1_{\lambda,c}$, we have that $z.
1_{\lambda,c}=\lambda\cdot 1_{\lambda,c}=\chi(z)\cdot
1_{\lambda,c}$. Since $z$ acts locally nilpotently on $U(\g)$, we
deduce that $z-\chi(z)$ acts locally nilpotently on $M$ for all
$z\in\m_\chi$. As a consequence, $M$ is an object of ${\mathcal
C}_\chi$. By our discussion in (3.1),  ${\mathrm Wh}(M)\ne 0$, the
algebra $H$ acts on $M$, and $M\cong Q_\chi\otimes_{H}\,{\rm
Wh}(M)$ as $\g$-modules.

Now observe that  $$z_k^*.\,m({\bf i, j, k},l)\,\in \,i_k\cdot
m({\bf i-e_{\it k}, j, k},l)+\text{span}\,\{m({\bf i', j',
k'},l')\,|\,\,|{\bf i}'|\ge |{\bf i}|\}$$ for all $k\le s$, and
$$(f-1).\,m({\bf i, j, k},l)\,\in \,2^l\cdot
m({\bf i, j, k},0)+\text{span}\,\{m({\bf i', j',
k'},l')\,|\,\,l'>0\}\quad\mbox{when }\, l>0.$$ From this it is
immediate that the map $\text{pr}\colon {\rm Wh}(M)\rightarrow
M_0$ is injective.

Note that $1_{\lambda,c}\in\text{Wh}(M)$, and for all $h\in\h_e$
we have
\begin{eqnarray*}
\Theta_h(1_{\lambda,c})&=&\Big(h+\frac{1}{2}\sum_{i=1}^{2s}\,
[h,z_i^*]z_i\Big)(1_{\lambda,c})\,=\,\Big(h+\frac{1}{2}\sum_{i=1}^s\,[[h,z_i^*],z_i]\Big)(1_{\lambda,c})\\
&=&\Big(\lambda(h)+\frac{1}{2}\sum_{i=1}^s\gamma_i^*(h)f\Big)\cdot
1_{\chi,c}\,=\,(\lambda+\delta)(h)\cdot 1_{\chi,c}.
\end{eqnarray*}
Suppose $x\in\z_\chi (0)$ is a root vector for $\h$ corresponding
to root $\gamma\in\Phi_{e,0}^+$. Then $x\in\n_\beta$ and
$[[x,z_i^*],z_i]\in\n_\beta$ for all $i\le s$. Therefore,
$$\Theta_x(1_{\lambda,c})\,=\,\Big(x+\frac{1}{2}\sum_{i=1}^{2s}\,
[x,z_i^*]z_i\Big)(1_{\lambda,c})\,=\,\Big(h+\frac{1}{2}\sum_{i=1}^s\,[[x,z_i^*],z_i]\Big)(1_{\lambda,c})=0.$$

Recall from (2.5) that for any positive root vector
$u\in\z_\chi(1)$ we have
\begin{eqnarray*}
z_u\,=\,-\frac{1}{3} \sum_{i=1}^{2s}\big\la
z_i^*,[u,[z_{i},z_i^*]]\big\ra z_i
=\frac{1}{3}\sum_{i=1}^{2s}\big\la z_i^*,[u,f]\big\ra
z_i\in\n_\beta.
\end{eqnarray*}
This implies that
\begin{eqnarray*}
\Theta_{u}(1_{\lambda,c})&=&\Big(u+\sum_{i=1}^{2s}\,[u,z_i^*]\,z_i+\frac{1}{3}
\sum_{i,j=1}^{2s}\,[uz_i^*z_j^*]\,z_jz_i+z_{u}
\Big)(1_{\lambda,c})\\
&=&\Big(\sum_{i=1}^{s}\,[[u,z_i^*],z_i]+\frac{1}{3}
\sum_{i,j=1}^{s}\,[uz_i^*z_j^*]\,z_jz_i +
\frac{1}{3}\sum_{i=1}^{s}\sum_{j=s+1}^{2s}\,[uz_i^*z_j^*]\,z_jz_i
\Big)(1_{\lambda,c})\\
&=&\frac{1}{3}\Big(\sum_{i,j=1}^{s}\,[[uz_i^*z_j^*],z_j]\,z_i +
\sum_{i,j=1}^{s}\,z_j\,[[uz_i^*z_j^*],z_i]-
\sum_{i=1}^{s}\,[uz_i^*z_i]\,[z_i^*,z_i] \Big)(1_{\lambda,c})\\
&=&\frac{1}{3}\Big(\sum_{i,j=1}^{s}\,[[[uz_i^*z_j^*],z_j],z_i] -
\sum_{i=1}^{s}\,[uz_i^*z_i]f \Big)(1_{\lambda,c})\,\in\,\n_\beta .
1_{\lambda,c}=0.
\end{eqnarray*}
Therefore, $\Theta_{u_i^*}(1_{\lambda,c})=0$ for all $i\le s$. Our
discussion in (4.4) shows that
\begin{eqnarray}\label{Cas}
C(1_{\lambda,c})&=&\Big(2e+\frac{h^2}{2}-(s+1)h+C_0+2\sum_{i=1}^{2s}\,[e,z_i^*]z_i\Big)
(1_{\lambda,c})\nonumber\\
&=&\Big(C_\beta-sh+C_0+2\sum_{i=1}^{s}\,[[e,z_i^*],z_i]\Big)
(1_{\lambda,c}).
\end{eqnarray}
As $[[[e,z_i^*],z_i],f]=[[[e,f],z_i^*],z_i]=[[h,z_i^*],z_i]=-f$ we
have $[[e,z_i^*],z_i]-\frac{1}{2}h\in\h_e$ for all $i\le s$. Let
$x$ be an arbitrary element in $\h_e$. Then $(x,h)=0,\,$
$\beta(x)=0$, and
$$\big(x,[[e,z_i^*],z_i]-\textstyle{\frac{1}{2}}h\big)\,=\,
\big([x,[e,z_i^*]],z_i]\big)\,=\,\gamma_i^*(x)\big([e,z_i^*],z_i\big)=\gamma_i^*(x),$$
that is $[[e,z_i^*],z_i]-\frac{1}{2}h=t_{\bar{\gamma}_i^*}$ for
all $i\le s$; see our discussion at the beginning of this
subsection. But then
\begin{eqnarray}\label{Cas'}
\Big(2\sum_{i=1}^{s}\,[[e,z_i^*],z_i]-sh\Big)(1_{\lambda,c})\,=\,
\Big(2\sum_{i=1}^s\,t_{\bar{\gamma}_i^*}\Big)(1_{\lambda,c})\,=\,4(\lambda,\bar{\delta})\cdot
1_{\lambda,c}.
\end{eqnarray}
Since $C_0=\sum a_ib_i$ is a Casimir element of $U(\z_\chi(0))$
and all positive root vectors in $\z_\chi(0)$ annihilate
$1_{\lambda,c}$, it is straightforward to see that
$C_0(1_{\lambda,c})=(\lambda,\lambda+2\bar{\rho}_0)\cdot
1_{\lambda,c}.$ In conjunction with (\ref{Cas}) and (\ref{Cas'})
this yields
$$C(1_{\lambda,c})\,=\,\big(c+(\lambda,\lambda+2\bar{\rho}_0)+4(\lambda,\bar{\delta})\big)\cdot
1_{\lambda,c}\,=\,(c+(\lambda,\lambda+2\bar{\rho}))\cdot
1_{\lambda,c}.$$ Put $\lambda':=\lambda+\bar{\delta}$ and
$c':=c+(\lambda,\lambda+2\bar{\rho}).$  Let $V_0$ denote the
$H$-submodule of $M$ generated by $1_{\lambda,c}$. The above
discussion shows that the left ideal $J_{\lambda',c'}$ of $H$
annihilates $1_{\lambda,c}$. Therefore, $V_0$ is a homomorphic
image of the Verma module $Z_H(\lambda',c')$.

We claim that the restriction of $\text{pr}\colon
M\twoheadrightarrow M_0$ to $V_0$ is surjective. Recall that $M_0$
is spanned by all $m({\bf 0, j, k},0$) with ${\bf j}\in \Z_+^t$
and ${\bf k}\in\Z_+^s$. Clearly, $m({\bf
0,0,0},0)=1_{\lambda,c}\in \text{pr}(V_0)$. Assume that all
vectors $m({\bf 0, j, k},0)$ of Kazhdan degree $2|{\bf j}|+3|{\bf
k}|<n$ are in $\text{pr}(V_0)$. Now let $m({\bf 0, a, b},0)\in
M_0$ be such that $2|{\bf a}|+3|{\bf b}|=n$ and $|{\bf a}|+|{\bf
b}|=k$, and denote by $M_{n,k}$  the span of all $m({\bf i, j,
k},l)$ of Kazhdan degree $n$ with $|{\bf i}|+|{\bf j}|+|{\bf
k}|+l>k$. Assume that all vectors $m({\bf 0, j, k},0)$ of Kazhdan
degree $n$ with $|{\bf j}|+|{\bf k}|>k$ are in $\text{pr}(V_0)$.
Since $M$ is a filtered $U(\g)$-module, it follows from
Lemmas~\ref{L3} and~\ref{L5} that
\begin{eqnarray}\label{ind}
\Theta_{y_1}^{a_1}\cdots\Theta_{y_t}^{a_t}
\Theta_{u_1}^{b_1}\cdots\Theta_{u_s}^{b_s}(1_{\lambda,c})\,\in\,
m({\bf 0, a, b},0)+M_{n,k}+M^{n-1}.
\end{eqnarray}
In view of our assumptions on $n$ and $k$ we get $m({\bf 0, a,
b},0)\in\text{pr}(V_0+M^{n-1}+M_{n,k})=\text{pr}(V_0)$. Our claim
now follows by double induction on $n$ and $k$. Since we have
already established that $\text{pr}\colon \text{Wh}(M)\rightarrow
M_0$ is injective, this yields $\text{Wh}(M)=V_0$.

Using (\ref{ind}) it is easy to observe that the vectors $
\Theta_{y_1}^{a_1}\cdots\Theta_{y_t}^{a_t}
\Theta_{u_1}^{b_1}\cdots\Theta_{u_s}^{b_s}(1_{\lambda,c})$ with
${\bf a}\in\Z_+^t$ and ${\bf b}\in \Z_+^s$ are linearly
independent over $\k$. Hence it follows from Lemma~\ref{L7} and
our discussion at the beginning of (7.2) that $V_0\cong
Z_H(\lambda',c')$ as $H$-modules.
\end{pf}
\begin{rem}\label{remlast}
Combined with Skryabin's equivalence and the main results of
Mili{\v c}i{\'c}--Soergel [\cite{McS}] and Backelin [\cite{Bac}],
Theorem~1.3 implies that the composition multiplicities of the
Verma modules $Z_H(\lambda,c)$ can be computed with the help of
inverse parabolic Kazhdan-Lusztig polynomials associated with the
coset $W/\langle s_\beta\rangle$. This confirms in the minimal
nilpotent case the Kazhdan-Lusztig conjecture for finite $\mathcal
W$-algebras as formulated by de~Vos and van~Driel in
[\cite{DvvD}]. Recall that our construction of $H_\chi$ is a
special instance of quantum Hamiltonian reduction where the
constraints imposed are read off the $\sl_2$-triple
$(e_\beta,h_\beta,e_{-\beta})$. In the physics literature the
algebra $H$ appears undercover under the name of a finite
$\mathcal W$-algebra associated with the minimal embedding ${\frak
sl}(2,\k)\hookrightarrow\g$.
\end{rem}

\begin{rem}
It would be interesting to relate the (Kazhdan) filtered algebra
$H$ with the Becchi--Rouet--Stora--Tyutin (BRST) quantisation of
the Poisson algebra $\text{gr}\,H$. We recall that the Poisson
structure on $\text{gr}\,H$ is determined in [\cite{GG}]. It would
be important (for the characteristic $p$ theory and possibly for
the theory of minimal representations of $p$-adic groups) to
determine all $(\lambda,c)\in\h_e^*\times\k$ such that the simple
$H$-module $L_H(\lambda,c)$ is finite dimensional; see
Proposition~\ref{P4}(iii).
\end{rem}
Given $t\in \k$ we let $H_t$ denote the factor algebra $H/(C-t)$
where $(C-t)$ is the two-sided ideal of $H$ generated by the
central element $C-t$. It is clear from the definition that each
$L(\lambda,t)$ is an $H_t$-module.
\begin{corollary}
If $\g$ is of type ${\mathrm C}_n$ or ${\mathrm G}_2$, then any
finite dimensional $H_t$-module is completely reducible.
\end{corollary}
\begin{pf}
It suffices to show that $\text{Ext}^1_{H_t}(M,N)=0$ for any two
finite dimensional simple $H_t$-modules $M$ and $N$. If the
$H_t$-modules $M$ and $N$ are not isomorphic, then they have
distinct central characters; see Theorem~\ref{CG}. Thus it remains
to show that $\text{Ext}^1_{H_t}(M,M)=0$. By
Proposition~\label{P4}(iii), $M\cong L_H(\lambda,t)$ for some
$\lambda\in \h_e^*$. Let $V$ be a finite dimensional $H_t$-module
containing $M$ as a submodule and such that $V/M\cong M$ as
$H_t$-modules. As $\Theta(\z_\chi(0))$ is a semisimple Lie
subalgebra of $H$, Weyl's theorem yields that $V$ decomposes into
weight spaces relative to $\Theta(\h_e)$, say
$V=\,\bigoplus_{\mu\in X(V)\,} V_\mu$. Moreover, the set of
$\Theta(\h_e)$-weights of $V$ coincide with that of
$L(\lambda,c)$, showing that $\mu\le\lambda$ for all $\mu\in
X(V)$. It follows from our discussion in (7.2) that $\dim
M_{\lambda}=1$ and $\dim V_\lambda=2$. Let $v\in
V_\lambda\setminus M_\lambda$ and let $M'$ denote the
$H$-submodule of $V$ generated by $v$. By construction, the left
ideal $J_{\lambda,t}$ of $H$ annihilates $v$, showing that $M'$ is
a homomorphic image of the Verma module $Z_H(\lambda,t)$. But then
the $\Theta(\h_e)$-weight space $M'_\lambda$ is $1$-dimensional,
implying $(M\cap M')_\lambda=M_\lambda\cap M'_\lambda=\,0$.
Consequently, $M\cap M'=0$. The irreducibility of $V/M$ now
entails that $V=M\oplus M'$ and $M'\cong M$. Then
$\text{Ext}^1_{H_t}(M,M)=0$, completing the proof.
\end{pf}

\end{document}